\documentclass[final,3p,11pt, times]{article}
\usepackage{authblk}
\usepackage{amsmath}
\usepackage{amssymb,latexsym}
\usepackage[german,english]{babel}
\usepackage{url}
\usepackage{graphicx}
\usepackage{gastex}
\usepackage{longtable}
\usepackage{lscape}
\usepackage{tabularx}
\usepackage{multicol}
\usepackage{verbatim}
\usepackage{multirow}
\usepackage{graphicx}
\usepackage{epsfig,graphics,graphicx}
\usepackage{geometry}
\usepackage{authblk}
\usepackage{verbatim}
\usepackage{multirow}
\usepackage{ragged2e}

\usepackage[T1]{fontenc}
\usepackage{authblk}

\usepackage{lineno}

\usepackage{lipsum}

\hbadness=\maxdimen

\newtheorem{lem}{{\bf Lemma}}[section]
\newtheorem{thm}{{\bf Theorem}}[section]
\newtheorem{prop}{{\bf Proposition}}[section]

\begin{document}
	
	\title{{2-semi equivelar maps on the torus and the Klein bottle with few vertices}}
	
	\author[1]{Anand Kumar Tiwari}
	\author[1]{Yogendra Singh}
	\author[2]{Amit Tripathi}
	
	\affil[1]{\small Department of Applied Science, Indian Institute of Information Technology, Allahabad 211\, 015, India, E-mail: anand@iiita.ac.in, tiwari.iitp@gmail.com}

	\affil[2]{Department of Applied Science \& Humanities, Rajkiya Engineering College, Banda -210201, India,\\ E-mail: amittripathi@recbanda.ac.in}

	\maketitle
	
	\hrule
	
	\begin{abstract}
		The $k$-semi equivelar maps, for $k \geq 2$, are generalizations of maps on the surfaces of Johnson solids to closed surfaces other than the 2-sphere. In the present study, we determine 2-semi equivelar maps of curvature 0 exhaustively on the torus and the Klein bottle. Furthermore, we classify (up to isomorphism) all these 2-semi equivelar maps on the surfaces with up to 12 vertices. 
	\end{abstract}

	\smallskip
	
	\textbf{MSC(2010):} Primary: 52B70; Secondary: 52C20.
	
	\textbf{Keywords:} 2-Semi equivelar map, Torus, Klein bottle.
	
	\smallskip
	
	\hrule
	
\section{Introduction} \label{s1}


A map $M$ on a surface $S$ is an embedding of a graph $G$ into $S$ such that (i) the closure of each component of $S \setminus G$ is topologically a $p$-gonal 2-disk $D_p$ ($p \geq 3$), which is called a face of $M$ and (ii) the non-empty intersection of any two distinct faces is either a vertex or an edge, see \cite{brehm(1990)}. The vertices and edges of the underlying graph $G$ in a map $M$ are called the vertices and edges of the map. Let $M_1$ and $M_2$ be maps with the vertex sets $V(M_1)$ and $V(M_2)$ respectively on a surface $S$. Then $M_1$ is isomorphic to $M_2$, denoted as $M_1 \cong M_2$, if there is a bijective map $f\,: V(M_1) \to V(M_2)$ such that $e$ is an edge in $M_1$ iff $f(e)$ is an edge in $M_2$, $D_p$ is a $p$-gonal face in $M_1$ iff $f(D_p)$ is a $p$-gonal face in $M_2$, and $f$ preserves the incidence of edges and faces. 




The face-sequence of a vertex $v \in V(M)$, denoted as $f_{seq}(v)$, is $f_{seq}(v)=(p_1^{n_1}. p_2^{n_2} \ldots p_k^{n_k})$ if the consecutive $n_1$ numbers of $p_1$ gons, $n_2$ numbers of $p_2$ gons, $\ldots$ , $n_k$ numbers of $p_k$ gons are incident at $v$ in the given cyclic order. The curvature of the vertex $v$, denoted as $\phi(v)$, is then defined as $\phi(v) = 1 - (\sum_{i=1}^{k} n_i)/2 + \sum_{i=1}^{k} (n_i/p_i)$. We say that a map $M$ has the combinatorial curvature $k$ if $\phi(v) = k$ for every $v \in V(M)$. Some maps of positive curvature are discussed in \cite{zhang, rbk}.

A map $M$ having $k$ distinct face-sequences, say $f_1, f_2, \ldots, f_k$, is called a $k$-semi equivelar map of type $[f_1:f_2: \cdots: f_k]$. The 11 Archimedean tilings, 20 2-uniform tilings, 61 3-uniform tilings, 151 4-uniform tilings, 332 5-uniform tilings, and 673 6-uniform tilings on the Euclidean plane provide 1-, 2-, 3-, 4-, 5-, and 6-semi equivelar maps on the plane respectively, see \cite{ch(1989),gs(1987)}. 

In the case when $k=1$, the map is referred to as a 1-semi equivelar map (or simply semi equivelar map). Datta and Maity \cite{dm(2022)} described all types of semi equivelar maps on the surfaces of Euler characteristic 2 (2-sphere $\mathbb{S}^2$) and 1 (projective plane $\mathbb{RP}^2$). The first author with Maity and Upadhyay classified some semi equivelar maps on the surface of Euler characteristic $-1$, \cite{tu(2016),utm(2012)}. Karab\'{a}\v{s} and Nedela  \cite{karabas(2007),karabas(2012)} have described some semi equivelar maps on the orientable surfaces of Euler characteristic $-2$ (double torus), $-4$, and $-6$.  Semi equivelar maps have been studied extensively for the surfaces of Euler characteristic 0, that is, on the torus and Klein bottle.
Kurth \cite{kur(1984)} as well as  Brehm and K\"{u}hnel \cite{brehm(2008)} have given a technique to enumerate semi equivelar maps of type $[3^6]$, $[4^4]$ and $[6^3]$ on the torus. Datta and Nilakantan \cite{datta(2001)} have classified these types of maps on the torus and Klein bottle for $\leq 11$ vertices. Further, Datta and upadhyay \cite{daak(2005)} have extended this classification for $\leq 15$ vertices. In \cite{tu(2017)}, the first author with upadhyay have classified 1-semi equivelar maps of types $[3^4.6]$, $[3^3.4^2]$, $[3^2.4.3.4]$, $[3.4.6.4]$, $[3.6.3.6]$, $[3.12^2]$, $[4.6.12]$, $[4.8^2]$ on the torus and Klein bottle on at most 20 vertices. Recently, Dipendu \cite{dipendu(2021)} has studied some 2-semi equivelar maps on the torus that are obtained by taking quotient of 2-uniform tilings of the Euclidean plane.






A process of subdividing a $p$-gonal face by introducing a new vertex and joining it to each vertex of the polygon by an edge is called stacking of the face. Evidently, stacking of each face of semi-equivelar maps of types $[3^6]$ and $[4^4]$ on the torus and Klein bottle lead to 2-semi equivelar maps of types $[3^3:3^{12}]$ and $[3^4: 3^8]$ on the surfaces respectively. In this article, we show: 




	\begin{thm}\label{t1}
		Let $T$ denote the type of a 2-semi equivelar map of curvature 0 on the torus or Klein bottle. Then $T \in \{[3^6:3^4.6]$, $[3^6:3^3.4^2]$, $[3^6:3^2.4.3.4]$, $[3^6:3^2.4.12]$, $[3^6:3^2.6^2]$, $[3^4.6:3.6.3.6]$, $[3^3.4^2:3^2.4.3.4]$, $[3^3.4^2:3.4.6.4]$, $[3^3.4^2:4^4]$, $[3^2.4.3.4:3.4.6.4]$, $[3^2.6^2:3^4.6]$, $[3^2.6^2:3.6.3.6]$, $[3.4^2.6:3.4.6.4]$, $[3.4^2.6:3.6.3.6]$, $[3.4.6.4:4.6.12]$,  $[3.12^2:3.4.3.12]\}$. 
	\end{thm}
	
	
Next, we classify and enumerate all these types of 2-semi equivelar maps on $\leq 12$ vertices and show:  
	
	\begin{thm}\label{t2}
		There are exactly 31 2-semi equivelar maps of curvature 0 on the surfaces of Euler characteristic 0 on at most $12$ vertices; 18 of which are on the torus and 13 are on the Klein bottle. These 18 are  $\mathcal{A}_2(T)$, $\mathcal{A}_3(T)$, 
		$\mathcal{B}_2(T)$, $\mathcal{C}_2(T)$, $\mathcal{E}_3(T), \mathcal{E}_4(T), \mathcal{E}_6(T), \mathcal{E}_8(T)$, $\mathcal{E}_{11}(T)$, $\mathcal{E}_{13}(T)$, $\mathcal{E}_{14}(T)$, $\mathcal{F}_2(T)$, $\mathcal{F}_3(T)$, $\mathcal{F}_5(T)$, $\mathcal{F}_6(T)$, $\mathcal{F}_7(T)$, $\mathcal{F}_9(T)$, and 13 are $\mathcal{A}_1(K)$, $\mathcal{B}_1(K)$
		$\mathcal{C}_1(K)$, $\mathcal{D}_1(K)$, $\mathcal{E}_1(K)$, $\mathcal{E}_2(K)$, $\mathcal{E}_5(K)$, $\mathcal{E}_7(K)$, $\mathcal{E}_{9}(K)$, $\mathcal{E}_{10}(K)$, $\mathcal{F}_1(K)$, $\mathcal{F}_4(K)$, $\mathcal{F}_8(K)$, given in example Section \ref{s3}.
	\end{thm}

This article is organized in the following manner. In Section \ref{s2}, we give some definitions and notations that are used in the further sections. In Section \ref{s3}, we present examples of 2 semi-equivelar maps on the torus and the Klein bottle. Further, in Section \ref{s4}, we describe proof of all the results. Next, we end up this article by presenting some concluding remarks in Section \ref{s5}.

\section{Definitions and notations} \label{s2}

Let $M$ be a map with the vertex set $V(M)$, edge set $E(M)$, and face set $F(M)$. For $v \in V(M)$, consider $K_v = \{f \in F(M) \,: v \in f\}$. Then the geometric carrier $|K_v|$ (the union of all the elements in $K_v$) is a 2-disk with the boundary cycle $C_n=C_n (v_1, v_2, \ldots, v_n)$. The cycle $C_n$ is called the link of the vertex $v$ and is denoted as ${\rm lk}(v)$. For example, if $v$ is a vertex with $f_{seq}(v)=(3^6)$, $(3^3.4^2)$ or $(3^2.4.3.4)$, then ${\rm lk}(v)=C_6(v_1,v_2,v_3,v_4,v_5,v_6)$, ${\rm lk}(v) = C_7 (v_1, \boldsymbol{v_2},v_3,v_4,v_5,v_6, \boldsymbol{v_7})$ or ${\rm lk}(v) = C_7(v_1,v_2,v_3, \boldsymbol{v_4}, v_5,v_6,\boldsymbol{v_7})$, respectively (see in the following figures).  Here, the bold appearance of some $v_i$'s means $v$ is not adjacent with these $v_i$'s. We use a similar notation in Section \ref{s4} to express the link of a vertex with a specific face-sequence.


\vspace{-1.2cm}

\begin{picture}(0,0)(-17,33)
\setlength{\unitlength}{10mm}

\drawpolygon(0,0)(1,0)(1.5,1)(1,2)(0,2)(-.5,1)

\drawline[AHnb=0](-.5,1)(1.5,1)
\drawline[AHnb=0](0,0)(1,2)
\drawline[AHnb=0](1,0)(0,2)

\put(.7,.8){\scriptsize $v$}

\put(-.1,-.2){\scriptsize $v_1$}
\put(.9,-.2){\scriptsize $v_2$}
\put(1.6,1){\scriptsize $v_3$}
\put(.9,2.1){\scriptsize $v_4$}
\put(-.1,2.1){\scriptsize $v_5$}
\put(-.8,1){\scriptsize $v_6$}

\put(-1,-.7){\scriptsize ${\rm lk}(v) = C_6(v_1,v_2,v_3,v_4,v_5,v_6)$}

\end{picture}

\begin{picture}(0,0)(-60,28)
\setlength{\unitlength}{10mm}

\drawpolygon(0,0)(2,0)(2,1)(1.5,2)(.5,2)(0,1)

\drawline[AHnb=0](0,1)(2,1)
\drawline[AHnb=0](1,0)(1,1)
\drawline[AHnb=0](1,1)(1.5,2)
\drawline[AHnb=0](1,1)(.5,2)

\put(.7,.8){\scriptsize $v$}

\put(-.1,-.2){\scriptsize $v_7$}
\put(.9,-.2){\scriptsize $v_1$}
\put(2,-.2){\scriptsize $v_2$}
\put(2.1,1){\scriptsize $v_3$}
\put(1.4,2.1){\scriptsize $v_4$}
\put(.3,2.1){\scriptsize $v_5$}
\put(-.3,1){\scriptsize $v_6$}

\put(-.6,-.7){\scriptsize ${\rm lk}(v) = C_7(v_1, \boldsymbol{v_2}, v_3,v_4,v_5,v_6, \boldsymbol{v_7})$}

\end{picture}

\begin{picture}(0,0)(-125,12)
\setlength{\unitlength}{10mm}

\drawpolygon(.5,-1)(1.5,-.5)(1,.5)(0,1)(-1,.5)(-1.5,-.5)(-.5,-1)

\drawline[AHnb=0](0,0)(0,1)
\drawline[AHnb=0](0,0)(1,.5)
\drawline[AHnb=0](0,0)(-1,.5)
\drawline[AHnb=0](0,0)(.5,-1)
\drawline[AHnb=0](0,0)(-.5,-1)

\put(.1,-.1){\scriptsize $v$}

\put(.3,-1.2){\scriptsize $v_1$}
\put(1.6,-.6){\scriptsize $v_2$}
\put(1,.6){\scriptsize $v_3$}
\put(-.1,1.1){\scriptsize $v_4$}
\put(-1.2,.6){\scriptsize $v_5$}
\put(-1.7,-.4){\scriptsize $v_6$}

\put(-.6,-1.2){\scriptsize $v_7$}

\put(-1.5,-1.7){\scriptsize ${\rm lk}(v) = C_7(v_1, \boldsymbol{v_2}, v_3,v_4,v_5, \boldsymbol{v_6}, v_7)$}

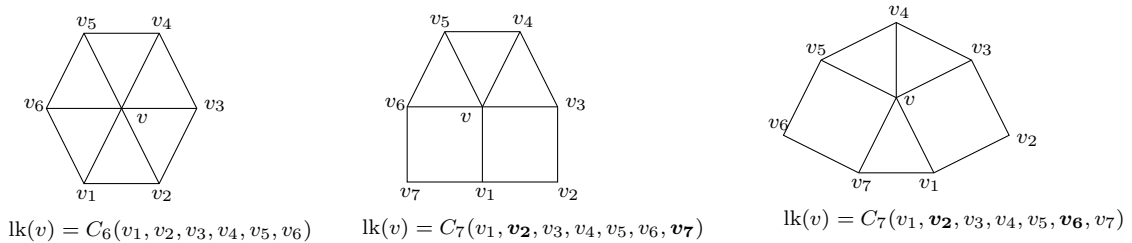
\put(-10,-2.3){\scriptsize Figure 2.1: vertex $v$ with face-sequences $(3^6)$, $(3^3.4^2)$ and $(3^2.4.3.4)$ (see from left)}

\end{picture}

	\newpage

\section{Examples: 2-semi-equivelar maps on the torus and Klein bottle}\label{s3}

Here, we present examples of 2-semi equivelar maps on the surfaces of Euler characteristic 0. The notations $\mathcal{A}_i(T)$-$[f_1:f_2]$ to $\mathcal{O}_i(T)$-$[g_1:g_2]$ represent maps on the torus of type $[f_1:f_2]$ to $[g_1:g_2]$ respectively. Similarly, the notations $\mathcal{A}_i(K)$-$[f_1:f_2]$ to $\mathcal{O}_i(K)$-$[g_1:g_2]$ represent maps on the Klein bottle of type $[f_1:f_2]$ to $[g_1:g_2]$ respectively.


	
\vspace{-.4cm}



\vspace{7.5cm}
	

Now, we show the following:


\begin{lem} \label{l0}
For the maps given in Examples \ref{s3}, we have:
\begin{enumerate}
\item[a:] $\mathcal{E}_3(T) \ncong \mathcal{E}_6(T)$. 
\item[b:] $\mathcal{E}_i(T) \ncong \mathcal{E}_j(T)$, for $i,j \in \{4,8,11,12,13\}$ and $ i \neq j$.
\item[c:] $\mathcal{F}_2(T) \ncong \mathcal{F}_5(T)$. 
\item[d:] $\mathcal{F}_i(T) \ncong \mathcal{F}_j(T)$, for $i,j \in \{3,6,7,9\}$ and $ i \neq j$. 
\item[e:] $\mathcal{E}_1(K) \ncong \mathcal{E}_5(K)$. 
\item[f:] $\mathcal{E}_i(K) \ncong \mathcal{E}_j(K)$, for $i,j \in \{2,7,9,10\}$ and $ i \neq j$. 
\item[i:] $\mathcal{F}_1(K) \ncong \mathcal{F}_8(K)$. 
\end{enumerate}
\end{lem}

\noindent{\bf Proof.} Note that in $\mathcal{E}_3(T)$, we see exactly one vertical non-contractible cycle of length 3 at each vertex (for example, we see $C_3(0,1,4)$ at vertex 4), while in $\mathcal{E}_6(T)$, there are two non-contractible cycles at each vertex (for example, we have $C_3(0,1,4)$ and $C_3(2,3,4)$ at vertex 4). This proves $a$. Following the same argument, we see $\mathcal{E}_1(K) \ncong \mathcal{E}_5(K)$. This proves part $e$. 

\smallskip

In $\mathcal{F}_5(K)$, at each vertex, we have a vertical non-contractible cycle of length 3 (for example, we see $C_5(0,1,5)$ at vertex 1), which is not true in $\mathcal{F}_2(K)$. This proves $c$. Following the same argument, we get $\mathcal{F}_1(K) \ncong \mathcal{F}_8(K)$. This proves $i$.

\smallskip

The following polynomials $p_{G(M)}(a)$ denote the characteristic polynomial (computed from MATLAB) of the adjacency matrix associated with the underlying graph $G$ of map $M$. We know that if two maps have different characteristic polynomials, then they are non-isomorphic. This proves, $b$, $d$, and $f$. 


\smallskip

\noindent $p_{G(\mathcal{E}_4(T))}(a) = a^{12} - 32a^{10} -40a^9 + 254a^8 + 440a^7 - 628a^6 - 1400a^5 + 105a^4 + 1000a^3 + 300a^2$.

\smallskip

\noindent $p_{G(\mathcal{E}_8(T))}(a) = a^{12} - 32a^{10} -48a^9 + 254a^8 + 656a^7 - 292a^6 - 2352a^5 - 2167a^4 + 624a^3 + 2044a^2 + 1120a+192$.

\smallskip

\noindent $p_{G(\mathcal{E}_{11}(T))}(a) = a^{12} - 31a^{10} -32a^9 + 222a^8 + 180a^7 - 746a^6 - 220a^5 + 1201a^4 - 228a^3 - 647a^2 + 322 a$.

\smallskip

\noindent $p_{G(\mathcal{E}_{12}(T))}(a) = a^{12} - 33a^{10} -44a^9 + 258a^8 + 432a^7 - 682a^6 - 1032a^5 + 957a^4 + 560a^3 - 789a^2 + 276 a-32$.

\smallskip

\noindent $p_{G(\mathcal{E}_{13}(T))}(a) = a^{12} - 33a^{10} -44a^9 + 252a^8 + 456a^7 - 568a^6 - 1296a^5 + 348a^4 + 1328a^3 + 108a^2 - 432 a-128$.

\smallskip

\noindent $p_{G(\mathcal{F}_{3}(T))}(a) = a^{12} - 26a^{10} -17a^9 + 176a^8 + 91a^7 - 505a^6 - 95a^5 + 590a^4 - 90a^3 - 118a^2 +24 a$.

\smallskip

\noindent $p_{G(\mathcal{F}_{6}(T))}(a) = a^{12} - 28a^{10} -24a^9 + 212a^8 + 280a^7 - 524a^6 - 976a^5 + 80a^4 + 860a^3 + 528a^2 +96 a$.

\smallskip

\noindent $p_{G(\mathcal{F}_{7}(T))}(a) = a^{12} - 27a^{10} -20a^9 + 201a^8 + 192a^7 - 532a^6 - 552a^5 + 492a^4 + 560a^3 - 84a^2 -192 a-44$.

\smallskip

\noindent $p_{G(\mathcal{F}_{9}(T))}(a) = a^{12} - 27a^{10} -20a^9 + 207a^8 + 168a^7 - 610a^6 - 288a^5 + 723a^4 - 136a^3 - 171a^2 + 84a-11$.

\smallskip

\noindent $p_{G(\mathcal{E}_{2}(K))}(a) = a^{12} - 32a^{10} -40a^9 + 254a^8 + 440a^7 - 644a^6 - 1400a^5 + 457a^4 + 1640a^3 + 156a^2 - 640a-192$.

\smallskip

\noindent $p_{G(\mathcal{E}_{7}(K))}(a) = a^{12} - 32a^{10} -48a^9 + 258a^8 + 640a^7 - 364a^6 - 220a^5 - 1635a^4 + 496a^3 + 684a^2 - 32a-64$.

\smallskip

\noindent $p_{G(\mathcal{E}_{9}(K))}(a) = a^{12} - 31a^{10} -39a^9 + 227a^8 + 377a^7 - 561a^6 - 1129a^5 + 416a^4 + 1283a^3 + 92a^2 - 492a-144$.

This proves the lemma. \hfill $\Box$

	
	\section{Proofs: Classification of 2-semi equivelar maps} \label{s4}
	
	In this section, we describe proof of the Theorems \ref{t1} and \ref{t2}. In our earlier study \cite{st(2022)}, we have shown the following: 
	
	\begin{prop}\label{p1}
		Let $v$ be a vertex with the face-sequence $f$ such that $\phi(v)=0$. Then $f \in S$, where $S= \{(3^3.4^2)$, $(3^6)$, $(3.4^2.6)$, $(3^2.6^2)$, $(3^4.6)$, $(3^2.4.3.4)$, $(3.6.3.6)$, $(4^4)$, $(3.4.6.4)$, $(3^2.4.12)$, $(4.8^2)$, $(3.12^2)$, $(6^3)$, $(5^2.10)$, $(3.8.24)$, $(3.9.18)$, $(3.10.15)$, $(4.5.20)$, $(3.7.42)$, $(4.6.12)$, $(3.4.3.12) \}$.
	\end{prop}

From the above proposition and the compatibility of face-sequences, we see that if a 2-semi equivelar map exists, then its possible type $T \in S_1= A \cup B$, where $A =  \{[3^3.4^2:3^2.6^2]$, $[3^3.4^2:3.4^2.6]$, $[3^3.4^2:3^4.6]$, $[3^3.4^2:3^2.4.12]$, $[3^3.4^2: 4.8^2]$, $[3^3.4^2:4.5.20]$, $[3^3.4^2:4.6.12]$, $[3^3.4^2:3.4.3.12]$, $[3.4^2.6:3^2.6^2]$, $[3.4^2.6:3^4.6]$, $[3.4^2.6:3^2.4.3.4]$, $[3.4^2.6:4^4]$, $[3.4^2.6:3^2.4.12]$, $[3.4^2.6:4.8^2]$, $[3.4^2.6:6^3]$, $[3.4^2.6:4.5.20]$, $3.4^2.6: 4.6.12]$, $[3.4^2.6:3.4.3.12]$, $[3^2.6^2:3^2.4.3.4]$, $[3^2.6^2:3^2.4.12]$, $[3^2.6^2:6^3]$, $[3^4.6:3^2.4.3.4]$, $[3^4.6: 3^2.4.12]$, $[3^4.6:6^3]$, $[3^4.6:4.6.12]$, $[3^2.4.3.4:4^4]$, $[3^2.4.3.4:3^2.4.12]$, $[3^2.4.3.4:3.4.3.12]$,  $[3^2.4.3.4:4.8^2]$, $[3.6.3.6:4.6.12]$, $[3.6.3.6:6^3]$, $[4^4:3.4.6.4]$, $[4^4:3^2.4.12]$, $[4^4:4.8^2]$, $[4^4:4.5.20]$, $[4^4:4.6.12]$, $[4^4:3.4.3.12]$, $[3.4.6.4:3^2.4.12]$, $[3.4.6.4:6^3]$, $[3.4.6.4:4.8^2]$, $[3.4.6.4:4.5.20]$, $[3.4.6.4:3.4.3.12]$, $[3^2.4.12:3.12^2]$, $[3^2.4.12:3.4.3.12]$, $[3^2.4.12:4.5.20]$, $[3^2.4.12:4.6.12]$, $[3^2.4.12:4.8^2]$,  $[4.8^2:4.5.20]$, $[4.8^2:4.6.12]$, $[4.8^2:3.4.3.12]$, $[3.12^2:3.4.3.12]$, $[6^3:4.6.12]$, $[5^2.10:4.5.20] \}$ and 

$B = \{[3^6:3^4.6]$, $[3^6:3^3.4^2]$, $[3^6:3^2.4.3.4]$, $[3^6:3^2.4.12]$, $[3^6:3^2.6^2]$, $[3^4.6:3.6.3.6]$, $[3^3.4^2:3^2.4.3.4]$, $[3^3.4^2:3.4.6.4]$, $[3^3.4^2:4^4]$, $[3^2.4.3.4:3.4.6.4]$, $[3^2.6^2:3^4.6]$, $[3^2.6^2:3.6.3.6]$, $[3.4^2.6:3.4.6.4]$, $[3.4^2.6:3.6.3.6]$, $[3.4.6.4:4.6.12]$,  $[3.12^2:3.4.3.12]\}$. 

\smallskip

We note the following:

\smallskip

\noindent{\bf Remark.} Let $M$ be a 2-semi equivelar map of type $[f_1:f_2]$ with the vertex set $V(M)$. Then for $f_1$ (or $f_2$), there is a vertex $v$ in $V(M)$ with the face-sequence $f_1$ (resp. $f_2$) such that ${\rm lk}(v)$ contains a vertex $u$ with the face-sequence $f_2$ (resp. $f_1$). Such vertex $v$ is called a critical vertex. Clearly, $M$ does not exist if it has no critical vertex.

Now we prove the following:

\smallskip

\noindent{\bf Proof of Theorem \ref{t1}.} As given above, a 2-semi equivelar map $M$ has possible type $T \in S_1$. Note that for each $T \in B$, there exists a 2-semi equivelar map of the type $T$, see in Section \ref{s3}. For $T \in A$, there exist no 2-semi-equivelar maps of type $T$. We establish this in the form of following claim:

\smallskip

\noindent{\bf Claim.} There exist no 2-semi equivelar maps of type $T \in A$.

Consider a map $M$ of type $T \in A$, where $T = [f_1=3^3.4^2: f_2]$. Let $x_0$ be a critical vertex with the face-sequence $(3^3.4^2)$ such that ${\rm lk}(x_0) = C_7(x_1,\boldsymbol{x_2},x_3,x_4,x_5,x_6,\boldsymbol{x_7})$. Then we have the following cases for $f_2$:

\smallskip 

\noindent{\bf Case 1.} $f_2=(3^2.6^2)$, i.e., $T=[f_1=3^3.4^2: f_2=3^2.6^2]$. Then at least one vertex in $\{x_4, x_5\}$ has the face-sequence $(3^2.6^2)$. If $f_{seq}(x_4) = (3^2.6^2)$ (or $f_{seq}(x_5)=(3^2.6^2)$), then $x_3$ (resp. $x_6$) has the face-sequence $f_3 \neq f_1, f_2$. This means $x_0$ is not a critical vertex. 

\smallskip

\noindent{\bf Case 2.} $f_2 = (3.4^2.6)$, i.e., $T= [f_1=3^3.4^2: f_2=3.4^2.6]$. Then at least one vertex in $\{x_1, x_2, x_3, x_6, x_7\}$ has the face-sequence $(3.4^2.6)$. If $f_{seq}(x_3) = (3.4^2.6)$, then we see consecutive 2 triangular faces and one hexagonal face at $x_4$, which implies $f_{seq}(x_4) \neq f_1, f_2$. Similarly, if $f_{seq}(x_6)=(3.4^2.6)$, then $f_{seq}(x_5) \neq f_1, f_2$.
If $f_{seq}(x_1)=(3.4^2.6)$, then ${\rm lk}(x_1) = C_9(x_0,\boldsymbol{x_3}, x_2,\boldsymbol{x_8},\boldsymbol{x_9},\boldsymbol{x_{10}},x_{11}$, $x_7, \boldsymbol{x_{6}})$
or ${\rm lk}(x_1) = C_9(x_0,\boldsymbol{x_6}, x_7,\boldsymbol{x_{11}},\boldsymbol{x_{10}},\boldsymbol{x_{9}},x_{8}, x_2, \boldsymbol{x_{3}})$. 
The first case of ${\rm lk}(x_1)$ implies $f_{seq}(x_{11})=(3.4^2.6)$. Now considering ${\rm lk}(x_{11})$, we get one triangular face adjacent with two quadrangular faces at $x_7$ (see Figure 4.1), which shows $f_{seq}(x_7) \neq f_1, f_2$. Similarly, for the later case of ${\rm lk}(x_1)$, we get a vertex $x$ such that $f_{seq}(x) \neq f_1, f_2$, (see Figure 4.2). So, $f_{seq}(x_1) \neq (3.4^2.6)$. 
If $f_{seq}(x_2)=(3.4^2.6)$, then ${\rm lk}(x_2) = C_9(x_3,\boldsymbol{x_{13}}, x_{12},\boldsymbol{x_{11}},\boldsymbol{x_{10}},\boldsymbol{x_{9}},x_{8}, x_1, \boldsymbol{x_{0}})$, which shows $f_{seq}(x_1)=(3^3.4^2)$ and then we get consecutive two triangular faces and one hexagonal face at $x_8$ (see Figure 4.3), which implies $f_{seq}(x_8) \neq f_1, f_2$. So, $f_{seq}(x_2) \neq f_2$. Similarly, we see that $f_{seq}(x_7) \neq f_2$. This means $x_0$ is not a critical vertex.

\smallskip



\vspace{-.4cm}

	\begin{picture}(0,0)(-20,27)
\setlength{\unitlength}{9mm}

\drawpolygon(1,0)(1,1)(.5,2)(-.5,2)(-1,1)(-1,0)(-2,-.5)(-1.5,-1.5)(-1,-2.5)(0,-2)(1,-2)(1.5,-1)


\drawline[AHnb=0](-1,0)(0,-2)
\drawline[AHnb=0](-1.5,-1.5)(-.5,-1)
\drawline[AHnb=0](0,0)(-.5,-1)
\drawline[AHnb=0](-1,0)(1,0)
\drawline[AHnb=0](-1,1)(1,1)

\drawline[AHnb=0](0,0)(0,1)
\drawline[AHnb=0](0,1)(-.5,2)
\drawline[AHnb=0](0,1)(.5,2)

\put(.2,.75){\scriptsize  $x_0$}
\put(.1,-.25){\scriptsize  $x_1$}
\put(1.15,-.2){\scriptsize  $x_2$}
\put(1.15,.8){\scriptsize  $x_3$}
\put(.5,2.1){\scriptsize  $x_4$}
\put(-.65,2.1){\scriptsize  $x_5$}
\put(-1.4,.9){\scriptsize  $x_6$}
\put(-1.4,.1){\scriptsize  $x_7$}
\put(1.6,-1){\scriptsize  $x_{8}$}
\put(1.1,-2.1){\scriptsize  $x_9$}
\put(0,-2.25){\scriptsize  $x_{10}$}

\put(-.4,-1.1){\scriptsize  $x_{11}$}

\put(-.3,-3){\scriptsize  Figure 4.1}
\end{picture}

\begin{picture}(0,0)(-80,11)
\setlength{\unitlength}{9mm}

\drawpolygon(1,0)(.5,1)(-.5,1)(-1,0)(-2,0)(-2,-1)(-3,-1.5)(-2,-3.5)(-1,-3)(0,-3)(1,-1)


\drawline[AHnb=0](-1,-3)(-2,-1)
\drawline[AHnb=0](-2,-1)(1,-1)
\drawline[AHnb=0](-1,0)(1,0)
\drawline[AHnb=0](0,0)(0,-1)

\drawline[AHnb=0](0,0)(-.5,1)
\drawline[AHnb=0](0,0)(.5,1)
\drawline[AHnb=0](-1,0)(-1,-1)

\drawline[AHnb=0](-1,-1)(-1.5,-2)
\drawline[AHnb=0](-2.5,-2.5)(-1.5,-2)
\drawline[AHnb=0](0,-1)(.5,-2)

\put(.2,.05){\scriptsize  $x_0$}
\put(-.3,-1.3){\scriptsize  $x_1$}
\put(1,-1.3){\scriptsize  $x_2$}
\put(1.1,-.3){\scriptsize  $x_3$}
\put(.6,1.1){\scriptsize  $x_4$}
\put(-1,1){\scriptsize  $x_5$}
\put(-1.4,.2){\scriptsize  $x_6$}
\put(-1,-1.3){\scriptsize  $x_7$}

\put(.6,-2.1){\scriptsize  $x_8$}

\put(-1.3,-2.1){\scriptsize  $x_{11}$}

\put(-1.1,-3.3){\scriptsize  $x_{10}$}

\put(-.1,-3.3){\scriptsize  $x_{9}$}

\put(-2.3,-1){\scriptsize  $x$}

\put(-1,-4.2){\scriptsize  Figure 4.2}
\end{picture}

\begin{picture}(0,0)(-120,11)
\setlength{\unitlength}{9mm}

\drawpolygon(1,0)(.5,1)(-.5,1)(-1,0)(-1,-1)(-.5,-2)(.5,-2)(1,-3)(2,-3)(2.5,-2)(2,-1)(2,0)


\drawline[AHnb=0](0,0)(0,-1)
\drawline[AHnb=0](-1,0)(1,0)
\drawline[AHnb=0](-1,-1)(2,-1)

\drawline[AHnb=0](0,0)(-.5,1)
\drawline[AHnb=0](0,0)(.5,1)

\drawline[AHnb=0](0,-1)(-.5,-2)
\drawline[AHnb=0](0,-1)(.5,-2)

\drawline[AHnb=0](.5,-2)(1,-1)

\drawline[AHnb=0](1,0)(1,-1)

\put(.1,-.2){\scriptsize  $x_0$}
\put(.1,-.9){\scriptsize  $x_1$}
\put(1,-1.3){\scriptsize  $x_2$}
\put(1.1,.1){\scriptsize  $x_3$}
\put(.5,1.1){\scriptsize  $x_4$}
\put(-.6,1.1){\scriptsize  $x_5$}
\put(-1.4,-.2){\scriptsize  $x_6$}
\put(-1.4,-1.1){\scriptsize  $x_7$}

\put(2,.1){\scriptsize  $x_{13}$}
\put(2.15,-1.1){\scriptsize  $x_{12}$}
\put(2.65,-2.1){\scriptsize  $x_{11}$}
\put(2.1,-3.1){\scriptsize  $x_{10}$}
\put(.5,-3){\scriptsize  $x_{9}$}
\put(.2,-2.3){\scriptsize  $x_{8}$}

\put(-.9,-3.7){\scriptsize  Figure 4.3}
\end{picture}

	\vspace{5cm}
\noindent{\bf Case 3.} $f_2 = (3^4.6)$, i.e., $T= [f_1=3^3.4^2: f_2=3^4.6]$. Then at least one vertex in $\{x_4, x_5\}$ has the face-sequence $(3^4.6)$. Without loss of generality let $f_{seq}(x_4) = (3^4.6)$.
Then ${\rm lk}(x_4) = C_8(x_3, x_0, x_5, x_8, x_9, \boldsymbol{x_{10}}, \boldsymbol{x_{11}}$, $\boldsymbol{x_{12}})$ or ${\rm lk}(x_4) = C_8(x_{12}, x_3, x_0, x_5, x_8, \boldsymbol{x_{9}}, \boldsymbol{x_{10}}, \boldsymbol{x_{11}})$ or ${\rm lk}(x_4) = C_8(x_5$, $x_0, x_3, x_8, x_9, \boldsymbol{x_{10}}, \boldsymbol{x_{11}}, \boldsymbol{x_{12}})$. In the first case of ${\rm lk}(x_4)$, we see that $f_{seq}(x_3) \neq f_1, f_2$, see Figure 4.4. In the second case of ${\rm lk}(x_4)$, considering successively ${\rm lk}(x_8)$ and ${\rm lk}(x_5)$, we see that $f_{seq}(x_6) \neq f_1, f_2$, see Figure 4.5. In the last case of ${\rm lk}(x_4)$, considering successively ${\rm lk}(x_9)$, ${\rm lk}(x_{10})$, ${\rm lk}(x_{11})$, ${\rm lk}(x_{12})$, ${\rm lk}(x_5)$, ${\rm lk}(x_6)$, and then ${\rm lk}(x)$, we get a vertex $y$ in ${\rm lk}(6)$ such that $f_{seq}(y) \neq f_1, f_2$, see Figure 4.6. This means $x_0$ is not a critical vertex.

	\begin{picture}(0,0)(-10,33)
	\setlength{\unitlength}{9mm}
	
	\drawpolygon(1,0)(2,0)(2.5,1)(2,2)(0,2)(-1,0)(-1,-1)(1,-1)


	\drawline[AHnb=0](-1,0)(1,0)
	\drawline[AHnb=0](0,0)(0,-1)
	\drawline[AHnb=0](0,0)(-.5,1)
	\drawline[AHnb=0](0,0)(1,2)
	
	\drawline[AHnb=0](-.5,1)(.5,1)
	\drawline[AHnb=0](1,0)(0,2)
	
	\put(.7,.9){\scriptsize  $x_4$}

\put(.2,.05){\scriptsize  $x_0$}
\put(-1.1,-1.3){\scriptsize  $x_7$}
\put(-.1,-1.3){\scriptsize  $x_1$}
\put(.9,-1.3){\scriptsize  $x_2$}
\put(1.1,-.3){\scriptsize  $x_3$}
\put(1.9,-.3){\scriptsize  $x_{12}$}
\put(2.6,.9){\scriptsize  $x_{11}$}
\put(1.9,2.1){\scriptsize  $x_{10}$}
\put(.9,2.1){\scriptsize  $x_9$}
\put(-.1,2.1){\scriptsize  $x_8$}
\put(-.9,1.1){\scriptsize  $x_5$}
\put(-1.4,-.1){\scriptsize  $x_6$}

\put(0,-2){\scriptsize  Figure 4.4}
	\end{picture}

	\begin{picture}(0,0)(-70,29)
\setlength{\unitlength}{9mm}

\drawpolygon(1,0)(2,2)(1.5,3)(-.5,3)(-1,2)(-2,2)(-2.5,1)(-2,0)(-1,0)(-1,-1)(1,-1)


\drawline[AHnb=0](-1,0)(.5,3)
\drawline[AHnb=0](0,0)(0,-1)
\drawline[AHnb=0](0,0)(-1,2)
\drawline[AHnb=0](0,0)(.5,1)

\drawline[AHnb=0](-.5,1)(1.5,1)
\drawline[AHnb=0](1,0)(-.5,3)
\drawline[AHnb=0](-1,0)(1,0)
\drawline[AHnb=0](-1,2)(0,2)

\put(.7,1.1){\scriptsize  $x_4$}

\put(.2,.05){\scriptsize  $x_0$}
\put(-1.1,-1.3){\scriptsize  $x_7$}
\put(-.1,-1.3){\scriptsize  $x_1$}
\put(.9,-1.3){\scriptsize  $x_2$}
\put(1.1,-.3){\scriptsize  $x_3$}
\put(1.7,.9){\scriptsize  $x_{12}$}
\put(2.1,1.9){\scriptsize  $x_{11}$}
\put(1.7,3.1){\scriptsize  $x_{10}$}
\put(.3,3.1){\scriptsize  $x_9$}
\put(.1,2){\scriptsize  $x_8$}
\put(-1,1){\scriptsize  $x_5$}
\put(-1.4,-.2){\scriptsize  $x_6$}

\put(-.5,-1.9){\scriptsize  Figure 4.5}
\end{picture}

	\begin{picture}(0,0)(-135,25)
\setlength{\unitlength}{9mm}

\drawpolygon(1,0)(2,2)(1,4)(-1,4)(-2,2)(-3,2)(-3.5,1)(-3,0)(-2,0)(-2,-1)(1,-1)


\drawline[AHnb=0](-2,0)(1,0)
\drawline[AHnb=0](-1,0)(-1,-1)
\drawline[AHnb=0](0,0)(0,-1)
\drawline[AHnb=0](-1.5,1)(1.5,1)
\drawline[AHnb=0](-2,2)(-1,2)
\drawline[AHnb=0](-1.5,3)(1.5,3)
\drawline[AHnb=0](1,2)(2,2)

\drawline[AHnb=0](-2,0)(0,4)
\drawline[AHnb=0](0,0)(1.5,3)
\drawline[AHnb=0](0,0)(-1.5,3)
\drawline[AHnb=0](-1,0)(-2,2)

\drawline[AHnb=0](1,0)(.5,1)
\drawline[AHnb=0](-1,4)(-.5,3)
\drawline[AHnb=0](0,4)(1.5,1)

\drawline[AHnb=0](.5,3)(1,4)
\drawline[AHnb=0](-1,0)(-.5,1)

\put(.7,1.1){\scriptsize  $x_4$}

\put(.2,.05){\scriptsize  $x_0$}
\put(-1.1,-1.3){\scriptsize  $x_7$}
\put(-.1,-1.3){\scriptsize  $x_1$}
\put(.9,-1.3){\scriptsize  $x_2$}
\put(1.1,-.3){\scriptsize  $x_3$}
\put(1.7,.9){\scriptsize  $x_{8}$}
\put(.5,2){\scriptsize  $x_9$}
\put(0,2.8){\scriptsize  $x_{10}$}
\put(-1.1,2.8){\scriptsize  $x_{11}$}
\put(-.9,1.95){\scriptsize  $x_{12}$}
\put(-1.8,1){\scriptsize  $x$}
\put(-2.25,-.2){\scriptsize  $y$}

\put(-.4,1.1){\scriptsize  $x_5$}
\put(-1.4,-.2){\scriptsize  $x_6$}

\put(-1,-1.9){\scriptsize  Figure 4.6}
\end{picture}

\vspace{4.8cm}


\noindent{\bf Case 4.} $f_2 = (3^2.4.12)$, i.e., $T= [f_1=3^3.4^2: f_2=3^2.4.12]$. Then at least one vertex in $\{x_2,x_3,x_4,x_5,x_6,x_7\}$ has the face-sequence $(3^2.4.12)$. Note that, the case $f_{seq}(x_2) = (3^2.4.12)$ is similar to the case $f_{seq}(x_7) = (3^2.4.12)$; the case $f_{seq}(x_3) = (3^2.4.12)$ is similar to the case $f_{seq}(x_6) = (3^2.4.12)$, and the case $f_{seq}(x_4) = (3^2.4.12)$ is similar to $f_{seq}(x_5) = (3^2.4.12)$. So, it is enough to see the cases when $f_{seq}(x_2) = (3^2.4.12)$, $f_{seq}(x_3) = (3^2.4.12)$, and $f_{seq}(x_4) = (3^2.4.12)$. For 
$f_{seq}(x_2) = (3^2.4.12)$, let 
${\rm lk}(x_2) = C_{14}(x_3, \boldsymbol{x_0}, x_1,x_8,x_9, \boldsymbol{x_{10}}, \boldsymbol{x_{11}}, \boldsymbol{x_{12}}$, $\boldsymbol{x_{13}}, \boldsymbol{x_{14}}, \boldsymbol{x_{15}}, \boldsymbol{x_{16}}, \boldsymbol{x_{17}}, \boldsymbol{x_{18}})$. This implies $f_{seq}(x_3) = (3^2.4.12)$ and $f_{seq}(x_i) = (3^2.4.12)$, for $9 \leq i \leq 18$, which further gives $f_{seq}(x_4) \neq f_1, f_2$. If $f_{seq}(x_3) = (3^2.4.12)$, then $f_{seq}(x_4) \neq f_1, f_2$. If $f_{seq}(x_4) = (3^2.4.12)$, then $f_{seq}(x_3) \neq f_1, f_2$. This means $x_0$ is not a critical vertex.

\smallskip

\noindent{\bf Case 5.} $f_2 = (4.8^2)$, i.e, $T= [f_1=3^3.4^2: f_2=4.8^2]$. Then at least one vertex in $\{x_2,x_7\}$ has the face-sequence $(4.8^2)$. Note that, the case $f_{seq}(x_2) = (4.8^2)$ is similar to the case $f_{seq}(x_7) = (4.8^2)$. If $f_{seq}(x_2) = (4.8^2)$, then $f_{seq}(x_3) \neq f_1,f_2$. This means $x_0$ is not a critical vertex.

\smallskip

\noindent{\bf Case 6.} $f_2 = (4.5.20)$, i.e., $T= [f_1=3^3.4^2: f_2=4.5.20]$. Then at least one vertex in $\{x_2,x_7\}$ has the face-sequence $(4.8^2)$. Note that, the case $f_{seq}(x_2) = (4.5.20)$ is similar to the case $f_{seq}(x_7) = (4.5.20)$. If $f_{seq}(x_2) = (4.5.20)$, then $f_{seq}(x_1) \neq f_1,f_2$. This means $x_0$ is not a critical vertex. 

\smallskip

\noindent{\bf Case 7.} $f_2 = (4.6.12)$, i.e., $T= [f_1=3^3.4^2: f_2=4.6.12]$. Then at least one vertex in $\{x_2,x_7\}$ has the face-sequence $(4.6.12)$. The case $f_{seq}(x_2) = (4.6.12)$ is similar to the case $f_{seq}(x_7) = (4.6.12)$. If $f_{seq}(x_2) = (4.6.12)$, then $f_{seq}(x_1) \neq f_1,f_2$. This means $x_0$ is not a critical vertex. 

\smallskip

\noindent{\bf Case 8.} $f_2 = (3.4.3.12)$, i.e., $T= [f_1=3^3.4^2: f_2=3.4.3.12]$. Then at least one vertex in $\{x_2,x_3,x_6,x_7\}$ has the face-sequence $(3.4.3.12)$. Note that, the cases $f_{seq}(x_2) = (3.4.3.12)$ and $f_{seq}(x_3) = (3.4.3.12)$  are similar to the cases $f_{seq}(x_7) = (3.4.3.12)$ and $f_{seq}(x_6) = (3.4.3.12)$ respectively. If $f_{seq}(x_2) = (3.4.3.12)$, then $f_{seq}(x_3) = (3.4.3.12)$. This gives $f_{seq}(x_4) \neq f_1, f_2$. Similarly, we see $f_{seq}(x_3) \neq f_1, f_2$. This means $x_0$ is not a critical vertex. 

Hence, in all the above cases for $f_2$, $x_0$ is not a critical vertex. Thus the map $M$ of type $T = [f_1=3^3.4^2: f_2]$ does not exist. By a similar computation, it can be shown easily that $M$ does not exist for the remaining type $T \in A$. Thus the claim and hence the lemma. \hfill $\Box$	
		
\smallskip

{\bf Further, we enumerate and classify 2-semi equivelar map $M$ of type $T \in B$ for the number of vertices $|V(M)| \leq 12$, by the Lemmas \ref{l2}-\ref{l7}. The classification is exhaustive search of all possible cases.  Let the vertex set $V(M) = \{0,1, \ldots, 11\}$}.

\begin{lem} \label{l2}
		There exists no 2-semi equivelar map of type $T \in S_2= \{[3^6:3^2.4.12]$, $[3^6:3^2.6^2]$, $[3^4.6:3.6.3.6]$, $[3^2.6^2:3^4.6]$, $[3^2.6^2:3.6.3.6]$, $[3.4^2.6:3.6.3.6]$, $[3.4^2.6:3.4.6.4]$, $[3.4.6.4:4.6.12]$, $[3.12^2:3.4.3.12]\}$ for the number of vertices $ \leq 12$. 
\end{lem}

\noindent{\bf Proof.}  Since one requires more than 12 vertices to complete the link of a vertex with the face-sequence $(3^2.4.12)$ or $(4.6.12)$ or $(3.4.3.12)$, $M$ of type $[3^6:3^2.4.12]$, $[3.4.6.4:4.6.12]$ or $[3.12^2:3.4.3.12]$ does not exist, for the given $V(M)$. For the remaining types, we have following cases:

\smallskip

\noindent {\bf Case 1.} If $M$ is of the type $[3^6:3^2.6^2]$. Without loss of generality, let $0$ be a critical vertex with the face-sequence $(3^2.6^2)$ and let ${\rm lk}(0) = C_{10}(1, \boldsymbol{2}, \boldsymbol{3}, \boldsymbol{4}, 5, 6, 7, \boldsymbol{8}, \boldsymbol{9}, \boldsymbol{10})$. Then the only vertex in ${\rm lk}(0)$ which can have face-sequence $(3^6)$ is 6. This gives ${\rm lk}(6) = C_6(5,0,7,x_1,x_2$, $x_3)$, where $(x_1, x_2, x_3) \in S_3=\{(2,1,10), (2,3,11), (3,2,11), (3,4,11), (4,3,2), (4,3,11), (11, 4,3)\}$. If $(x_1,x_2,x_3)=(2,1,10)$, then ${\rm lk}(6) = C_6(5,0,7,2,1,10)$. This implies  ${\rm lk}(2) = C_{10}(3, \boldsymbol{4}, \boldsymbol{5}, \boldsymbol{0}, 1, 6, 7$, $\boldsymbol{8}, \boldsymbol{x_5}, \boldsymbol{x_4})$. Now by the fact that two distinct hexagonal faces can not share more than 2 vertices, we get $x_4=11$ and then $x_5$ has no value in $V(M)$ so that ${\rm lk}(2)$ can be completed. So, $(x_1,x_2,x_3) \neq (2,1,10)$. By a similar argument, we see that $M$ does not exist for the remaining $(x_1, x_2,x_3) \in S_3$.  

\smallskip

\noindent {\bf Case 2.}  If $M$ is of the type $[3^4.6:3.6.3.6]$. Then, observe that, two distinct hexagonal faces can share at most one vertex. Now, without loss of generality, let $f_{seq}(0)=(3.6.3.6)$ and ${\rm lk}(0) =C_9(1, \boldsymbol{2}, \boldsymbol{3}, \boldsymbol{4}, 5,6, \boldsymbol{7}, \boldsymbol{8},\boldsymbol{9}, 10)$. Then $f_{seq}(1)=(3.6.3.6)$ or $(3^4.6)$.

If $f_{seq}(1)=(3.6.3.6)$, then  ${\rm lk}(1) =C_9(2, \boldsymbol{3}, \boldsymbol{4}, \boldsymbol{5}, 0,10, \boldsymbol{x_1}, \boldsymbol{x_2},\boldsymbol{x_3}, x_4)$, where $x_1, x_2, x_3, x_4 \in V(M)$. 
Note that $x_1 =11$ and then $x_4 \in \{7,8,9\}$. But for all these values of $x_4$, we see that two distinct hexagonal faces share two vertices. So $f_{seq}(1)\neq (3.6.3.6)$.

On the other hand, if $f_{seq}(1)=(3^4.6)$, then ${\rm lk}(1)=C_8(2, \boldsymbol{3}, \boldsymbol{4}, \boldsymbol{5}, 0, 10, 11, x_1)$, where $x_1 \in \{7,8,9\}$. If $x_1=7$, then ${\rm lk}(7)=C_8(8, \boldsymbol{9}, \boldsymbol{10}, \boldsymbol{0}, 6, 11, 1, 2)$ or 
${\rm lk}(7)=C_8(8, \boldsymbol{9}, \boldsymbol{10}, \boldsymbol{0}, 6, 2, 1, 11)$, which implies 
${\rm lk}(2)=C_8(3, \boldsymbol{4}, \boldsymbol{5}, \boldsymbol{0}, 1, 7, 8, x_2)$ or ${\rm lk}(2)=C_8(3, \boldsymbol{4}, \boldsymbol{5}, \boldsymbol{0}, 1, 7, 6, x_2)$ respectively. Observe that for both the cases of ${\rm lk}(2)$, $x_2$ has no value in $V(M)$. So $x_1 \neq 7$. Proceeding similrly, we see that $x_1 \neq 8$ or 9.

\smallskip

\noindent {\bf Case 3.}  If $M$ is of the type $[3^2.6^2:3.6.3.6]$. Then, observe that two distinct hexagonal faces can share at most two vertices. Now, without loss of generality, let 0 be a critical vertex with $f_{seq}(0)=(3^2.6^2)$ such that ${\rm lk}(0) =C_{10}(1, \boldsymbol{2}, \boldsymbol{3}, \boldsymbol{4}, 5,6,7, \boldsymbol{8}, \boldsymbol{9},\boldsymbol{10})$. Then $f_{seq}(6)=(3^6)$. This implies ${\rm lk}(6) = C_6(5,0,7,x_1,x_2,x_3)$. It is easy to see that $(x_1, x_2, x_3) \in \{(2,1,10)$, $(2,3,11)$, $(3,2,11)$, $(3,4,11)$, $(4,3,11)$, $(11,8,9)$, $(11,9,8)$, $(11,9,10)$, $(11,10,9)\}$. If $(x_1,x_2,x_3)=(2,1,10)$, i.e., ${\rm lk}(6) = C_6(5,0,7,2,1,10)$, then ${\rm lk}(2) =C_{10}(3, \boldsymbol{4}, \boldsymbol{5}, \boldsymbol{0}, 1,6,7, \boldsymbol{8}, \boldsymbol{x_4},\boldsymbol{x_5})$. Observe that $x_4=11$ and $x_5$ has no value in $V(M)$. So for $(x_1, x_2, x_3) = (2,1,10)$, we do not get $M$. Proceeding similarly, we see that $M$ does not exist for the remaining values of $(x_1,x_2,x_3)$.

\smallskip

\noindent {\bf Case 4.}  If $M$ is of the type $[3.4^2.6:3.6.3.6]$. Then, observe that, two distinct hexagonal faces can share at most one vertex. Let $f_{seq}(0)=(3.6.3.6)$ and ${\rm lk}(0) =C_{10}(1, \boldsymbol{2}, \boldsymbol{3}, \boldsymbol{4}, 5,6, \boldsymbol{7}, \boldsymbol{8}, \boldsymbol{9},10)$. Then $f_{seq}(1)=(3.6.3.6)$ or $(3.4^2.6)$. In the first case, we get ${\rm lk}(1) =C_{10}(2, \boldsymbol{3}, \boldsymbol{4}, \boldsymbol{5}, 0,10, \boldsymbol{x_1}, \boldsymbol{x_2}, \boldsymbol{x_3},x_4)$, for $x_1,x_2,x_3,x_4 \in V(M)$. Now following the fact that two distinct hexagonal faces can share at most one vertex, we see that ${\rm lk}(1)$ can not be completed. On the other hand, if $f_{seq}(1)=(3.4^2.6)$, then 
${\rm lk}(1) =C_9(0, \boldsymbol{5}, \boldsymbol{4}, \boldsymbol{3}, 2,\boldsymbol{x_3},x_2, \boldsymbol{x_1},10)$. In this case, it is easy to see that $x_2=11$, this gives $x_1=9$ and $x_3 \in \{6,7,8\}$. If $x_3=6$, then ${\rm lk}(11) =C_9(x_4, \boldsymbol{x_5}, \boldsymbol{x_6}, \boldsymbol{x_7}, 6,\boldsymbol{2},1, \boldsymbol{10},9)$ or ${\rm lk}(11) =C_9(x_4, \boldsymbol{x_5}, \boldsymbol{x_6}, \boldsymbol{x_7}, 9,\boldsymbol{10},1, \boldsymbol{2},6)$. Again by the fact that two distinct quadrangular faces can share at most 1 vertex, we see that in both the cases ${\rm lk}(11)$ can not be completed. Proceeding similarly for $x_3=7$ and $6$, we see that ${\rm lk}(11)$ can not be completed. 

\smallskip

\noindent {\bf Case 5.}  If $M$ is of the type $[3.4^2.6:3.4.6.4]$. Then, observe that, two distinct hexagonal faces can not share any vertex. Let $f_{seq}(0)=(3.4^2.6)$ and ${\rm lk}(0) =C_{9}(1, \boldsymbol{2}, \boldsymbol{3}, \boldsymbol{4}, 5, \boldsymbol{6},7, \boldsymbol{8},9)$. This implies 
${\rm lk}(7) =C_{9}(x_1, \boldsymbol{x_2}, \boldsymbol{x_3}, \boldsymbol{x_4}, 6, \boldsymbol{5},0, \boldsymbol{9},8)$ or ${\rm lk}(7) =C_{9}(x_1, \boldsymbol{x_2}, \boldsymbol{x_3}, \boldsymbol{x_4}, 8, \boldsymbol{9},0, \boldsymbol{5},6)$. Then by the fact that two distinct hexagonal faces are disjoint, successively, we see $x_1=10$ and $x_2=11$. Now observe that $x_3,x_4$ have no values in $V(M)$ so that ${\rm lk}(7)$ can be completed.

\smallskip

\noindent {\bf Case 6.}  If $M$ is of the type $[3^4.6:3^2.6^2]$. Then, observe that, two distinct hexagonal faces can share at most two adjacent vertices (an edge). Let $f_{seq}(0)=(3^2.6^2)$ and ${\rm lk}(0) =C_{10}(1, \boldsymbol{2}, \boldsymbol{3}, \boldsymbol{4}, 5, 6, 7, \boldsymbol{8}, \boldsymbol{9}$, $\boldsymbol{10})$. This implies ${\rm lk}(1) =C_{10}(0, \boldsymbol{5}, \boldsymbol{4}, \boldsymbol{3}, 2, 11, 10, \boldsymbol{9}, \boldsymbol{8}, \boldsymbol{7})$. Then $f_{seq}(6) = (3^2.6^2)$ or $(3^4.6)$. If $f_{seq}(6)=(3^2.6^2)$, then by the above observation, we see easily that ${\rm lk}(6)$ can not be completed. on the other hand, If $f_{seq}(6)=(3^4.6)$, then 
${\rm lk}(6) =C_8(5, \boldsymbol{x_1}, \boldsymbol{x_2}, \boldsymbol{x_3}, x_4, x_5, 7,0)$ or 
${\rm lk}(6) =C_8(7, \boldsymbol{x_1}, \boldsymbol{x_2}$, $\boldsymbol{x_3}, x_4, x_5, 5,0)$ or 
${\rm lk}(6) =C_8(x_1, \boldsymbol{x_2}, \boldsymbol{x_3}, \boldsymbol{x_4}, x_5, 7, 0,5)$. The first case of ${\rm lk}(6)$ implies $x_1=4$. Then there exist two distinct hexagonal faces at $5$ but $C_6(0,4,x_2,x_3,x_4,6) \subseteq {\rm lk}(5)$, which can not be true. Similarly, we see that the second case of ${\rm lk}(6)$ is also not possible. In the last case of ${\rm lk}(6)$, it is easy to see that $(x_1,x_2,x_3,x_4,x_5) \in \{(8,9,3,2,11)$, $(8,9,4,3,11)$, $(8,9,11,2,3)$, $(8,9,11,3,4)$, $(8,9, 11,4,3)$, $(9,8,3,2,11)$, $(9,8,11,2,3)$, $(9,8,11,3,4)$, $(9,8,11,4,3)$, $(11,10,9,3,4)$, $(11,10,9,4,3))\}$. If $(x_1,x_2,x_3,x_4,x_5)=(8,9,3,2,11)$, then ${\rm lk}(8) =C_{10}(9, \boldsymbol{10}, \boldsymbol{1}, \boldsymbol{0}, 7, 4, 6, \boldsymbol{11}, \boldsymbol{2}, \boldsymbol{3})$ and ${\rm lk}(9) =C_{10}(8, \boldsymbol{7}, \boldsymbol{0}, \boldsymbol{1}, 10, 4, 3, \boldsymbol{2}, \boldsymbol{11}, \boldsymbol{6})$. This implies $C_9(0,1,2,11,6,8,9,4,5) \subseteq {\rm lk}(3)$. Similarly, we see that if $(x_1,x_2,x_3,x_4,x_5)=(8,9,4,3,11)$, then considering ${\rm lk}(8)$, ${\rm lk}(9)$ can not be completed. If $(x_1,x_2,x_3,x_4,x_5)=(8,9,11,2,3)$, then considering ${\rm lk}(8)$, we get $C_9(0,1,11,9,8,6,3,4,5) \subseteq {\rm lk}(2)$.
If $(x_1,x_2,x_3,x_4,x_5) = (8,9,11,3,4)$, then considering ${\rm lk}(8)$, we see that ${\rm lk}(9)$ can not be completed. If $(x_1,x_2,x_3,x_4,x_5) = (8,9,11,4,3)$, then considering ${\rm lk}(8)$ and ${\rm lk}(9)$ successively, we see that ${\rm lk}(4) =C_{10}(3, \boldsymbol{2}, \boldsymbol{1}, \boldsymbol{0}, 5, 10, 11, \boldsymbol{9}, \boldsymbol{8}, \boldsymbol{6})$, this implies ${\rm lk}(5)$ can not be completed. 
If $(x_1,x_2,x_3,x_4,x_5)=(9,8,3,2,11)$, then considering ${\rm lk}(9)$, we see ${\rm lk}(8)$ can not be completed. If $(x_1,x_2,x_3,x_4,x_5)=(9,8,11,2,3)$, then ${\rm lk}(9) =C_{10}(8, \boldsymbol{7}, \boldsymbol{0}, \boldsymbol{1}, 10, 5, 6, \boldsymbol{3}, \boldsymbol{2}, \boldsymbol{11})$, ${\rm lk}(8) =C_{10}(9$, $\boldsymbol{10}, \boldsymbol{1}, \boldsymbol{0}, 7, 4, 11, \boldsymbol{2}, \boldsymbol{3}, \boldsymbol{6})$. This implies  ${\rm lk}(7) =C_8(0, \boldsymbol{1}, \boldsymbol{10}, \boldsymbol{9}, 8, 4, 3,6)$ and ${\rm lk}(4) =C_8(3, \boldsymbol{2}, \boldsymbol{1}, \boldsymbol{0}, 5, 11$, $8,7)$. Now we observe that ${\rm lk}(5)$ can not be completed. If $(x_1,x_2,x_3,x_4,x_5)=(9,8,11,3,4)$, then ${\rm lk}(9) =C_{10}(8, \boldsymbol{7}, \boldsymbol{0}, \boldsymbol{1}, 10, 5, 6, \boldsymbol{4}, \boldsymbol{3}, \boldsymbol{11})$, ${\rm lk}(8) =C_{10}(9$, $\boldsymbol{10}, \boldsymbol{1}, \boldsymbol{0}, 7, 2, 11, \boldsymbol{3}, \boldsymbol{4}, \boldsymbol{6})$. This implies  ${\rm lk}(2) =C_8(1, \boldsymbol{0}, \boldsymbol{5}, \boldsymbol{4}, 3, 7, 8,11)$. Now we see that ${\rm lk}(7)$ can not be completed. 
If $(x_1,x_2,x_3,x_4,x_5)=(11,10,9,3,4)$, then ${\rm lk}(4) =C_{10}(3, \boldsymbol{2}, \boldsymbol{1}, \boldsymbol{0}, 5, 7, 6, \boldsymbol{11}, \boldsymbol{10}, \boldsymbol{9})$ and ${\rm lk}(3) =C_{10}(4, \boldsymbol{5}, \boldsymbol{0}, \boldsymbol{1}, 2, 8, 9, \boldsymbol{10}, \boldsymbol{11}$, $\boldsymbol{6})$. This implies $C_9(0,1,10,11,6,4,3,8,7) \subseteq {\rm lk}(9)$. If $(x_1,x_2,x_3,x_4,x_5)=(11,10,9,4,3)$, then ${\rm lk}(3) =C_{10}(4, \boldsymbol{5}, \boldsymbol{0}, \boldsymbol{1}, 2, 7, 6, \boldsymbol{11}, \boldsymbol{10}, \boldsymbol{9})$ and ${\rm lk}(4) =C_{10}(3, \boldsymbol{2}, \boldsymbol{1}, \boldsymbol{0}, 5, 8, 9, \boldsymbol{10}, \boldsymbol{11}, \boldsymbol{6})$. This implies $C_9(0,1,10$, $11,6,3,4,8,7) \subseteq {\rm lk}(9)$. If $(x_1,x_2,x_3,x_4,x_5)=(9,8,11,4,3)$, then ${\rm lk}(9) =C_{10}(8, \boldsymbol{7}, \boldsymbol{0}, \boldsymbol{1}$, $10, 5, 6, \boldsymbol{3}, \boldsymbol{4}, \boldsymbol{11})$, ${\rm lk}(8) =C_{10}(9$, $\boldsymbol{10}, \boldsymbol{1}, \boldsymbol{0}, 7, 2, 11, \boldsymbol{4}, \boldsymbol{3}, \boldsymbol{6})$. Now successively, we get ${\rm lk}(7) =C_8(0, \boldsymbol{1}$, $\boldsymbol{10}, \boldsymbol{9}, 8, 2, 3,6)$, ${\rm lk}(2) =C_8(1, \boldsymbol{0}, \boldsymbol{5}$, $\boldsymbol{4}, 3, 7, 8,11)$, ${\rm lk}(11) =C_8(4, \boldsymbol{3}, \boldsymbol{6}, \boldsymbol{9}, 8, 2, 1,10)$, ${\rm lk}(10) =C_8(1, \boldsymbol{0}$, $\boldsymbol{7}, \boldsymbol{8}, 9, 5, 4,11)$, and ${\rm lk}(5) =C_8(0, \boldsymbol{1}$, $\boldsymbol{2}, \boldsymbol{3}, 4, 10, 9,6)$. Then observe that ${\rm lk}(3)$ and ${\rm lk}(4)$ can not be completed. This completes the proof. \hfill$\Box$

\begin{lem} \label{l3}
	Let $M$ be a 2-semi equivelar map of type $[3^6:3^4.6]$ on $\leq 12$ vertices. Then $M$ is isomorphic to one of $\mathcal{A}_1(K)$, $\mathcal{A}_2(T)$ or $\mathcal{A}_3(T)$ given in Section \ref{s3}. 
\end{lem}

\noindent{\bf Proof.} Without loss of generality, let $0$ be a vertex of face-sequence $(3^4.6)$ such that ${\rm lk}(0)=C_8(1, \boldsymbol{2}, \boldsymbol{3}, \boldsymbol{4}$, $5,6,7,8)$. Since the vertices $1,2,3,4,5$ appear in hexagonal face $[0,1,2,3,4,5]$, the vertices have the face-sequence $(3^4.6)$. This gives ${\rm lk}(1)=C_8(2, \boldsymbol{3}, \boldsymbol{4}, \boldsymbol{5},0,8,x_1,x_2)$. It is easy to see that $(x_1,x_2) \in \{(6,7)$, $(6,9)$, $(9,6)$, $(9,10)$, $(9,7)\}$. The case $(6, 9) \cong (9,7)$ by the map $(0,1)(2,5)(3,4)(6,7,9)$. So, we need not consider the last case.

\smallskip

\noindent{\bf Claim 1.} $(x_1, x_2) = (9, 6)$ or $(9, 10)$.

If $(x_1,x_2)=(6,7)$, then ${\rm lk}(1)=C_8(2, \boldsymbol{3}, \boldsymbol{4}, \boldsymbol{5},0,8,6,7)$ and this gives $f_{seq}(6)$ is $(3^6)$ or $(3^4.6)$. By the fact that two distinct hexagonal faces are disjoint, we see that for the given number of vertices, ${\rm lk}(6)$ can not be completed if $f_{seq}(6) = (3^4.6)$. So, $f_{seq}(6)=(3^6)$. This gives ${\rm lk}(6)=C_6(0,5,9,8,1,7)$. Again by the same fact as above, we see $f_{seq}(7)=(3^6)$ and then we get ${\rm lk}(7)=C_6(0,6,1,2,x_3,8)$, where $x_3 \in \{9,10\}$. If $x_3=9$, then $C_5(0,1,6,9,7) \subseteq {\rm lk}(8)$, which can not be true. If $x_3=10$, then ${\rm lk}(7)=C_6(0,6,1,2,10,8)$. This implies ${\rm lk}(8)=C_6(0,1,6,9,10,7)$, ${\rm lk}(2)=C_8(3, \boldsymbol{4}, \boldsymbol{5}, \boldsymbol{0},1,7,10,11)$. Now, considering the same fact as above, we see $f_{seq}(10)=(3^6)$. This gives ${\rm lk}(10)=C_6(9,8,7,2,11,x_4)$. But observe that $x_4$ has no value in $V(M)$ so that ${\rm lk}(10)$ can be completed. So $(x, y) \neq (6, 7)$.

If $(x_1,x_2)=(6,9)$, then ${\rm lk}(1)=C_8(2, \boldsymbol{3}, \boldsymbol{4}, \boldsymbol{5},0,8,6,9)$,  ${\rm lk}(6)=C_6(0,5,8,1,9,7)$, ${\rm lk}(5)=C_8(0, \boldsymbol{1}, \boldsymbol{2}$, $\boldsymbol{3},4,10,8,6)$, ${\rm lk}(8)=C_6(0,1,6,5,10,7)$. This gives $f_{seq}(7)=(3^6)$, which implies ${\rm lk}(7)=C_6(9,6,0$, $8,10,x_3)$, where $x_3 \in \{3,4,11\}$. If $x_3=a$, for $a \in \{3,11\}$, we see that the $f_{seq}(9) = (3^6)$, otherwise the vertices 2,3 lie in two distinct hexagonal faces which is not allowed. So ${\rm lk}(9)=C_6(2,1,6,7,l,x_4)$. Now observe that $x_4$ has no value in $V(M)$ so that ${\rm lk}(9)$ can be completed. If $x_4=4$, then ${\rm lk}(7)=C_6(9,6,0,8,10,4)$ and ${\rm lk}(4)=C_8(5, \boldsymbol{0}, \boldsymbol{1}, \boldsymbol{2},3,9,7,10)$. This implies $C_6(0,1,2,9,4,5) \subseteq{\rm lk}(3)$. This can not be true as $f_{seq}(3)=(3^4.6)$. This proves the claim.

\smallskip

\noindent{\bf Case 1.} $(x_1,x_2)=(9,6)$, i.e., ${\rm lk}(1)=C_8(2, \boldsymbol{3}, \boldsymbol{4}, \boldsymbol{5},0,8,9,6)$. Then  ${\rm lk}(6)=C_6(0,5,9,1,2,7)$ and ${\rm lk}(2)=C_8(3, \boldsymbol{4}, \boldsymbol{5}, \boldsymbol{0},1,6,7,x_3)$, where $x_3 \in \{8, 10\}$. 
If $x_3=8$, then ${\rm lk}(2)=C_8(3, \boldsymbol{4}, \boldsymbol{5}, \boldsymbol{0}$, $1,6,7,8)$, ${\rm lk}(8)=C_6(0,1,9,3,2,7)$, ${\rm lk}(9)=C_6(1,6,5,10,3,8)$, ${\rm lk}(3)=C_8(4, \boldsymbol{5}, \boldsymbol{0}, \boldsymbol{1},2,8,9,10)$, ${\rm lk}(5)=C_8(0, \boldsymbol{1}, \boldsymbol{2}, \boldsymbol{3},4,10,9,6)$. This implies $C_4(3,4,5,9) \subseteq {\rm lk}(10)$. So $x_3=10$. Then ${\rm lk}(2)=C_8(3, \boldsymbol{4}, \boldsymbol{5}, \boldsymbol{0},1,6,7,10)$ and ${\rm lk}(7)=C_6(8,0,6,2,10,x_4)$, where $x_4 \in \{4,11\}$. In case $x_4=11$, for ${\rm lk}(8)=C_6(9,1,0,7,11,x_5)$, we get $x_5=4$. But then considering ${\rm lk}(9)$ and ${\rm lk}(3)$, we see that $\deg(5)>5$. So, $x_4=4$. Then it is easy to see that ${\rm lk}(8)=C_6(9,1,0,7,4,3)$.  Now completing successively, we get ${\rm lk}(4)=C_8(5, \boldsymbol{0}, \boldsymbol{1}, \boldsymbol{2},3,8,7,10)$, ${\rm lk}(5)=C_8(0, \boldsymbol{1}, \boldsymbol{2}, \boldsymbol{3},4,10,9,6)$, ${\rm lk}(9)=C_6(1,6,5,10,3,8)$, ${\rm lk}(3)=C_8(4, \boldsymbol{5}, \boldsymbol{0}, \boldsymbol{1},2,10,9,8)$, and ${\rm lk}(10)=C_6(2,3,9,5,4,7)$. Then $M \cong \mathcal{A}_1 (K)$ by the identity map. 
 
\smallskip

\noindent{\bf Case 2.} $(x_1,x_2)=(9,10)$, i.e., ${\rm lk}(1)=C_8(2, \boldsymbol{3}, \boldsymbol{4}, \boldsymbol{5},0,8,9,10)$. This implies ${\rm lk}(2)=C_8(3, \boldsymbol{4}, \boldsymbol{5}, \boldsymbol{0},1$, $10,x_3,x_4)$. It is easy to see that $(x_3,x_4) \in \{(6,7)$, $(6,11)$, $(7,6)$, $(7,8)$, $(7,9)$, $(7,11)$, $(8,7)$, $(8,9)$, $(11,6)$, $(11,7)$, $(11,9)\}$. Here, the cases $(6,7) \cong (7,6)$ and $(6, 11) \cong (11,7)$ by the maps $(0,2)(3,5)(6,7)(8,10)$ and $(0,2)(3,5)(6,7,11)(8,10)$ respectively. So, we need not to discuss the cases $(x_3, x_4) = (6,7)$ or $(11,7)$.

\smallskip

\noindent{\bf Claim 2.} $(x_3, x_4) = (6,11)$, $(7,6)$ or $(7,8)$. 

If $(x_3,x_4)=(7,9)$, then ${\rm lk}(9)=C_6(1,8,7,2,3,10)$, ${\rm lk}(7)=C_6(0,6,10,2,9,8)$, ${\rm lk}(10)=C_6(1,2$, $7,6,3,9)$, and ${\rm lk}(3)=C_8(4, \boldsymbol{5}, \boldsymbol{0}, \boldsymbol{1},2,9,10,6)$. Now considering ${\rm lk}(6)$, we see that the set $\{3,4,5\}$ forms a triangular face, which is not true in ${\rm lk}(0)$. 

If $(x_3,x_4)=(7,11)$, then ${\rm lk}(2)=C_8(3, \boldsymbol{4}, \boldsymbol{5}, \boldsymbol{0},1,10,7,11)$. This implies ${\rm lk}(7)=C_6(0,6,10,2$, $11,8)$ or ${\rm lk}(7)=C_6(0,6,11,2,10,8)$. In the first case, ${\rm lk}(10)=C_6(1,2,7,6,4,9)$, ${\rm lk}(4)=C_8(5, \boldsymbol{0}, \boldsymbol{1}$, $\boldsymbol{2},3,6,10,9)$. This implies ${\rm lk}(6)=C_6(0,5,3,4,10,7)$, which contradicts the fact that $\{5,6\}$ forms a non-edge in ${\rm lk}(0)$. 
So ${\rm lk}(7)=C_6(0,6,11,2,10,8)$. Then ${\rm lk}(10)=C_6(1,2,7,8,4,9)$, ${\rm lk}(8)=C_6(0,1,9,4,10,7)$, this gives $C_4(1,8,4$, $10) \subseteq {\rm lk}(9)$. 

If $(x_3,x_4)=(8,7)$, then ${\rm lk}(2)=C_8(3, \boldsymbol{4}, \boldsymbol{5}, \boldsymbol{0},1,10,8,7)$. This implies ${\rm lk}(8)=C_6(0,1,9,10,2,7)$, and we get triangular face $[8,9,10]$, which is not true if we consider ${\rm lk}(1)$. 

If $(x_3,x_4)=(8,9)$, then ${\rm lk}(2)=C_8(3, \boldsymbol{4}, \boldsymbol{5}, \boldsymbol{0},1,10,8,9)$. This implies ${\rm lk}(8)=C_6(0,1,9,2,10,7)$, and ${\rm lk}(10)=C_6(7,8,2,1,9,x_5)$, where $x_5 \in \{4,11\}$. If $x_5=4$, then ${\rm lk}(9)=C_6(1,8,2,3,4,10)$ and we get $C_6(0,1,2,9,4,5) \subseteq {\rm lk}(3)$. If $x_5=11$, then successively, we get ${\rm lk}(9)=C_6(1,8,2,3,11,10)$, ${\rm lk}(7)=C_6(0,6,4,11,10,8)$. Observe that ${\rm lk}(4)=C_8(5, \boldsymbol{0}, \boldsymbol{1}, \boldsymbol{2},3,6,7,11)$, this implies ${\rm lk}(3)=C_8(4, \boldsymbol{5}, \boldsymbol{0}, \boldsymbol{1},2,9,11,6)$, and we see $\deg(11)>6$.

If $(x_3,x_4)=(11,6)$, then ${\rm lk}(2)=C_8(3, \boldsymbol{4}, \boldsymbol{5}, \boldsymbol{0},1,10,11,6)$. This implies ${\rm lk}(6)=C_6(0,5,11,2,3,7)$, ${\rm lk}(3)=C_8(4, \boldsymbol{5}, \boldsymbol{0}, \boldsymbol{1},2,6,7,9)$, and ${\rm lk}(7)=C_6(8,0,6,3,9,x_5)$, where $x_5 \in \{4,10\}$. If $x_5=4$, considering ${\rm lk}(4)$ and ${\rm lk}(8)$, we see $\deg(9)>6$. If $x_5=10$, then considering ${\rm lk}(10)$, ${\rm lk}(8)$ and ${\rm lk}(9)$ successively, we see $\deg(3)>5$.


If $(x_3,x_4)=(11,9)$, then ${\rm lk}(2)=C_8(3, \boldsymbol{4}, \boldsymbol{5}, \boldsymbol{0},1,10,11,9)$. This gives ${\rm lk}(9)=C_6(1,8,11,2,3,10)$, and ${\rm lk}(10)=C_6(3,9,1,2,11,x_5)$. Here $x_5 \in \{4,6,7\}$. If $x_5=4$, then $C_7(0,1,2,9,10,4,5) \subseteq{\rm lk}(3)$. If $x_5 = a \in \{6,7\}$, then completing ${\rm lk}(6)$, we see that $\deg(a) >6$. This proves the claim. 

\smallskip

\noindent{\bf Subcase 2.1.} If $(x_3,x_4)=(6,11)$, then  ${\rm lk}(6)=C_6(0,5,10,2,11,7)$ or ${\rm lk}(6)=C_6(0,5,11,2,10,7)$. In the first case, ${\rm lk}(10)$ can not be completed. So, ${\rm lk}(6)=C_6(0,5,11,2,10,7)$. Then ${\rm lk}(10)=C_6(1,2,6,7,4,9)$ and ${\rm lk}(7)=C_6(0,6,10,4,3,8)$. Completing successively, we get ${\rm lk}(4)=C_8(5, \boldsymbol{0}, \boldsymbol{1}$, $\boldsymbol{2},3,7,10,9)$, ${\rm lk}(5)=C_8(0, \boldsymbol{1}, \boldsymbol{2}, \boldsymbol{3},4,9,11,6)$, ${\rm lk}(3)=C_8(4, \boldsymbol{5}, \boldsymbol{0}, \boldsymbol{1},2,11,8,7)$, ${\rm lk}(8)=C_6(0,1,9,11$, $3,7)$, ${\rm lk}(9)=C_6(1,8,11,5,4,10)$. Then $M \cong \mathcal{A}_2 (T)$ by the identity map.

\smallskip

\noindent{\bf Subcase 2.2.} If $(x_3,x_4)=(7,6)$, then  ${\rm lk}(2)=C_8(3, \boldsymbol{4}, \boldsymbol{5}, \boldsymbol{0},1,10,7,6)$ and ${\rm lk}(6)=C_6(3,2,7,0,5,x_5)$, where $x_5 \in \{9, 11\}$. If $x_5=11$, then ${\rm lk}(3)=C_8(4, \boldsymbol{5}, \boldsymbol{0}, \boldsymbol{1},2,6,11,9)$ and we see that ${\rm lk}(7)$ can not be completed. So $x_5=9$. Then ${\rm lk}(6)=C_6(3,2,7,0,5,9)$, this implies, ${\rm lk}(9)=C_6(1,8,5,6,3,10)$ or ${\rm lk}(9)=C_6(1,8,3,6,5,10)$. In the first case of ${\rm lk}(9)$, completing successively we get ${\rm lk}(5)=C_8(0, \boldsymbol{1}, \boldsymbol{2}, \boldsymbol{3},4,8,9,6)$, ${\rm lk}(8)=C_6(0,1,9,5,4,7)$, ${\rm lk}(7)=C_6(0,6,2,10,4,8)$, ${\rm lk}(4)=C_8(5, \boldsymbol{0}, \boldsymbol{1}, \boldsymbol{2},3$, $10,7,8)$, ${\rm lk}(3)=C_8(4, \boldsymbol{5}, \boldsymbol{0}, \boldsymbol{1},2,6,9,10)$, and ${\rm lk}(10)=C_6(1,2,7,4,3,9)$. Then $M \cong \mathcal{C}_1 (K)$ by the map $(1,5)(2,4)(6,8)$. Also, when ${\rm lk}(9)=C_6(1,8,3,6,5,10)$, completing successively we get ${\rm lk}(5)=C_8(0, \boldsymbol{1}, \boldsymbol{2}, \boldsymbol{3},4,10,9,6)$, ${\rm lk}(10)=C_6(1,2,7,4,5,9)$, ${\rm lk}(7)=C_6(0,6,2,10,4,8)$, ${\rm lk}(4)=C_8(5, \boldsymbol{0}, \boldsymbol{1}, \boldsymbol{2},3,8,7,10)$, ${\rm lk}(3)=C_8(4, \boldsymbol{5}, \boldsymbol{0}, \boldsymbol{1},2,6,9,8)$. Then $M \cong \mathcal{A}_3 (T)$ by the identity map.

\smallskip
\noindent{\bf Subcase 2.3.} If $(x_3,x_4)=(7,8)$, then ${\rm lk}(8)=C_6(0,1,9,3,2,7)$ and ${\rm lk}(3)=C_8(4, \boldsymbol{5}, \boldsymbol{0}, \boldsymbol{1},2,8,9,x_5)$. Observe that $x_5=6$. Now completing successively, we get  ${\rm lk}(3)=C_8(4, \boldsymbol{5}, \boldsymbol{0}, \boldsymbol{1},2,8,9,6)$, ${\rm lk}(6)=C_6(0,5,9,3,4,7)$, ${\rm lk}(4)=C_8(5, \boldsymbol{0}, \boldsymbol{1}, \boldsymbol{2},3,6,7,10)$, ${\rm lk}(10)=C_6(1,2,7,4,5,9)$, ${\rm lk}(5)=C_8(0, \boldsymbol{1}, \boldsymbol{2}, \boldsymbol{3},4$, $10,9,6)$, ${\rm lk}(9)=C_6(1,8,3,6,5,10)$, and ${\rm lk}(7)=C_6(0,6,4,10,2,8)$. Then $M \cong \mathcal{A}_1 (K)$ by the map $(0,5,4,3,2,1)(6,10,8)(7,9)$. Thus the proof. \hfill$\Box$

\begin{lem}\label{l3.5}
	Let $M$ be a 2-semi equivelar map of type $[3^3.4^2:3.4.6.4]$ on $\leq 12$ vertices. Then $M$ is isomorphic to $\mathcal{B}_1(K)$ or $ \mathcal{B}_2(T)$ given in example Section \ref{s3}. 
\end{lem}

\noindent {\bf Proof.} If $M$ is of the type $[3^3.4^2:3.4.6.4]$. Then, observe that: (i) if $i$ is a vertex with the face-sequence $(3^3.4^2)$ or $(3^2.4.3.4)$, then $\deg(i)=5$ or 4 respectively (ii) two distinct hexagonal faces do not share any vertex, and (iii) the non-empty intersection of a hexagonal face and a quadrangular face is an edge. Without loss of generality, let $f_{seq}(0)=(3.4.6.4)$ and ${\rm lk}(0) =C_9(1, \boldsymbol{2}, \boldsymbol{3}, \boldsymbol{4}, 5, \boldsymbol{6}, 7,8, \boldsymbol{9})$.

\smallskip

\noindent{\bf Claim 1.} The vertices $6,7,8,9$ can not have the face-sequence $(3.4.6.4)$.

Note that the cases, $f_{seq}(6)=(3.4.6.4)$ and $f_{seq}(7)=(3.4.6.4)$ are similar to the cases $f_{seq}(9)=(3.4.6.4)$ and $f_{seq}(8)=(3.4.6.4)$ respectively. So, first assume that $f_{seq}(6)=(3.4.6.4)$. Then ${\rm lk}(6) =C_9(5, \boldsymbol{y_1},\boldsymbol{y_2}, \boldsymbol{y_3}, y_4, \boldsymbol{y_5}, y_6, 7,\boldsymbol{0})$ or ${\rm lk}(6) =C_9(7, \boldsymbol{y_1}, \boldsymbol{y_2}, \boldsymbol{y_3}, y_4, \boldsymbol{y_5}, y_6, 5, \boldsymbol{0})$, where $y_1, y_2, y_3, y_4,y_5$, $y_6 \in V(M)$. In the first case of ${\rm lk}(6)$, we see that 5 appears in two distinct hexagonal faces $[0,1,2,3,4,5]$ and $[5,y_1,y_2,y_3,y_4]$, which is not allowed. While for the later case of ${\rm lk}(6)$, we see that $y_1, y_2, y_3, y_4 \neq 8$, as if $y_1=8$, then it contradicts the fact that $[0,7,8]$ is a triangular face and if any $x_j = 8$, for $j \in \{2,3,4\}$, then it contradicts the fact that $7\,8$ is an edge. Now, observe that we have four vertices $y_1, y_2, y_3, y_4$ in ${\rm lk}(6)$ which have possible choices from the set $\{9, 10,11\}$. This is not possible. So $f_{seq}(6) \neq (3.4.6.4)$. Similarly, we see that $f_{seq}(7) \neq (3.4.6.4)$. Thus the claim.

\smallskip 

From ${\rm lk}(0)$, we get ${\rm lk}(1) =C_9(2, \boldsymbol{3}, \boldsymbol{4}, \boldsymbol{5}, 0, \boldsymbol{8}, 9, x_1, \boldsymbol{x_2})$, where $x_1, x_2 \in V(M)$. It is easy to see that $(x_1, x_2) \in \{(6,7)$, $(6,10)$, $(7,6)$, $(7,10)$, $(10,6)$, $(10,7)$, $(10,11)\}$. 

\smallskip

\noindent{\bf Claim 2.} $(x_1, x_2) = (10,11)$.

By  Claim 1, we know if there exist two quadrangular faces at 6, 7, 8 or 9, then these faces share an edge. It follows that $(x_1, x_2) \not \in \{(6, 10)$, $(7, 10)$, $(10, 6)$, $(10, 7)\}$. If $(7, 6)$, then ${\rm lk}(7) = C_7(6, \boldsymbol{2}, 1, 9, 8, 0, \boldsymbol{5})$. This implies $C_4(0,1,9,7) \subseteq {\rm lk}(8)$. If $(x_1, x_2)=(6, 7)$, then ${\rm lk}(6) = C_7(7, \boldsymbol{2}, 1, 9$, $10, 5, \boldsymbol{0})$, ${\rm lk}(7) = C_7(6, \boldsymbol{1}, 2, 11, 8, 0, \boldsymbol{5})$ and ${\rm lk}(2) =C_9(1, \boldsymbol{0}, \boldsymbol{5}, \boldsymbol{4}, 3, \boldsymbol{10}, 11, 7, \boldsymbol{6})$. This implies ${\rm lk}(10) = C_7(3, \boldsymbol{2}, 11, 9, 6, 5, \boldsymbol{x_3})$ or ${\rm lk}(10) = C_7(11, \boldsymbol{2}, 3, 5, 6, 9, \boldsymbol{x_3})$. In the first case of ${\rm lk}(10)$, we see that the faces $[0,1,2,3,4,5]$ and $[3,10,5,x_3]$ share two non-adjacent vertices $\{3,5\}$, which contradicts (iii) given above. In the later case of ${\rm lk}(10)$, we see that the set $\{3,5\}$ forms an edge and non-edge both. This proves the claim. 

Let $(x_1, x_2)= (10,11)$, i.e., ${\rm lk}(1) =C_9(2, \boldsymbol{3}, \boldsymbol{4}, \boldsymbol{5}, 0, \boldsymbol{8}, 9, 10, \boldsymbol{11})$, then 
${\rm lk}(2) =C_9(3, \boldsymbol{4}, \boldsymbol{5}, \boldsymbol{0}, 1, \boldsymbol{10}$, $11, x_3, \boldsymbol{x_4})$. Considering the above facts, we observe that $(x_3, x_4) \in \{(6,7)$, $(7,6)$, $(8, 9)$, $(9, 8)\}$. 

\smallskip

\noindent{\bf Claim 3.} $(x_3, x_4) = (7,6)$ or $(8, 9)$. 

If $(x_3,x_4)=(6,7)$, then ${\rm lk}(6) = C_7(7, \boldsymbol{0}, 5, 10, 11, 2, \boldsymbol{3})$. This implies ${\rm lk}(7) = C_7(6, \boldsymbol{2}, 3, x_5, 8, 0$, $\boldsymbol{5})$. Now observe that $x_5$ has no value in $V(M)$. Similarly if $(x_3, x_4) = (9, 8)$, then ${\rm lk}(9) = C_7(8, \boldsymbol{0}, 1, 6, 11, 2, \boldsymbol{3})$. This implies ${\rm lk}(8) = C_7(9, \boldsymbol{1}, 0, 7, x_5, 3, \boldsymbol{2})$. Now again we see that $x_5$ has no value in $V(M)$. This proves the claim. 

\smallskip

\noindent{\bf Case 1.} If $(x_3, x_4) = (8,9)$, then completing successively, we get 
${\rm lk}(8) = C_7(9, \boldsymbol{1}, 0, 7, 11, 2, \boldsymbol{3})$, 
${\rm lk}(9) = C_7(8, \boldsymbol{0}, 1, 10, 6, 3, \boldsymbol{2})$, 
${\rm lk}(6) = C_7(7, \boldsymbol{0}, 5, 10, 9, 3, \boldsymbol{4})$, 
${\rm lk}(7) = C_7(6, \boldsymbol{3}, 4, 11, 8, 0, \boldsymbol{5})$, 
${\rm lk}(10) = C_7(11, \boldsymbol{2}, 1, 9, 6, 5, \boldsymbol{4})$, 
${\rm lk}(11) = C_7(10, \boldsymbol{1}, 2, 8, 7, 4, \boldsymbol{5})$, 
${\rm lk}(3) =C_9(2, \boldsymbol{1}, \boldsymbol{0}, \boldsymbol{5}, 4, \boldsymbol{7}, 6, 9, \boldsymbol{8})$, 
${\rm lk}(4) =C_9(3$, $\boldsymbol{2}, \boldsymbol{1}, \boldsymbol{0}, 5, \boldsymbol{10}, 11, 7, \boldsymbol{6})$. Then $M \cong \mathcal{B}_1(K)$ by the identity map.  

\smallskip

\noindent{\bf Case 2.} If $(x_3, x_4) = (7,6)$, then ${\rm lk}(7) = C_7(6, \boldsymbol{3}, 2, 11, 8, 0, \boldsymbol{5})$. Observe that $(x_5, x_6) \in \{(9, 10)$, $(10, 9)\}$. If $(x_5, x_6) = (9, 10)$, then completing successively, we get 
${\rm lk}(6) = C_7(7, \boldsymbol{0}, 5, 9, 10, 3, \boldsymbol{2})$, 
${\rm lk}(9) = C_7(8, \boldsymbol{0}, 1, 10, 6, 5, \boldsymbol{4})$,  
${\rm lk}(8) = C_7(9, \boldsymbol{1}, 0, 7, 11, 4, \boldsymbol{5})$,
${\rm lk}(10) = C_7(11, \boldsymbol{2}, 1, 9, 6, 3, \boldsymbol{4})$, 
${\rm lk}(11) = C_7(10, \boldsymbol{1}, 2, 7, 8, 4, \boldsymbol{3})$, 
${\rm lk}(3) =C_9(2, \boldsymbol{1}, \boldsymbol{0}, \boldsymbol{5}, 4, \boldsymbol{11}, 10, 6, \boldsymbol{7})$, 
${\rm lk}(4) =C_9(3, \boldsymbol{2}, \boldsymbol{1}, \boldsymbol{0}, 5, \boldsymbol{9}, 8, 11, \boldsymbol{10})$. Then $M \cong \mathcal{B}_1(K)$ by the map $(0,2,4)(1,3,5)(6,10)(7,11)$. 

On the other hand, if $(x_5, x_6) = (10, 9)$, then completing successively, we get 
${\rm lk}(6) = C_7(7, \boldsymbol{0}, 5$, $10, 9, 3, \boldsymbol{2})$, 
${\rm lk}(9) = C_7(8, \boldsymbol{0}, 1, 10, 6, 3, \boldsymbol{4})$,  
${\rm lk}(8) = C_7(9, \boldsymbol{1}, 0, 7, 11, 4, \boldsymbol{3})$,
${\rm lk}(10) = C_7(11, \boldsymbol{2}, 1, 9, 6, 5, \boldsymbol{4})$, 
${\rm lk}(11) = C_7(10, \boldsymbol{1}, 2, 7, 8, 4, \boldsymbol{5})$, 
${\rm lk}(3) =C_9(2, \boldsymbol{1}, \boldsymbol{0}, \boldsymbol{5}, 4, \boldsymbol{8}, 9, 6, \boldsymbol{7})$, 
${\rm lk}(4) =C_9(3, \boldsymbol{2}, \boldsymbol{1}, \boldsymbol{0}, 5, \boldsymbol{10}, 11, 8, \boldsymbol{9})$. Then $M \cong \mathcal{B}_2(T)$ by the identity map. This proves the lemma. \hfill $\Box$

\begin{lem}  \label{l4}
	Let $M$ be a 2-semi equivelar map of the type $[3^2.4.3.4:3.4.6.4]$ on $\leq 12$ vertices. Then $M$ is isomorphic to $\mathcal{C}_1(K)$ or $\mathcal{C}_2(T)$ given in example Section \ref{s3}. 
\end{lem}

\noindent{\bf Proof:}  Without loss of generality let $f_{seq}(0)=(3.4.6.4)$ and ${\rm lk}(0) = C_9(1, \boldsymbol{2}, \boldsymbol{3}, \boldsymbol{4}, 5, \boldsymbol{6}, 7, 8, \boldsymbol{9})$. Since, the vertices 1, 2, 3, 4, and 5 lie on the hexagonal face $[0,1,2,3,4,5]$, $f_{seq}(1)=f_{seq}(2)=f_{seq}(3)=f_{seq}(4)=f_{seq}(5)=(3.4.6.4)$. This implies ${\rm lk}(1) = C_9(2, \boldsymbol{3}, \boldsymbol{4}, \boldsymbol{5}, 0, \boldsymbol{8}, 9, x_1, \boldsymbol{x_2})$, where $x_1, x_2 \in V(M)$. Observe that, $x_1 \in \{6,7,10\}$. If $x_1=6$, then $x_2 \in \{7, 10\}$, but for $x_2=7$ we see that the edge 67 appears in two distinct quadrangular faces $[0,5,6,7]$ and $[1,2,7,6]$, which is not possible. If $x_1=7$, then $x_2=6$ and we see again that two distinct quadrangular faces share more than one vertices. If $x_1=10$, then $x_2 \in \{6,7,11\}$. This gives $(x_1,x_2) \in \{(6,10), (10,6), (10, 7),(10,11)\}$. Since $(6,10) \cong (10,7)$ by the map $(0,1)(2,5)(3,4)(6,7,10)(8,9)$, we need not consider the case $(6,10)$. 

\smallskip

\noindent{\bf Claim 1.} $(x_1, x_2) = (10, 6)$ or $(10,7)$.

For $(x_1,x_2)=(10,11)$, we get ${\rm lk}(1) = C_9(2, \boldsymbol{3}, \boldsymbol{4}, \boldsymbol{5}, 0, \boldsymbol{8}, 9, 10, \boldsymbol{11})$. This implies ${\rm lk}(2) = C_9(3, \boldsymbol{4}$, $\boldsymbol{5}, \boldsymbol{0}, 1, \boldsymbol{10}, 11, x_3, \boldsymbol{x_4})$, where $(x_3,x_4) \in \{(6,8), (6,9), (7,8), (8,7)\}$. If $(x_3,x_4)=(6,8)$, then ${\rm lk}(6) = C_7(2,11,7, \boldsymbol{0}, 5, 8, \boldsymbol{3})$, which implies $\deg(8)>5$. If $(x_3,x_4)=(6,9)$, then ${\rm lk}(6) = C_7(2,11,7, \boldsymbol{0}, 5, 9, \boldsymbol{3})$ or ${\rm lk}(6) = C_7(2,11,5, \boldsymbol{0}, 7, 9, \boldsymbol{3})$. But for both the cases, we see $\deg(9)>5$. If $(x_3,x_4)=(8,7)$, then ${\rm lk}(7) = C_7(3,x_5,6, \boldsymbol{5}, 0, 8, \boldsymbol{2})$. Now observe that $x_5$ has no value in $V(M)$. If $(x_3,x_4)=(7,8)$, then ${\rm lk}(7) = C_7(2,11,6, \boldsymbol{5}, 0, 8, \boldsymbol{3})$ and ${\rm lk}(3) = C_9(4, \boldsymbol{5}, \boldsymbol{0}, \boldsymbol{1}, 2, \boldsymbol{7}, 8, x_5, \boldsymbol{x_6})$, where $(x_5,x_6) \in \{(6,11), (10,9)\}$. In case $(x_5,x_6)=(6,11)$ (or $(10,9)$), considering ${\rm lk}(8)$, we see $\deg(6)>5$ (resp. $C_4(3,4,9,8) \subseteq {\rm lk}(10)$). This proves the claim.


\noindent{\bf Case 1.} $(x_1, x_2)=(10,6)$, i.e., ${\rm lk}(1) = C_9(2, \boldsymbol{3}, \boldsymbol{4}, \boldsymbol{5}, 0, \boldsymbol{8}, 9, 10, \boldsymbol{6})$ then ${\rm lk}(2) = C_9(3, \boldsymbol{4}, \boldsymbol{5}, \boldsymbol{0}, 1$, $\boldsymbol{10}, 6, x_3, \boldsymbol{x_4})$, where $(x_3,x_4) \in \{(7,9), (8,7), (11,7), (11,8), (11,9)\}$. If $(x_3,x_4)=(7,9)$ then considering ${\rm lk}(7)$ we get $C_4(0,1,9,7) \subseteq {\rm lk}(8)$. If $(x_3, x_4)=(8,7)$, then considering ${\rm lk}(6)$, we get three quadrangular faces $[6,5,0,7]$, $[6,2,1,10]$ and $[6,9,x_5,x_6]$ at $6$, which is not allowed. If $(x_3,x_4)=(11,7)$ or $(11,9)$ then ${\rm lk}(6) = C_7(5,11,2, \boldsymbol{1}, 10, 7, \boldsymbol{0})$ or ${\rm lk}(6) = C_7(7,11,2, \boldsymbol{1}, 10, 5, \boldsymbol{0})$. In the first case $\deg(7)>5$ and in the second case we get triangular face $[6,7,11]$ which is not possible, see ${\rm lk}(2)$. Now we search for the remaining cases. 

If $(x_3,x_4)=(11,8)$, then ${\rm lk}(2) = C_9(3, \boldsymbol{4}, \boldsymbol{5}, \boldsymbol{0}, 1, \boldsymbol{10}, 6, 11, \boldsymbol{8})$ and ${\rm lk}(3) = C_9(4, \boldsymbol{5}, \boldsymbol{0}, \boldsymbol{1}, 2, \boldsymbol{11}$, $8, x_5, \boldsymbol{x_6})$. Here, $(x_5, x_6) \in \{(7,9), (7,10), (9,7)\}$. In case $(x_5,x_6)=(7,9)$, considering ${\rm lk}(7)$, we get $\deg(6)>5$. In case $(x_5,x_6)=(9,7)$, considering ${\rm lk}(9)$, we get $\deg(7)>5$. If $(x_5,x_6)=(7,10)$, then ${\rm lk}(7) = C_7(0,8,3, \boldsymbol{4}, 10, 6, \boldsymbol{5})$, ${\rm lk}(6) = C_7(2,11,5, \boldsymbol{0}, 7, 10, \boldsymbol{1})$. Now completing successively, we get ${\rm lk}(5) = C_9(0, \boldsymbol{1}, \boldsymbol{2}, \boldsymbol{3}, 4, \boldsymbol{9}, 11, 6, \boldsymbol{7})$, ${\rm lk}(11) = C_7(2,6,5, \boldsymbol{4}, 9, 8, \boldsymbol{3})$, ${\rm lk}(8) = C_7(0,7,3, \boldsymbol{2}, 11, 9, \boldsymbol{1})$, ${\rm lk}(4) = C_9(5, \boldsymbol{0}, \boldsymbol{1}, \boldsymbol{2}, 3, \boldsymbol{7}, 10, 9, \boldsymbol{11})$, ${\rm lk}(9) = C_7(1,10,4, \boldsymbol{5}, 11, 8, \boldsymbol{0})$, and ${\rm lk}(10) = C_7(1,9,4, \boldsymbol{3}, 7, 6, \boldsymbol{2})$. Then $M \cong \mathcal{C}_2(T)$ by the identity map.

\smallskip 

If $(x_3,x_4)=(11,9)$, then ${\rm lk}(2) = C_9(3, \boldsymbol{4}, \boldsymbol{5}, \boldsymbol{0}, 1, \boldsymbol{10}, 6, 11, \boldsymbol{9})$ and ${\rm lk}(3) = C_9(4, \boldsymbol{5}, \boldsymbol{0}, \boldsymbol{1}, 2, \boldsymbol{11}$, $9, x_5, \boldsymbol{x_6})$. Here, $(x_5, x_6) \in \{(8,10), (10,7), (10,8)\}$. If $(x_5,x_6)=(8,10)$, we get ${\rm lk}(8) = C_7(0,7,10$, $\boldsymbol{4}, 3, 9, \boldsymbol{1})$, which implies ${\rm lk}(4)$ can not be completed. If $(x_5,x_6)=(10,8)$, then considering ${\rm lk}(9)$ and ${\rm lk}(10)$, we get $\deg(6)>5$. If $(x_5,x_6)=(10,7)$, then ${\rm lk}(4) = C_9(5, \boldsymbol{0}, \boldsymbol{1}, \boldsymbol{2}, 3, \boldsymbol{10}, 7, 8, \boldsymbol{11})$, ${\rm lk}(8) = C_7(0,7,4, \boldsymbol{5}, 11, 9, \boldsymbol{1})$, ${\rm lk}(7) = C_7(0,8,4, \boldsymbol{3}, 10, 6, \boldsymbol{5})$, ${\rm lk}(6) = C_7(2,11,5, \boldsymbol{0}, 7, 10, \boldsymbol{1})$. Now completing successively, we get ${\rm lk}(11) = C_7(2,6,5, \boldsymbol{4}, 8, 9, \boldsymbol{3})$, ${\rm lk}(5) = C_9(0, \boldsymbol{1}, \boldsymbol{2}, \boldsymbol{3}, 4, \boldsymbol{8}, 11, 6, \boldsymbol{7})$, ${\rm lk}(9) = C_7(1,10,3, \boldsymbol{2}, 11, 8, \boldsymbol{0})$. Then $M \cong \mathcal{C}_1(K)$ by the identity map. 
 

\smallskip

\noindent{\bf Case 2.} $(x_1,x_2)=(10,7)$, i.e., ${\rm lk}(1) = C_9(2, \boldsymbol{3}, \boldsymbol{4}, \boldsymbol{5}, 0, \boldsymbol{8}, 9, 10, \boldsymbol{7})$. Then ${\rm lk}(7) = C_7(0,8,10, \boldsymbol{1}$, $2, 6, \boldsymbol{5})$ or ${\rm lk}(7) = C_7(0,8,2, \boldsymbol{1}, 10, 6, \boldsymbol{5})$. In the first case of ${\rm lk}(7)$, we get ${\rm lk}(2) = C_9(3, \boldsymbol{4}, \boldsymbol{5}, \boldsymbol{0}, 1, \boldsymbol{10}, 7$, $6, \boldsymbol{x_3})$, where $x_3 \in \{9,11\}$. But for both the cases of $x_3$, ${\rm lk}(6)$ can not be completed. While for ${\rm lk}(7) = C_7(0,8,2, \boldsymbol{1}, 10, 6, \boldsymbol{5})$, we get ${\rm lk}(2) = C_9(3, \boldsymbol{4}, \boldsymbol{5}, \boldsymbol{0}, 1, \boldsymbol{10}, 7, 8, \boldsymbol{x_3})$, where $x_3 \in \{6,11\}$. If $x_3=6$, considering ${\rm lk}(2)$ and ${\rm lk}(8)$, we see $\deg(6)>5$. So, let $x_3=11$. Then ${\rm lk}(2) = C_9(3, \boldsymbol{4}, \boldsymbol{5}, \boldsymbol{0}, 1, \boldsymbol{10}, 7, 8, \boldsymbol{11})$, ${\rm lk}(8) = C_7(0,7,2, \boldsymbol{3}, 11, 9, \boldsymbol{1})$. Note that ${\rm lk}(3) = C_9(4, \boldsymbol{5}, \boldsymbol{0}, \boldsymbol{1}, 2, \boldsymbol{8},11, 6$, $\boldsymbol{10})$. Now completing successively, we get ${\rm lk}(6) = C_7(3,11,5, \boldsymbol{0}, 7, 10, \boldsymbol{4})$, ${\rm lk}(10) = C_7(1,9,4, \boldsymbol{3},6, 7, \boldsymbol{2})$, ${\rm lk}(9) = C_7(1,10,4, \boldsymbol{5}, 11, 8, \boldsymbol{0})$, ${\rm lk}(4) = C_9(5, \boldsymbol{0}, \boldsymbol{1}, \boldsymbol{2}, 3, \boldsymbol{6}, 10, 9, \boldsymbol{11})$, ${\rm lk}(5) = C_9(0, \boldsymbol{1}, \boldsymbol{2}, \boldsymbol{3}, 4, \boldsymbol{9}, 11, 6, \boldsymbol{7})$. Then $M \cong \mathcal{C}_1$ by the map $(1,5)(2,4)(6,9)(7,8)(10,11)$. This proves the lemma. \hfill $\Box$

\begin{lem} \label{l5}
	Let $M$ be a 2-semi equivelar map of type $[3^6:3^2.4.3.4]$ on $\leq 12$. Then $M$ is isomorphic to $\mathcal{D}_1(K)$ given in example Section \ref{s3}.
\end{lem}

\noindent{\bf Proof:} Without loss of generality, let $0$ be a critical vertex in $M$ with the face-sequence $(3^6)$ and ${\rm lk}(0) = C_6(1,2,3,4,5,6)$. Since, each vertex in the ${\rm lk}(0)$ has two triangular faces, without loss of generality, let $f_{seq}(1) = (3^2.4.3.4)$. This gives ${\rm lk}(1)=C_7(2,0,6,\boldsymbol{x_4}, x_3, x_2, \boldsymbol{x_1})$. It is easy to see that $(x_1, x_2, x_3, x_4) \in \{(5,4,7,8), (8,7,4,3), (9,7,8,10)\}$.

\smallskip

\noindent{\bf Case 1:} $(x_1, x_2, x_3, x_4) = (5,4,7,8)$, i.e., ${\rm lk}(1)=C_7(2,0,6,\boldsymbol{8}, 7, 4, \boldsymbol{5})$. This implies ${\rm lk}(4)=C_7(3,0,5,\boldsymbol{2}, 1, 7, \boldsymbol{9})$ and ${\rm lk}(3)=C_7(2,0,4,\boldsymbol{7}, 9, x_5, \boldsymbol{x_6})$, where $(x_5, x_6) \in \{(10, 8), (10, 11), (8, 10)\}$. If $(x_5, x_6) = (10, 8)$, then ${\rm lk}(2)=C_7(1,0,3,\boldsymbol{10}, 8, 5, \boldsymbol{4})$. Now considering ${\rm lk}(5)$ and ${\rm lk}(1)$, we see that the set $\{6, 8\}$ forms both an edge and a non-edge. If $(x_5, x_6)=(10,11)$, then ${\rm lk}(2)=C_7(1,0,3,\boldsymbol{10}, 11, 5, \boldsymbol{4})$, ${\rm lk}(5)=C_7(4,0,6,\boldsymbol{9}, 11, 2, \boldsymbol{1})$, ${\rm lk}(6)=C_7(1,0,5,\boldsymbol{11}, 9, 8, \boldsymbol{7})$, and ${\rm lk}(7)=C_7(8,x_7,9,\boldsymbol{3}, 4, 1, \boldsymbol{6})$. Now observe that $x$ has no value in $V(M)$ so that ${\rm lk}(7)$ can be completed. If $(x_5,x_6)=(8,10)$, then ${\rm lk}(2)=C_7(1,0,3,\boldsymbol{8}, 10, 5, \boldsymbol{4})$, ${\rm lk}(5)=C_7(6,0,4,\boldsymbol{1}, 2$, $10, \boldsymbol{x_7})$. It is easy to see that $x_7=9$, completing successively, we get ${\rm lk}(5)=C_7(4,0,6,\boldsymbol{9}, 10, 2,\boldsymbol{1})$, ${\rm lk}(6)=C_7(1,0,5,\boldsymbol{10}, 9, 8, \boldsymbol{7})$, ${\rm lk}(9)=C_7(3,8,6,\boldsymbol{5}, 10, 7, \boldsymbol{4})$, ${\rm lk}(7)=C_7(8,10,9,\boldsymbol{3}, 4, 1, \boldsymbol{6})$, ${\rm lk}(10)$ $=C_7(8,7$, $9,\boldsymbol{6}, 5, 2, \boldsymbol{3})$ and ${\rm lk}(8)=C_7(3,9,6,\boldsymbol{1}, 7, 10, \boldsymbol{2})$. Then $M \cong \mathcal{D}_1(K)$ by the identity map. 

\smallskip

\noindent{\bf Case 2:} $(x_1, x_2, x_3, x_4)=(8,7,4,3)$, i.e., ${\rm lk}(1)=C_7(2,0,6,\boldsymbol{3}, 4, 7, \boldsymbol{8})$. Then ${\rm lk}(4)=C_7(3,0,5,\boldsymbol{9}, 7$, $1, \boldsymbol{6})$ and  ${\rm lk}(7)=C_7(8,10,9,\boldsymbol{5}, 4, 1, \boldsymbol{2})$. This implies ${\rm lk}(2)=C_7(3,0,1,\boldsymbol{7}, 8, x_5, \boldsymbol{x_6})$, where, $(x_5, x_6) \in \{(11,9)$, $(9, 10)\}$. If $(x_5,x_6)=(11,9)$, then
${\rm lk}(3)=C_7(2,0,4,\boldsymbol{1}, 6, 9, \boldsymbol{11})$, ${\rm lk}(6)=C_7(1,0,5,\boldsymbol{10}, 9, 3$, $\boldsymbol{4})$, and ${\rm lk}(9)=C_7(10,x_7,11,\boldsymbol{2}, 3, 6, \boldsymbol{5})$, but observe that $x_7$ has no value in $V(M)$. On the other hand, if $(x_5, x_6)=(9, 10)$, completing successively, we get ${\rm lk}(2)=C_7(1,0,3,\boldsymbol{10}, 9, 8, \boldsymbol{7})$, ${\rm lk}(3)=C_7(2,0,4,\boldsymbol{1}, 6, 10, \boldsymbol{9})$, ${\rm lk}(10)=C_7(8,7,9,\boldsymbol{2}, 3, 6, \boldsymbol{5})$, ${\rm lk}(5)=C_7(4,0,6,\boldsymbol{10}, 8, 9, \boldsymbol{7})$, ${\rm lk}(8)=C_7(2,9,5,\boldsymbol{6}, 10, 7, \boldsymbol{1})$, ${\rm lk}(9)=C_7(2,8,5,\boldsymbol{4}, 7, 10, \boldsymbol{3})$. Then $M \cong \mathcal{D}_1(K)$ by the map $(2,6)(3,5)$.

\smallskip

\noindent{\bf Case 3:} $(x_1, x_2, x_3, x_4)=(9,7,8,10)$, i.e., ${\rm lk}(1)=C_7(2,0,6,\boldsymbol{10}, 8, 7, \boldsymbol{9})$. This implies ${\rm lk}(2)=C_7(3,0,1,\boldsymbol{7}, 9, x_5, \boldsymbol{x_6})$, where $(x_5, x_6) \in \{(10, 11), (11,8), (11,10), (5,6)\}$. If $(x_5,x_6)=(10, 11)$, then ${\rm lk}(10)=C_7(2,9, 8,\boldsymbol{1}, 6, 11, \boldsymbol{3})$ or ${\rm lk}(10)=C_7(2,9, 6,\boldsymbol{1}, 8, 11, \boldsymbol{3})$. In the first case, we get ${\rm lk}(6)=C_7(1,0, 5,\boldsymbol{7}, 11, 10, \boldsymbol{8})$ and ${\rm lk}(11)=C_7(3,x_7, 7,\boldsymbol{5}, 6, 10, \boldsymbol{2})$, but observe that $x_7$ has no value in $V(M)$. On the other hand when, ${\rm lk}(10)=C_7(2,9, 6,\boldsymbol{1}, 8, 11, \boldsymbol{3})$, then considering ${\rm lk}(10)$, ${\rm lk}(6)$, and ${\rm lk}(10)$ successively, we see that ${\rm lk}(11)$ can not be completed. If $(x_5, x_6)=(11,8)$, then ${\rm lk}(8)=C_7(1,7, 3,\boldsymbol{2}, 11, 10, \boldsymbol{6})$ or ${\rm lk}(8)=C_7(1,7,11,\boldsymbol{2}, 3, 10, \boldsymbol{6})$. Now, as in previous case, we see that link of all vertices can not be completed. If $(x_5,x_6)=(11,10)$, then ${\rm lk}(2)=C_7(1,0,3,\boldsymbol{10}, 11, 9, \boldsymbol{7})$. This implies ${\rm lk}(3)=C_7(4,0, 2,\boldsymbol{11}, 10, x_7, \boldsymbol{x_8})$, where $(x_7,x_8) \in \{(7,8), (8,7)\}$. But for both the cases, ${\rm lk}(7)$ can not be completed. If $(x_5,x_6)=(5,6)$, then ${\rm lk}(2)=C_7(1,0,3,\boldsymbol{6}, 5, 9, \boldsymbol{7})$, and ${\rm lk}(5)=C_7(4,0,6,\boldsymbol{3}, 2, 9, \boldsymbol{x_7})$, where $x_7 \in \{8, 10\}$. If $x_7=10$, then considering ${\rm lk}(5)$ and  ${\rm lk}(6)$ we see that ${\rm lk}(3)$ can not be completed. So, $x_7=8$. Completing successively, we get ${\rm lk}(5)=C_7(4,0,6,\boldsymbol{3}, 2, 9, \boldsymbol{8})$, ${\rm lk}(6)=C_7(1,0,5,\boldsymbol{2}, 3, 10, \boldsymbol{8})$, ${\rm lk}(3)=C_7(2,0,4,\boldsymbol{7}, 10, 5, \boldsymbol{6})$, ${\rm lk}(4)=C_7(3,0,5,\boldsymbol{9}, 8, 7, \boldsymbol{10})$, ${\rm lk}(7)=C_7(1,8,4,\boldsymbol{3}, 10, 9, \boldsymbol{2})$, ${\rm lk}(8)=C_7(1,7,4,\boldsymbol{5}, 9, 10, \boldsymbol{6})$, ${\rm lk}(9)=C_7(7,10,8,\boldsymbol{4}, 5, 2, \boldsymbol{1})$. Then $M \cong \mathcal{D}_1(K)$ by the map $(1,6,5,4,3$, $2)(7,8,9)$. This proves the lemma. \hfill$\Box$

\begin{lem}  \label{l6}
Let $M$ be a 2-semi equivelar map of type $[3^6:3^3.4^2]$ on $\leq 12$ vertices. Then $M$ is isomorphic to one of $\mathcal{E}_1(K)$, $\mathcal{E}_2(K)$, $\mathcal{E}_3(T)$, $\mathcal{E}_4(T)$, $\mathcal{E}_5(K)$, $\mathcal{E}_6(T)$, $\mathcal{E}_7(K)$, $\mathcal{E}_8(T)$, $\mathcal{E}_9(K)$, $\mathcal{E}_{10}(K)$, $\mathcal{E}_{11}(T)$, $\mathcal{E}_{12}(T)$ or $\mathcal{E}_{13}(T)$ given in example Section \ref{s3}.
\end{lem}
	
\noindent{\bf Proof:} Without loss of generality, let $f_{seq}(0)= (3^3.4^2)$ and ${\rm lk}(0) = C_7(1, \boldsymbol{2}, 3, 4, 5, 6, \boldsymbol{7})$. Then the vertices 1, 2, 3, 6 and 7 have the face-sequence $(3^3.4^2)$. Therefore, ${\rm lk}(1) = C_7(0, \boldsymbol{6}, 7, x_1, x_2, 2$, $\boldsymbol{3})$, where $x_1, y_1 \in V(M)$. Then we see that $(x_1,x_2) \in \{(4,5), (4,8), (5,4), (5,8), (8,4), (8,5)$, $(8,9)\}$. Observe that the case $(4,8) \cong (8,5)$ by the map $(0,1)(2,3)(4,5,8)(6,7)$. So, we need not discuss the case $(x_1, x_2) = (8,5)$.  
	
\smallskip
	
\noindent{\bf Case 1.} $(x_1, x_2) = (4,5)$, i.e., ${\rm lk}(1) = C_7(0, \boldsymbol{6}, 7, 4, 5, 2, \boldsymbol{3})$. This implies ${\rm lk}(4) = C_6(0,3,8,7,1,5)$ and 
	${\rm lk}(3) = C_7(2, \boldsymbol{1}, 0, 4, 8, x_3, \boldsymbol{x_4})$, where $(x_3, x_4) \in \{(6,7), (9,10)\}$. 
	In case $(x_3, x_4) = (6,7)$, completing successively, we get 
	${\rm lk}(3) = C_7(2, \boldsymbol{1}, 0, 4, 8, 6, \boldsymbol{7})$, 
	${\rm lk}(7) = C_7(6, \boldsymbol{0}, 1, 4, 8, 2, \boldsymbol{3})$,
	${\rm lk}(2) = C_7(3, \boldsymbol{0}, 1, 5, 8, 7, \boldsymbol{6})$,
	${\rm lk}(5) = C_6(0,4,1,2,8,6)$,
	${\rm lk}(6) = C_7(7, \boldsymbol{1}, 0, 5, 8, 3, \boldsymbol{2})$, 
	${\rm lk}(8) = C_6(2,5,6,3$, $4,7)$. Then $M \cong \mathcal{E}_1(K)$ by the identity map. 
	
	On the other hand, when $(x_3, x_4) = (9,10)$, then 
	${\rm lk}(2) = C_7(3, \boldsymbol{0}, 1, 5, x_3, 10, \boldsymbol{9})$, where $x_3 \in \{6, 11\}$. If $x_3= 6$, we get $C_5(0,4,1,2,6) \subseteq {\rm lk}(5)$. So, $x_3=11$. Now completing successively, we get 
	${\rm lk}(5) = C_6(0,4,1,2,11,6)$, ${\rm lk}(6) = C_7(7, \boldsymbol{1}, 0, 5, 11, 9, \boldsymbol{10})$, 
	${\rm lk}(9) = C_7(10, \boldsymbol{2}, 3, 8, 11, 6, \boldsymbol{7})$, 
	${\rm lk}(8) = C_6(3,4,7,10,11,9)$, 
	${\rm lk}(10) = C_7(9, \boldsymbol{3}$, $2, 11, 8, 7, \boldsymbol{6})$, 
	${\rm lk}(7) = C_7(6, \boldsymbol{0}, 1, 4, 8, 10, \boldsymbol{9})$, and 
	${\rm lk}(11) = C_6(2,5,6,9,8,10)$. Then $M \cong \mathcal{E}_2(K)$  by the identity map.
	
\smallskip	
	
\noindent{\bf Case 2.} $(x_1, x_2) = (4,8)$.  Then ${\rm lk}(4) = C_6(0,3,7,1,8,5)$ or 
	${\rm lk}(4) = C_6(0,3,8,1,7,5)$. If ${\rm lk}(4) = C_6(0,3,7,1,8,5)$, then ${\rm lk}(3) = C_7(2, \boldsymbol{1}, 0, 4, 7, 9, \boldsymbol{10})$, ${\rm lk}(2) = C_7(3, \boldsymbol{0}, 1, 8, x_3, 10, \boldsymbol{9})$. Observe that the possible value of $x_3=5$, but then ${\rm lk}(5) = C_6(0,4,8,2,10,6)$ and we get $C_4(1,2,5,4)$. So ${\rm lk}(4) = C_6(0,3,8,1,7,5)$. Then 
	${\rm lk}(3) = C_7(2, \boldsymbol{1}, 0, 4, 8, x_3, \boldsymbol{x_4})$, where $(x_3, x_4) \in \{(6,7), (9,10)\}$.  
	If $(x_3, x_4) = (6,7)$, completing successively, we get 
	${\rm lk}(3) = C_7(2, \boldsymbol{1}, 0, 4, 8, 6, \boldsymbol{7})$,
	${\rm lk}(6) = C_7(7, \boldsymbol{1}, 0, 5, 8$, $3, \boldsymbol{2})$, 
	${\rm lk}(7) = C_7(6, \boldsymbol{0}, 1, 4, 5, 2, \boldsymbol{3})$,
	${\rm lk}(2) = C_7(3, \boldsymbol{0}, 1, 8, 5, 7, \boldsymbol{6})$,
	${\rm lk}(5) = C_6(0,4,7,2,8,6)$, and ${\rm lk}(8) = C_6(1,2,5,6,3,4)$. Then $M \cong \mathcal{E}_3(T)$  by the identity map.

	On the other hand when $(x_3, x_4) = (9,10)$, then 
	${\rm lk}(2) = C_7(3, \boldsymbol{0}, 1, 8, 11, 10, \boldsymbol{9})$, completing successively, we get ${\rm lk}(8) = C_6(1,2,11,9,3,4)$,
	${\rm lk}(9) = C_7(10, \boldsymbol{2}, 3, 8, 11, 6, \boldsymbol{7})$,
	${\rm lk}(6) = C_7(7, \boldsymbol{1}, 0, 5, 11$, $9, \boldsymbol{10})$,
	${\rm lk}(7) = C_7(6, \boldsymbol{0}, 1, 4, 5, 10, \boldsymbol{9})$,
	${\rm lk}(5) = C_6(0,4,7,10,11,6)$, 
	${\rm lk}(10) = C_7(9, \boldsymbol{3}, 2$, $11, 5, 7, \boldsymbol{6})$, and ${\rm lk}(11) = C_6(2,8,9,6,5,10)$. Then $M \cong \mathcal{E}_4(T)$  by the identity map.
	
\smallskip
	
\noindent{\bf Case 3.} 	$(x_1, x_2) = (5,4)$.  Then ${\rm lk}(5) = C_6(0,4,1,7,8,6)$. This implies ${\rm lk}(4) = C_6(0,3,8,2,1,5)$ or 
	${\rm lk}(4) = C_6(0,3,9,2,1,5)$. In the first case, when ${\rm lk}(4) = C_6(0,3,8,2,1,5)$, then ${\rm lk}(8) = C_6(2,4,3,7,5,6)$ or ${\rm lk}(8) = C_6(2,4,3,6,5,7)$. In case ${\rm lk}(8) = C_6(2,4,3,7,5,6)$, completing successively, we get 
	${\rm lk}(7) = C_7(6, \boldsymbol{0}, 1, 5, 8, 3, \boldsymbol{2})$, 
	${\rm lk}(6) = C_7(7, \boldsymbol{1}, 0, 5, 8, 2, \boldsymbol{3})$, 
	${\rm lk}(2) = C_7(3, \boldsymbol{0}, 1, 4, 8, 6, \boldsymbol{7})$, and ${\rm lk}(3) = C_7(2, \boldsymbol{1}, 0, 4, 8, 7, \boldsymbol{6})$. Then $M \cong \mathcal{E}_5(K)$  by the identity map. While, for 
	${\rm lk}(8) = C_6(2,4,3,6,5,7)$, completing successively, we get 
	${\rm lk}(7) = C_7(6, \boldsymbol{0}, 1, 5, 8, 2, \boldsymbol{3})$, 
	${\rm lk}(6) = C_7(7, \boldsymbol{1}, 0, 5, 8$, $3, \boldsymbol{2})$, 
	${\rm lk}(2) = C_7(3, \boldsymbol{0}, 1, 4, 8, 7, \boldsymbol{6})$. Then $M \cong \mathcal{E}_6(T)$ by the identity map.
	
	On the other hand, when ${\rm lk}(4) = C_6(0,3,9,2,1,5)$, then 
	${\rm lk}(2) = C_7(3, \boldsymbol{0}, 1, 4, 9, x_3, \boldsymbol{x_4})$, where we see easily that $(x_3, x_4) = (10,11)$. Then ${\rm lk}(3) = C_7(2, \boldsymbol{1}, 0, 4, 9, 11, \boldsymbol{10})$, ${\rm lk}(9) = C_6(2,4,3,11,8,10)$, and ${\rm lk}(10) = C_7(11, \boldsymbol{3}, 2, 9, 8, x_5, \boldsymbol{x_6})$, where $(x_5, x_6) \in \{(6,7), (7,6)\}$. If $(x_5, x_6) = (6,7)$, completing successively, we get ${\rm lk}(6) = C_7(7, \boldsymbol{1}, 0, 5, 8, 10, \boldsymbol{11})$, 
	${\rm lk}(7) = C_7(6, \boldsymbol{0}, 1, 5$, $8, 11, \boldsymbol{10})$, 
	${\rm lk}(8) = C_6(5,6,10,9,11,7)$. Then $M \cong \mathcal{E}_7(K)$  by the identity map.
	Also, if  $(x_5, x_6) = (7,6)$, completing successively, we get 
	${\rm lk}(7) = C_7(6, \boldsymbol{0}, 1, 5, 8, 10, \boldsymbol{11})$, 	
	${\rm lk}(6) = C_7(7, \boldsymbol{1}, 0, 5, 8, 11, \boldsymbol{10})$
	${\rm lk}(8) = C_6(5,6,11,9,10,7)$. Then $M \cong \mathcal{E}_8(T)$  by the identity map.
	
\smallskip
	
\noindent{\bf Case 4.} 	$(x_1, x_2) = (5,8)$. Then ${\rm lk}(5) = C_6(0,4,7,1,8,6)$. This implies ${\rm lk}(6) = C_7(7, \boldsymbol{1}, 0, 5, 8, x_3$, $\boldsymbol{x_4})$, where $(x_3, x_4) \in \{(2,3), (3,2), (9,10)\}$. If $(x_3, x_4) = (2,3)$, then  ${\rm lk}(7) = C_7(6, \boldsymbol{0}, 1, 5, 4, 3, \boldsymbol{2})$ and we get $C_4(0,3,7,5) \subseteq {\rm lk}(4)$. 
	If $(x_3, x_4) = (3,2)$, then ${\rm lk}(7) = C_7(6, \boldsymbol{0}, 1, 5, 4, 2, \boldsymbol{3})$ and ${\rm lk}(4) = C_6(0,3,8,2,7,5)$. Now completing successively, we get 
	${\rm lk}(8) = C_6(1,2,4,3,6,5)$, 
	${\rm lk}(6) = C_7(7, \boldsymbol{1}, 0, 5, 8, 3, \boldsymbol{2})$, 
	${\rm lk}(7) = C_7(6, \boldsymbol{0}, 1, 5, 4, 2, \boldsymbol{3})$, 
	${\rm lk}(2) = C_7(3, \boldsymbol{0}, 1, 8, 4, 7, \boldsymbol{6})$. Then $M \cong \mathcal{E}_1(K)$ by the map $(0,3)(1,2)(5,8)$. 
	
	If $(x_3, x_4) = (9,10)$, then ${\rm lk}(7) = C_7(6, \boldsymbol{0}, 1, 5, 4, 10, \boldsymbol{9})$ and ${\rm lk}(4) = C_6(0,3,11,10,7$, $5)$. Now completing successively, we get ${\rm lk}(10) = C_7(9, \boldsymbol{3}, 2, 11, 4, 7, \boldsymbol{6})$,
	${\rm lk}(3) = C_7(2, \boldsymbol{1}, 0, 4, 11, 9, \boldsymbol{10})$,
	${\rm lk}(9) = C_7(10, \boldsymbol{2}, 3, 11, 8, 6, \boldsymbol{7})$,
	${\rm lk}(2) = C_7(3, \boldsymbol{0}, 1, 8, 11, 10, \boldsymbol{9})$,
	${\rm lk}(8) = C_6(1$, $2,11,9,6,5)$, ${\rm lk}(11) = C_6(2,8,9$, $3,4,10)$. Then $M \cong \mathcal{E}_{9}(K)$  by the identity map. 
	
	
\smallskip
	
\noindent{\bf Case 5.} 	$(x_1, x_2) = (8,4)$. Then ${\rm lk}(4) = C_6(0,3,8,1,2,5)$. This implies ${\rm lk}(2) = C_7(3, \boldsymbol{0}, 1, 4, 5, x_3$, $ \boldsymbol{x_4})$, where we see easily that $(x_3, x_4) \in \{(7,6), (9,10)\}$. If $(x_3, x_4) = (7,6)$, then completing successively, we get 
	${\rm lk}(3) = C_7(2, \boldsymbol{1}, 0, 4, 8, 6, \boldsymbol{7})$,
	${\rm lk}(6) = C_7(7, \boldsymbol{1}, 0, 5, 8, 3, \boldsymbol{2})$, 
	${\rm lk}(7) = C_7(6, \boldsymbol{0}, 1, 8, 5, 2, \boldsymbol{3})$,
	${\rm lk}(5) = C_6(0,4,2,7,8,6)$. 
	Then $M \cong \mathcal{E}_{1}(K)$ by the map$(0,6)(1,7)(4,8)$.


	If $(x_3, x_4) = (9,10)$, then ${\rm lk}(3) = C_7(2, \boldsymbol{1}, 0, 4, 8, 10, \boldsymbol{9})$ and ${\rm lk}(8) = C_6(1,4,3,10,11,7)$. Now completing successively, we get 
	${\rm lk}(7) = C_7(6, \boldsymbol{0}, 1, 8, 11, 9, \boldsymbol{10})$,
	${\rm lk}(9) = C_7(10, \boldsymbol{3}, 2, 5, 11, 7, \boldsymbol{6})$,
	${\rm lk}(10) = C_7(9, \boldsymbol{2}, 3, 8, 11$, $6, \boldsymbol{7})$,
	${\rm lk}(5) = C_6(0,4,2,9,11,6)$,
	${\rm lk}(6) = C_7(7, \boldsymbol{1}, 0$, $5, 11, 10, \boldsymbol{9})$. Then $M \cong \mathcal{E}_{9}(K)$ by the map $(0,9,1,10)(4,11$, $5,8)$.
	. 	
\smallskip
	
\noindent{\bf Case 6.} $(x_1, x_2) = (8,9)$. Then
	${\rm lk}(1) = C_7(0, \boldsymbol{3}, 2, 9, 8, 7, \boldsymbol{6})$. This implies ${\rm lk}(2) = C_7(3, \boldsymbol{0}, 1, 9, x_3$, $x_4, \boldsymbol{x_5})$. Then it is easy to see that $(x_3, x_4, x_5) \in \{(4,5,8), (4,5,10), (4,6,7), (4,7,6), (4,10,5), (4,10$, $8), (4,10,11), (5,6,7), (5,7,6), (5,10,8), (5,10,11), (6,5,8), (6,5,10), (7,8,5), (7,8,10), (10,6,7), (10$, $7,6), (10,11,5), (10,11,8)\}$.
	
	\smallskip
	
	\noindent{\bf Claim 1.} $(x_3,x_4,x_5)=(4,10,11)$, $(5, 10, 11)$ or $(10,7,6)$. 
	
	If $(x_3, x_4, x_5)=(4,5,8)$ (or $(4,5,10)$), then 
	${\rm lk}(5) = C_7(8, \boldsymbol{3}, 2, 4, 0, 6, \boldsymbol{7})$ (resp. ${\rm lk}(5) = C_7(10, \boldsymbol{3}, 2$, $4, 0, 6, \boldsymbol{7})$), which implies $C_5(0,1,7,8,5) \subseteq {\rm lk}(6)$ (resp. $C_5(0,1,7,10,5) \subseteq {\rm lk}(6)$). If $(x_3, x_4, x_5)=(4,6,7)$, then considering ${\rm lk}(6)$, we get $C_3(0,4,5) \subseteq {\rm lk}(5)$. If $(x_3, x_4, x_5)=(4,7,6)$, considering ${\rm lk}(7)$, we get $\deg (4) >6$. If $(x_3, x_4, x_5)=(4,10,5)$, then considering ${\rm lk}(5)$, we see that ${\rm lk}(4)$ can not be completed. If $(x_3, x_4, x_5)=(4,10,8)$, then ${\rm lk}(4) = C_6(0,3,10,2,9,5)$ or ${\rm lk}(4) = C_6(0,3,9,2,10,5)$, but for both the cases, ${\rm lk}(9)$ can not be completed. By a similar computation we see easily that $M$ does not exists for $(x_3, x_4, x_5) \in \{(5,6,7), (5,7,6), (5, 10, 8), (7, 8, 10)\}$.
	
	If $(x_3,x_4,x_5)=(6,5,8)$, then ${\rm lk}(6) = C_7(7, \boldsymbol{1}, 0, 5, 2, 9, \boldsymbol{x_6})$, where $x_6 \in \{4, 10\}$. If $x_6=10$, then ${\rm lk}(9) = C_7(10, \boldsymbol{5}, 8, 1, 2, 6, \boldsymbol{x_7})$, and we see that ${\rm lk}(7)$ can not be completed. If $x_6=4$, then completing successively, we get 
	${\rm lk}(9) = C_7(4, \boldsymbol{5}, 8, 1, 2, 6, \boldsymbol{7})$, ${\rm lk}(4) = C_7(9, \boldsymbol{6}, 7, 3, 0, 5, \boldsymbol{8})$, 
	${\rm lk}(7) = C_7(6, \boldsymbol{0}, 1, 8, 3, 4, \boldsymbol{9})$,
	${\rm lk}(3) = C_7(2, \boldsymbol{1}, 0, 4, 7, 8, \boldsymbol{5})$,
	${\rm lk}(5) = C_7(8, \boldsymbol{3}, 2, 6, 0, 4, \boldsymbol{9})$,
	${\rm lk}(8) = C_7(5, \boldsymbol{2}, 3, 7, 1, 9, \boldsymbol{4})$. But, this gives a semi-equivelar map.  
	
	If $(x_3,x_4,x_5)=(6,5,10)$, then ${\rm lk}(6) = C_7(7, \boldsymbol{1}, 0, 5, 2, 9, \boldsymbol{x_6})$, where $x_6 \in \{4, 11\}$. If $x_6=11$, then we see that ${\rm lk}(9)$ can not be completed. If $x_6=4$, completing as in previous case, we get a semi-equivelar map. 
	
	If $(x_3,x_4,x_5)=(7,8,5)$, then ${\rm lk}(7) = C_7(6, \boldsymbol{0}, 1, 8, 2, 9, \boldsymbol{x_6})$, where $x_6=4$. Then, completing as in previous case, we get a semi-equivelar map.
	
	If $(x_3,x_4,x_5)=(7,8,10)$, then ${\rm lk}(8) = C_7(10, \boldsymbol{3}, 2, 7, 1, 9, \boldsymbol{x_6})$, where $x_6\in \{5,11\}$. If $x_6=5$, then ${\rm lk}(9) = C_7(5, \boldsymbol{6}, 7, 2, 1, 8, \boldsymbol{10})$ and we get $C_5(0,1,7,9,5) \subseteq {\rm lk}(6)$. On the other hand if $x_6=11$, then considering ${\rm lk}(9)$, ${\rm lk}(7)$, ${\rm lk}(6)$, and ${\rm lk}(10)$ successively, we get $C_5(6,7,9,8,10) \subseteq {\rm lk}(5)$. This proves the claim.
	
	\smallskip
	
	\noindent{\bf Subcase 6.1.} $(x_3,x_4,x_5)=(4,10,11)$. Completing successively, we get 
	${\rm lk}(2) = C_7(3, \boldsymbol{0}, 1, 9, 4, 10$, $\boldsymbol{11})$, ${\rm lk}(4) = C_6(0,3,9,2,10,5)$,
	${\rm lk}(3) = C_7(2, \boldsymbol{1}, 0, 4, 9, 11, \boldsymbol{10})$, ${\rm lk}(9) = C_6(1,2,4,3,11,8)$, 
	${\rm lk}(10) = C_7(11, \boldsymbol{3}, 2, 4, 5, 7, \boldsymbol{6})$, 
	${\rm lk}(11) = C_7(10, \boldsymbol{2}, 3, 9, 8, 6, \boldsymbol{7})$, 
	${\rm lk}(6) = C_7(7, \boldsymbol{1}, 0, 5, 8, 11, \boldsymbol{10})$, 
	${\rm lk}(5) = C_6(0,4$, $10,7,8,6)$, ${\rm lk}(8) = C_6(1,7,5,6,11,9)$. 
	Then $M \cong \mathcal{E}_{2}(K)$ by the map $(0,6,9,4,5,11,3)(1,7,10,2)$.
	
	\smallskip
	
	\noindent{\bf Subcase 6.2.} $(x_3,x_4,x_5)=(5,10,11)$.  Then ${\rm lk}(5) = C_6(0,4,9,2,10,6)$ or ${\rm lk}(5) = C_6(0,4,10,2,9$, $6)$. In the first case, considering ${\rm lk}(9)$, we see that ${\rm lk}(4)$ can not be completed. So, ${\rm lk}(5) = C_6(0,4,10,2,9,6)$. Completing successively, we get ${\rm lk}(9) = C_6(1,2,5,6,11,8)$, ${\rm lk}(6) = C_7(7, \boldsymbol{1}, 0, 5$, $9, 11, \boldsymbol{10})$, ${\rm lk}(10) = C_7(11, \boldsymbol{3}, 2, 5, 4, 7, \boldsymbol{6})$, ${\rm lk}(11) = C_7(10, \boldsymbol{2}$, $3, 8, 9, 6, \boldsymbol{7})$, ${\rm lk}(3) = C_7(2, \boldsymbol{1}, 0, 4, 8, 11$, $\boldsymbol{10})$, ${\rm lk}(4) = C_6(0,3,8,7,10,5)$. 
	Then $M \cong \mathcal{E}_{11}(T)$ by the identity map.
	
	\smallskip
	
	\noindent{\bf Subcase 6.3.} $(x_3,x_4,x_5)=(10,7,6)$. Then ${\rm lk}(2) = C_7(3, \boldsymbol{0}, 1, 9, 10, 7, \boldsymbol{6})$, ${\rm lk}(7) = C_7(6, \boldsymbol{0}, 1, 8, 10$, $2, \boldsymbol{3})$,  ${\rm lk}(6) = C_7(7, \boldsymbol{1}, 0, 5, 11, 3, \boldsymbol{2})$,  ${\rm lk}(3) = C_7(2, \boldsymbol{1}, 0, 4, 11, 6, \boldsymbol{7})$. This implies ${\rm lk}(10) = C_6(8,7,2,9$, $x_6,x_7)$ or ${\rm lk}(10) = C_7(x_7, \boldsymbol{x_8}, 8, 7, 2, 9, \boldsymbol{x_6})$. 
	
	\noindent{\bf Subcase 6.3.1.} ${\rm lk}(10) = C_6(8,7,2,9,x_6,x_7)$. Then $(x_6, x_7) \in \{(4,5), (5,4), (11,4), (11,5)\}$. If $(x_6,x_7)=(4,5)$, completing successively, we get ${\rm lk}(4) = C_6(0,3,11,9,10,5)$, 
	${\rm lk}(5) = C_6(0,4,10,8$, $11,6)$, ${\rm lk}(8) = C_6(1,7,10,5,11,9)$, 
	${\rm lk}(9) = C_6(1,2,10,4,11,8)$, and ${\rm lk}(11) = C_6(3,4,9$, $8,5,6)$. 
	Then $M \cong \mathcal{E}_{10}(K)$ by the identity map.

	If $(x_6,x_7)=(5,4)$, completing successively, we get ${\rm lk}(4) = C_6(0,3,11,8,10$, $5)$, ${\rm lk}(5) = C_6(0,4$, $10,9,11,6)$, ${\rm lk}(9) = C_6(1,2,10,5,11,8)$, ${\rm lk}(8) = C_6(1,7,10,4,11,9)$, 
	and ${\rm lk}(11) = C_6(3,4,8,9,5,6)$. Then $M \cong \mathcal{E}_{12}(T)$ by the identity map. 
	
	If $(x_6,x_7)=(11,4)$, completing successively, we get ${\rm lk}(4) = C_6(0,3,11,10,8,5)$, ${\rm lk}(8) = C_6(1,7,10,4,5,9)$, ${\rm lk}(5) = C_6(0,4,8,9,11,6)$, ${\rm lk}(9) = C_6(1,2,10,11,5,8)$,
	and ${\rm lk}(11) = C_6(3,4,10$, $9,5,6)$. Then $M \cong \mathcal{E}_{10}(K)$ by the map $(0,6)(1,7)(4,11)(9,10)$.

	If $(x_6,x_7)=(11,5)$, completing successively, we get ${\rm lk}(5) = C_6(0,4,8,10$, $11,6)$, ${\rm lk}(8) = C_6(1,7,10,5,4,9)$, ${\rm lk}(4) = C_6(0,3,11,9,8,5)$, ${\rm lk}(9) = C_6(1,2,10,11,4,8)$, and ${\rm lk}(11) = C_6(3,4,9$, $10,5,6)$. Then $M \cong \mathcal{E}_{13}(T)$ by the identity map. 
	
	\noindent{\bf Subcase 6.3.2.} On the other hand when ${\rm lk}(10) = C_7(x_7, \boldsymbol{x_8}, 8, 7, 2, 9, \boldsymbol{x_6})$, then it is easy to see that $(x_6,x_7,x_8) \in \{(5,4,11), (11,4,5), (4,5,11), (11,5,4), (4,11,5), (5,11,4)\}$. But for all these cases of $(x_6,x_7,x_8)$ we see at least one quadrangular face incident at each vertex $i$, for $0 \leq i \leq 11$, which means there is no vertex in $M$ with face-sequence $(3^6)$. Thus for these cases, $M$ does not exists. This proves the lemma. \hfill$\Box$

	\begin{lem} \label{l7}
		Let $M$ be a 2-semi equivelar map of type $[3^3.4^2:4^4]$ on $\leq 12$ vertices. Then $M$ is isomorphic to one of $\mathcal{F}_1(K)$, $\mathcal{F}_2(T)$, $\mathcal{F}_3(T)$, $\mathcal{F}_4(K)$, $\mathcal{F}_5(T)$, $\mathcal{F}_6(T)$, $\mathcal{F}_7(T)$, $\mathcal{F}_8(K)$ or $\mathcal{F}_9(T)$ given in example Section \ref{s3}.
	\end{lem}
	
	{\bf Proof:} Assume that $f_{seq}(0)= (4^4)$ and ${\rm lk}(0) = C_8(1, \boldsymbol{2}, 3, \boldsymbol{4}, 5, \boldsymbol{6}, 7, \boldsymbol{8})$ in $M$.
	Then we have two cases, either $f_{seq}(1) = (3^3.4^2)$ or $f_{seq}(1) = (4^4)$. 
	
	\noindent{\bf Case 1:} $f_{seq}(1) = (3^3.4^2)$. Then ${\rm lk}(1) = C_7(0, \boldsymbol{7}, 8, x_1, x_2, 2, \boldsymbol{3})$. It is easy to see that $(x_1, x_2) \in \{(4,5), (4,6), (4, 9), (5,6), (5,4), (6,4), (6,5), (6,9), (9,4), (9$, $6), (9, 10)\}$. Here $(4,6) \cong (6,4)$ by the map $(2,8)(3,7)(4,6)$. So, we need not consider the case $(6, 4)$. 
	
	\smallskip	
	
\noindent{\bf Subcase 1.1.}	$(x_1, x_2)=(4,5)$. Then ${\rm lk}(1) = C_7(0, \boldsymbol{7}, 8, 4,5, 2, \boldsymbol{3})$, ${\rm lk}(5) = C_7(0, \boldsymbol{3}, 4, 1, 2, 6, \boldsymbol{7})$. This implies ${\rm lk}(2) = C_7(3, \boldsymbol{0}, 1, 5, 6, x_3, \boldsymbol{x_4})$, where $(x_3, x_4) \in \{(9, 4), (9, 10)\}$. If $(x_3, x_4)=(9,4)$, then considering ${\rm lk}(4)$, we get $C_6(0,1,2,9,4,5) \subseteq {\rm lk}(3)$. If $(x_3, x_4)=(9,10)$, then it is easy to see that 
${\rm lk}(3) = C_8(0, \boldsymbol{1}, 2, \boldsymbol{9}, 10, \boldsymbol{11}, 4, \boldsymbol{5})$, ${\rm lk}(4) = C_7(3, \boldsymbol{0}, 5, 1, 8, 11, \boldsymbol{10})$ and ${\rm lk}(8) = C_7(7, \boldsymbol{0}, 1, 4, 11, x_5, \boldsymbol{x_6})$, where $(x_5, x_6)=(9, 10)$. Completing successively, we get ${\rm lk}(8) = C_7(7, \boldsymbol{0}, 1, 4, 11, 9, \boldsymbol{10})$, ${\rm lk}(9) = C_7(10, \boldsymbol{3}, 2, 6, 11, 8, \boldsymbol{7})$, ${\rm lk}(6) = C_7(7, \boldsymbol{0}, 5, 2, 9, 11, \boldsymbol{10})$, ${\rm lk}(11) = C_7(10, \boldsymbol{3}, 4, 8, 9, 6, \boldsymbol{7})$,
${\rm lk}(7) = C_8(0$, $\boldsymbol{1}, 8, \boldsymbol{9}, 10, \boldsymbol{11}, 6, \boldsymbol{5})$.  Then $M \cong \mathcal{F}_1(K)$ by the identity map.  
	
	\smallskip

\noindent{\bf Subcase 1.2.}	$(x_1, x_2)=(4,6)$.  Then ${\rm lk}(1) = C_7(0, \boldsymbol{7}, 8, 4,6, 2, \boldsymbol{3})$. This implies ${\rm lk}(4) = C_7(3, \boldsymbol{0}, 5, 8, 1$, $6, \boldsymbol{x_3})$, where $x_3=7$. Completing successively, we get ${\rm lk}(6) = C_7(7, \boldsymbol{0}, 5, 2, 1, 4, \boldsymbol{3})$, ${\rm lk}(2) = C_7(3, \boldsymbol{0}, 1, 6, 5, 8, \boldsymbol{7})$, ${\rm lk}(8) = C_7(7, \boldsymbol{0}, 1, 4, 5, 2, \boldsymbol{3})$, ${\rm lk}(7) = C_8(0, \boldsymbol{1}, 8, \boldsymbol{2}, 3, \boldsymbol{4}, 6, \boldsymbol{5})$, ${\rm lk}(5) = C_7(0, \boldsymbol{3}, 4, 8$, $2, 6, \boldsymbol{7})$, ${\rm lk}(3) = C_8(0, \boldsymbol{1}, 2, \boldsymbol{8}, 7, \boldsymbol{6}, 4, \boldsymbol{5})$. Then $M \cong \mathcal{F}_2(T)$ by the identity map.
	
\smallskip	

\noindent{\bf Subcase 1.3.}	$(x_1, x_2) = (4,9)$. Then ${\rm lk}(4) = C_7(3, \boldsymbol{0}, 5, 9, 1, 8, \boldsymbol{x_3})$ or ${\rm lk}(4) = C_7(5, \boldsymbol{0}, 3, 9, 1, 8, \boldsymbol{x_3})$ or ${\rm lk}(4) = C_7(5, \boldsymbol{0}, 3, 8, 1, 9, \boldsymbol{x_3})$, or ${\rm lk}(4) = C_7(3, \boldsymbol{0}, 5, 8, 1, 9, \boldsymbol{x_3})$, where $x_3 \in V(M)$. If ${\rm lk}(4) = C_7(3, \boldsymbol{0}, 5, 9, 1, 8, \boldsymbol{x_3})$ (or ${\rm lk}(4) = C_7(5, \boldsymbol{0}, 3, 9, 1, 8, \boldsymbol{x_3})$), then $x_3=7$ and we get $C_5(0,1,4,3,7) \subseteq {\rm lk}(8)$ (resp. $C_5(0,1,4,5,7) \subseteq {\rm lk}(8)$). If ${\rm lk}(4) = C_7(5, \boldsymbol{0}, 3, 8, 1, 9, \boldsymbol{x_3})$, then $x_3=10$. This implies ${\rm lk}(5) = C_8(0, \boldsymbol{3}, 4, \boldsymbol{9}, 10, \boldsymbol{11}, 6, \boldsymbol{7})$, but then ${\rm lk}(9)$ can not be completed.  If ${\rm lk}(4) = C_7(3, \boldsymbol{0}, 5, 8, 1, 9, \boldsymbol{x_3})$, then ${\rm lk}(5) = C_7(0, \boldsymbol{3}, 4, 8, 10, 6, \boldsymbol{7})$ and ${\rm lk}(8) = C_7(7, \boldsymbol{0}, 1, 4, 5, 10, \boldsymbol{x_4})$, where $x_4 \in \{9, 11\}$. In case $x_4=9$, ${\rm lk}(7)$ can not be completed. So, $x_4=11$. Completing successively, we get ${\rm lk}(8) = C_7(7, \boldsymbol{0}, 1, 4, 5, 10, \boldsymbol{11})$, ${\rm lk}(7) = C_8(0, \boldsymbol{1}, 8, \boldsymbol{10}, 11, \boldsymbol{9}, 6, \boldsymbol{5})$, ${\rm lk}(6) = C_7(7, \boldsymbol{0}, 5, 10, 2, 9, \boldsymbol{11})$, ${\rm lk}(2) = C_7(3, \boldsymbol{0}, 1, 9, 6, 10, \boldsymbol{11})$, ${\rm lk}(3) = C_8(0, \boldsymbol{1}, 2, \boldsymbol{10}, 11, \boldsymbol{9}, 4, \boldsymbol{5})$, ${\rm lk}(4) = C_7(3, \boldsymbol{0}, 5, 8, 1, 9, \boldsymbol{11})$, ${\rm lk}(9) = C_7(11$, $\boldsymbol{3}, 4, 1, 2, 6, \boldsymbol{7})$, ${\rm lk}(10) = C_7(11, \boldsymbol{3}, 2, 6, 5, 8, \boldsymbol{7})$. Then $M \cong \mathcal{F}_{3}(T)$ by the identity map.

\smallskip	

\noindent{\bf Subcase 1.4.}	$(x_1, x_2) = (5,6)$, i.e., ${\rm lk}(1) = C_7(0, \boldsymbol{3}, 2, 6, 5, 8, \boldsymbol{7})$. Then ${\rm lk}(5) = C_7(0, \boldsymbol{3}, 4, 8, 1, 6, \boldsymbol{7})$ and ${\rm lk}(6) = C_7(7, \boldsymbol{0}, 5, 1, 2, x_3, \boldsymbol{x_4})$. Here, we see that $(x_3, x_4) \in \{(4,3), (9, 10)\}$. If $(x_3, x_4) = (4,3)$, completing successively, we get 
${\rm lk}(4) = C_7(3, \boldsymbol{0}, 5, 8, 2, 6, \boldsymbol{7})$, ${\rm lk}(8) = C_7(7, \boldsymbol{0}, 1, 5, 4, 2, \boldsymbol{3})$, ${\rm lk}(2) = C_7(3, \boldsymbol{0}, 1, 6, 4, 8, \boldsymbol{7})$, ${\rm lk}(3) = C_8(0, \boldsymbol{1}, 2, \boldsymbol{8}, 7, \boldsymbol{6}, 4, \boldsymbol{5})$. Then $M \cong \mathcal{F}_4(K)$ by the identity map.
	
On the other hand, if $(x_3, x_4) = (9,10)$, completing successively, we get ${\rm lk}(2) = C_7(3, \boldsymbol{0}, 1, 6, 9, 11$, $\boldsymbol{10})$, ${\rm lk}(3) = C_8(0, \boldsymbol{1}, 2, \boldsymbol{11}, 10, \boldsymbol{9}, 4, \boldsymbol{5})$, ${\rm lk}(9) = C_7(10, \boldsymbol{3}, 4, 11, 2, 6, \boldsymbol{7})$, ${\rm lk}(4) = C_7(3, \boldsymbol{0}, 5, 8, 11, 9, \boldsymbol{10})$, ${\rm lk}(7) = C_8(0, \boldsymbol{1}, 8, \boldsymbol{11}, 10, \boldsymbol{9}, 6, \boldsymbol{5})$, ${\rm lk}(8) = C_7(7, \boldsymbol{0}, 1, 5, 4, 11, \boldsymbol{10})$, ${\rm lk}(11) = C_7(10, \boldsymbol{3}, 2, 9, 4, 8, \boldsymbol{7})$. Then $M \cong \mathcal{F}_1(K)$ by the map $(2,8)(3,7)(4,6)(9,11)$.

\smallskip	

\noindent{\bf Subcase 1.5.} $(x_1, x_2)=(5,4)$, i.e., ${\rm lk}(1) = C_7(0, \boldsymbol{3}, 2, 4, 5, 8, \boldsymbol{7})$. Then ${\rm lk}(5) = C_7(0, \boldsymbol{3}, 4, 1, 8, 6, \boldsymbol{7})$ and ${\rm lk}(8) = C_7(7, \boldsymbol{0}, 1, 5, 6$, $x_3, \boldsymbol{x_4})$, where $(x_3, x_4) \in \{(2,3), (9, 10)\}$. If $(x_3, x_4)=(2,3)$, completing successively, we get ${\rm lk}(8) = C_7(7, \boldsymbol{0}, 1, 5, 6, 2, \boldsymbol{3})$, ${\rm lk}(2) = C_7(3, \boldsymbol{0}, 1, 4, 6, 8, \boldsymbol{7})$, ${\rm lk}(4) = C_7(3, \boldsymbol{0}, 5, 1, 2$, $6, \boldsymbol{7})$, ${\rm lk}(6) = C_7(7, \boldsymbol{0}, 5, 8, 2, 4, \boldsymbol{3})$, ${\rm lk}(7) = C_8(0, \boldsymbol{1}, 8, \boldsymbol{2}, 3, \boldsymbol{4}, 6, \boldsymbol{5})$, ${\rm lk}(3) = C_8(0, \boldsymbol{1}, 2, \boldsymbol{8}, 7, \boldsymbol{6}, 4, \boldsymbol{5})$. Then $M \cong \mathcal{F}_5(T)$ by the identity map.
	
On the other hand, if $(x_3, x_4)=(9,10)$, then ${\rm lk}(8) = C_7(7, \boldsymbol{0}, 1, 5, 6, 9, \boldsymbol{10})$. This implies ${\rm lk}(6) = C_7(7, \boldsymbol{0}, 5, 8, 9$, $x_5, \boldsymbol{x_6})$, where $(x_5, x_6) \in \{(2,3), (11, 10)\}$. If $(x_5, x_6)=(2,3)$, then ${\rm lk}(6) = C_7(7, \boldsymbol{0}, 5, 8, 9, 2, \boldsymbol{3})$, ${\rm lk}(2) = C_7(3, \boldsymbol{0}, 1, 4, 9, 6, \boldsymbol{7})$, ${\rm lk}(9) = C_7(10, \boldsymbol{3}, 4, 2, 6, 8, \boldsymbol{7})$, which implies that ${\rm lk}(3)$ can not be completed. So $(x_5, x_6)=(11, 10)$. Completing successively, we get ${\rm lk}(6) = C_7(7, \boldsymbol{0}, 5, 8, 9, 11, \boldsymbol{10})$, ${\rm lk}(7) = C_8(0, \boldsymbol{1}, 8, \boldsymbol{9}, 10, \boldsymbol{11}, 6, \boldsymbol{5})$, ${\rm lk}(9) = C_7(10, \boldsymbol{3}, 2, 11, 6, 8, \boldsymbol{7})$, ${\rm lk}(2) = C_7(3, \boldsymbol{0}, 1, 4, 11, 9, \boldsymbol{10})$, ${\rm lk}(3) = C_8(0, \boldsymbol{1}, 2, \boldsymbol{9}, 10, \boldsymbol{11}, 4, \boldsymbol{5})$, ${\rm lk}(4) = C_7(3, \boldsymbol{0}, 5, 1, 2, 11, \boldsymbol{10})$, ${\rm lk}(11) = C_7(10, \boldsymbol{3}, 4, 2, 9, 6, \boldsymbol{7})$. Then $M \cong \mathcal{F}_6(T)$ by the identity map.

\smallskip	

\noindent{\bf Subcase 1.6.}	$(x_1, x_2) = (6,5)$, i.e., ${\rm lk}(1) = C_7(0, \boldsymbol{3}, 2, 5, 6, 8, \boldsymbol{7})$. Then ${\rm lk}(5) = C_7(0, \boldsymbol{3}, 4, 2, 1, 6, \boldsymbol{7})$. This implies 
${\rm lk}(6) = C_7(7, \boldsymbol{0}, 5, 1, 8, x_3, \boldsymbol{x_4})$, where $(x_3, x_4) \in \{(2,3), (4,3), (9, 10)\}$. If $(x_3, x_4)=(2,3)$, then considering ${\rm lk}(6)$ and ${\rm lk}(2)$, we see that ${\rm lk}(8)$ can not be completed. If $(x_3, x_4)=(4,3)$, then completing successively we get ${\rm lk}(4) = C_7(3, \boldsymbol{0}, 5, 2, 8, 6, \boldsymbol{7})$, ${\rm lk}(2) = C_7(3, \boldsymbol{0}, 1, 5, 4, 8, \boldsymbol{7})$, ${\rm lk}(8) = C_7(7, \boldsymbol{0}, 1, 6, 4, 2, \boldsymbol{3})$, ${\rm lk}(3) = C_8(2, \boldsymbol{1}$, $0, \boldsymbol{5}, 4, \boldsymbol{6}, 7, \boldsymbol{8})$. Then $M \cong \mathcal{F}_5(T)$ by the  map $(0,3)(1,2)(4,5)$.

On the other hand, if $(x_3, x_4)=(9, 10)$, then ${\rm lk}(6) = C_7(7, \boldsymbol{0}, 5, 1, 8, 9, \boldsymbol{10})$, ${\rm lk}(7) = C_8(0, \boldsymbol{1}, 8, \boldsymbol{11}$, $10, \boldsymbol{9}, 6, \boldsymbol{5})$, ${\rm lk}(8) = C_7(7, \boldsymbol{0}, 1, 6, 9, 11, \boldsymbol{10})$. This implies ${\rm lk}(9) = C_7(10, \boldsymbol{7}, 6, 8, 11, x_5, \boldsymbol{x_6})$, where $(x_5, x_6) \in \{(3,2), (4,3)\}$. In case $(x_5, x_6) = (3,2)$, considering ${\rm lk}(9)$ and ${\rm lk}(2)$, we see that ${\rm lk}(10)$ can not be completed. If $(x_5, x_6) = (4,3)$, completing successively, we get ${\rm lk}(9) = C_7(10, \boldsymbol{3}, 4, 11, 8, 6, \boldsymbol{7})$, ${\rm lk}(4) = C_7(3, \boldsymbol{0}, 5, 2, 11, 9, \boldsymbol{10})$, ${\rm lk}(3) = C_8(0, \boldsymbol{1}, 2, \boldsymbol{11}, 10, \boldsymbol{9}, 4, \boldsymbol{5})$, ${\rm lk}(2) = C_7(3, \boldsymbol{0}, 1, 5, 4, 11 \boldsymbol{10})$, ${\rm lk}(11) = C_7(10, \boldsymbol{3}, 2, 4, 9, 8, \boldsymbol{7})$. Then $M \cong \mathcal{F}_6(T)$ by the map $(0,10)(1,9,5,11)$.
	
\smallskip	

\noindent{\bf Subcase 1.7.}	$(x_1, x_2) = (6,9)$, i.e., ${\rm lk}(1) = C_7(0, \boldsymbol{3}, 2, 9, 6, 8, \boldsymbol{7})$. Then ${\rm lk}(6) = C_7(7, \boldsymbol{0}, 5, 8, 1, 9, \boldsymbol{x_3})$, where $x_3 \in \{4, 10\}$. If $x_3=4$, then considering ${\rm lk}(5)$ and ${\rm lk}(8)$, we see that ${\rm lk}(7)$ can not be completed. If $x_3=10$, then  ${\rm lk}(7) = C_8(0, \boldsymbol{1}, 8, \boldsymbol{11}, 10, \boldsymbol{9}, 6, \boldsymbol{5})$, ${\rm lk}(8) = C_7(7, \boldsymbol{0}, 1, 6, 5, 11, \boldsymbol{10})$, ${\rm lk}(5) = C_7(0, \boldsymbol{3}, 4, 11, 8, 6, \boldsymbol{7})$. This implies ${\rm lk}(11) = C_7(10, \boldsymbol{3}, 2, 4, 5, 8, \boldsymbol{7})$, completing successively, we get ${\rm lk}(2) = C_7(3, \boldsymbol{0}, 1, 9, 4, 11, \boldsymbol{10})$, ${\rm lk}(3) = C_8(0, \boldsymbol{1}, 2, \boldsymbol{11}, 10, \boldsymbol{9}, 4, \boldsymbol{5})$, ${\rm lk}(4) = C_7(3, \boldsymbol{0}, 5, 11, 2, 9, \boldsymbol{10})$, ${\rm lk}(9) = C_7(10, \boldsymbol{3}, 4, 2, 1, 6, \boldsymbol{7})$. Then $M \cong \mathcal{F}_{1}(K)$ by the map $(0,7,10)(1,8,9,4,5,6,11,2)$.

\smallskip	

\noindent{\bf Subcase 1.8.} $(x_1, x_2)=(9,4)$. Then ${\rm lk}(4) = C_7(3, \boldsymbol{0}, 5, 2, 1, 9, \boldsymbol{x_3})$, where $x_3 \in \{6, 10\}$. The case $x_3=6$ implies ${\rm lk}(3) = C_8(0, \boldsymbol{1}, 2, \boldsymbol{x_4}, 6, \boldsymbol{9}, 4, \boldsymbol{5})$. Observe that $x_4 \in \{7, 10\}$. If $x_4=7$, then considering ${\rm lk}(3)$, we see that ${\rm lk}(2)$ can not be completed. While if $x_4=10$, considering ${\rm lk}(2)$ and ${\rm lk}(5)$, we see that ${\rm lk}(6)$ can not be completed. On the other hand, if $x_3=10$, then ${\rm lk}(3) = C_8(0, \boldsymbol{1}, 2, \boldsymbol{x_4}, 10, \boldsymbol{9}, 4, \boldsymbol{5})$. Observe that $x_4 \in \{6, 11\}$. As above, we see that $x_4 \neq 6$. So $x_4=11$. Now completing successively, we get ${\rm lk}(3) = C_8(0, \boldsymbol{1}, 2, \boldsymbol{11}, 10, \boldsymbol{9}, 4, \boldsymbol{5})$, ${\rm lk}(6) = C_7(7, \boldsymbol{0}, 5, 11, 8, 9, \boldsymbol{10})$, ${\rm lk}(9) = C_7(10, \boldsymbol{3}, 4, 1, 8, 6, \boldsymbol{7})$, ${\rm lk}(8) = C_7(7, \boldsymbol{0}, 1, 9, 6, 11, \boldsymbol{10})$, ${\rm lk}(7) = C_8(0, \boldsymbol{1}, 8, \boldsymbol{11}, 10, \boldsymbol{9}, 6, \boldsymbol{5})$, ${\rm lk}(11)$ $= C_7(10, \boldsymbol{3}, 2, 5, 6, 8, \boldsymbol{7})$. Then $M \cong \mathcal{F}_{1}(K)$ by the map $(0,3)(1,2)(4,5)(6,11)(7,10)$.

\smallskip	

\noindent{\bf Subcase 1.9.}	$(x_1, x_2) = (9,6)$, i.e., ${\rm lk}(1) = C_7(0, \boldsymbol{3}, 2, 6, 9, 8, \boldsymbol{7})$. Then ${\rm lk}(6) = C_7(5, \boldsymbol{0}, 7, 2, 1, 9, \boldsymbol{x_3})$ or ${\rm lk}(6) = C_7(7, \boldsymbol{0}, 5, 2, 1, 9, \boldsymbol{x_3})$, for $x_3 \in V(M)$. In the first case, we see $x_3 \in \{4, 10\}$. But, a small computation shows that for both of these values $M$ does not exist. On the other hand, when ${\rm lk}(6) = C_7(7, \boldsymbol{0}, 5, 2, 1, 9, \boldsymbol{x_3})$, then $x_3 \in \{4, 10\}$. If $x_3=4$, then considering ${\rm lk}(6)$ and ${\rm lk}(9)$, we see ${\rm lk}(4)$ can not be completed. If $x_3=10$, then ${\rm lk}(6) = C_7(7, \boldsymbol{0}, 5, 2, 1, 9, \boldsymbol{10})$ and ${\rm lk}(7) = C_8(8, \boldsymbol{1}, 0, \boldsymbol{5}, 6, \boldsymbol{9}, 10, \boldsymbol{x_6})$, where $x_6 \in \{4, 11\}$. If $x_6=4$, considering ${\rm lk}(7)$ and ${\rm lk}(5)$, we see that ${\rm lk}(4)$ can not be completed. If $x_6=11$, completing successively, we get ${\rm lk}(8) = C_7(7, \boldsymbol{0}, 1, 9, 4, 11, \boldsymbol{10})$, ${\rm lk}(9) = C_7(10, \boldsymbol{3}, 4, 8, 1, 6, \boldsymbol{7})$, ${\rm lk}(4) = C_7(3, \boldsymbol{0}, 5, 11, 8, 9, \boldsymbol{10})$, ${\rm lk}(3) = C_8(0, \boldsymbol{1}$, $2, \boldsymbol{11}, 10, \boldsymbol{9}, 4, \boldsymbol{5})$, ${\rm lk}(11) = C_7(10, \boldsymbol{3}, 2, 5, 4, 8, \boldsymbol{7})$, ${\rm lk}(2) = C_7(3, \boldsymbol{0}, 1, 6, 5, 11, \boldsymbol{10})$, ${\rm lk}(10) = C_8(3, \boldsymbol{2}, 11$, $\boldsymbol{8}, 7, \boldsymbol{6}, 9, \boldsymbol{4})$. Then $M \cong \mathcal{F}_{3}(T)$ by the map $(0,3)(1,2)(4,5)(6,9)(7,11,8,10)$.

\smallskip	

\noindent{\bf Subcase 1.10.} $(x_1, x_2)=(9,10)$, i.e., ${\rm lk}(1) = C_7(0, \boldsymbol{3}, 2, 10, 9, 8, \boldsymbol{7})$. Then ${\rm lk}(2) = C_7(3, \boldsymbol{0}, 1, 10, x_3$, $x_4, \boldsymbol{x_5})$. Observe that, $(x_3, x_4, x_5) \in \{(4,8,7), (6, 11, 9), (8,9,6), (8,9,11), (11,6,7), (11,8,7)\}$. A small computation shows that $M$ does not exist for these values of $(x_3,x_4,x_5)$, except $(11,8,7)$. If $(x_3,x_4,x_5)=(11,8,7)$, then ${\rm lk}(2) = C_7(3, \boldsymbol{0}, 1, 10, 11, 8, \boldsymbol{7})$,  ${\rm lk}(8) = C_7(7, \boldsymbol{0}, 1, 9, 11, 2, \boldsymbol{3})$, ${\rm lk}(3) = C_8(0, \boldsymbol{1}, 2, \boldsymbol{8}, 7$, $\boldsymbol{6}, 4, \boldsymbol{5})$ and ${\rm lk}(7) = C_8(0, \boldsymbol{1}, 8, \boldsymbol{2}, 3, \boldsymbol{4}, 6, \boldsymbol{5})$. This implies ${\rm lk}(9) = C_7(x_7, \boldsymbol{x_8}, 11, 8, 1$, $10, \boldsymbol{x_6})$. It is easy to see that $(x_6, x_7, x_8) \in \{(4,5,6), (4,6,5), (5,4,6), (5,6,4), (6,4,5), (6,5,4)\}$. Note that $(5,4,6) \cong (4,6,5)$ by the map $(0,7,3)(1,8,2)(4,5,6)(9,11,10)$, $(5,6,4) \cong (4,5,6)$ by the map $(0,7)(1,8)(5,6)(10,11)$, and $(6,5,4) \cong (4,6,5)$ by the map $(0,7)(1,8)(5,6)(10,11)$. So, we search for $(x_6, x_7, x_8) \in \{(4,5,6), (4,6, 5), (6,4,5)\}$.

\smallskip	

\noindent{\bf Subcase 1.10.1.} If $(x_6, x_7, x_8) = (4,5,6)$, then completing successively, we get ${\rm lk}(9) = C_7(5, \boldsymbol{4}, 10$, $1, 8, 11, \boldsymbol{6})$, ${\rm lk}(5) = C_8(0, \boldsymbol{3}, 4, \boldsymbol{10}, 9, \boldsymbol{11}, 6, \boldsymbol{7})$, ${\rm lk}(4) = C_8(3, \boldsymbol{0}, 5, \boldsymbol{9}, 10, \boldsymbol{11}, 6, \boldsymbol{7})$, ${\rm lk}(6) = C_8(4, \boldsymbol{3}, 7, \boldsymbol{0}$, $5, \boldsymbol{9}, 11, \boldsymbol{10})$, ${\rm lk}(10) = C_7(4, \boldsymbol{5}, 9, 1, 2, 11, \boldsymbol{6})$, ${\rm lk}(11) = C_7(6, \boldsymbol{4}, 10, 2, 8, 9, \boldsymbol{5})$. Then $M \cong \mathcal{E}_7(T)$. 

\smallskip	

\noindent{\bf Subcase 1.10.2.} If $(x_6, x_7, x_8) = (4,6,5)$, then completing successively, we get ${\rm lk}(9) = C_7(6, \boldsymbol{4}, 10$, $1, 8, 11, \boldsymbol{5})$, ${\rm lk}(6) = C_8(4, \boldsymbol{3}, 7, \boldsymbol{0}, 5, \boldsymbol{11}, 9, \boldsymbol{10})$, ${\rm lk}(4) = C_8(3, \boldsymbol{0}, 5, \boldsymbol{11}, 10, \boldsymbol{9}, 6, \boldsymbol{7})$, ${\rm lk}(10) = C_7(4, \boldsymbol{5}, 11$, $2, 1, 9, \boldsymbol{6})$, ${\rm lk}(11) = C_7(5, \boldsymbol{4}, 10, 2, 8, 9, \boldsymbol{6})$, ${\rm lk}(5) = C_8(0, \boldsymbol{3}, 4, \boldsymbol{10}, 11, \boldsymbol{9}, 6, \boldsymbol{7})$. Then $M \cong \mathcal{F}_8(K)$. 

\smallskip	

\noindent{\bf Subcase 1.10.3.} If $(x_6, x_7, x_8) = (6,4,5)$, then completing successively, we get ${\rm lk}(9) = C_7(4, \boldsymbol{6}, 10$, $1, 8, 11, \boldsymbol{5})$, ${\rm lk}(4) = C_8(3, \boldsymbol{0}, 5, \boldsymbol{11}, 9, \boldsymbol{10}, 6, \boldsymbol{7})$, ${\rm lk}(5) = C_8(0, \boldsymbol{3}, 4, \boldsymbol{9}, 11, \boldsymbol{10}, 6, \boldsymbol{7})$, ${\rm lk}(11) = C_7(5, \boldsymbol{4}, 9, 8$, $2, 10, \boldsymbol{6})$, ${\rm lk}(10) = C_7(6, \boldsymbol{4}, 9, 1, 2, 11, \boldsymbol{5})$. Then $M \cong \mathcal{F}_9(T)$. 

\smallskip

\noindent{\bf Case 2:} If $f_{seq}(1) = (4^4)$, then ${\rm lk}(1) = C_8(2, \boldsymbol{3}, 0, \boldsymbol{7}, 8, \boldsymbol{x_1}, x_2, \boldsymbol{x_3})$. We see that, $(x_1, x_2, x_3) \in \{(4,5,6)$, $(4,6,5)$, $(4,9,6)$, $(4,9,10)$, $(5,4,6)$, $(5,4,9)$, $(6,4,5)$, $(6,5,4)$, $(6,9,4)$, $(6,9,10)$, $(9,4,5)$, $(9,6,5)$, $(9,10,4)$, $(9,10,6)$, $(9,10,11)\}$. Here $(4,5,6) \cong (6,5,4)$ and $(4,9,6) \cong (6,9,4)$ by the map $(2,8)(3,7)(4,6)$. Hence, we search for the following cases:

\noindent{\bf Subcase 2.1} If $(x_1,x_2,x_3)=(4,5,6)$, then ${\rm lk}(1) = C_8(2, \boldsymbol{3}, 0, \boldsymbol{7}, 8, \boldsymbol{4}, 5, \boldsymbol{6})$ and ${\rm lk}(5) = C_8(0, \boldsymbol{7}, 6, \boldsymbol{2}$, $1, \boldsymbol{8}, 4, \boldsymbol{3})$. This gives $f_{seq}(2)=(3^3.4^2)$ or $f_{seq}(2)=(4^4)$.  
	
If $f_{seq}(2)=(3^3.4^2)$, then ${\rm lk}(2) = C_7(1, \boldsymbol{0}, 3, x_5, x_4,6, \boldsymbol{5})$, where $(x_4,x_5) \in \{(4,8), (9,10)\}$. If $(x_4,x_5)=(4,8)$, then considering ${\rm lk}(4)$, we get $\deg(6)>5$. If $(x_4,x_5)=(9,10)$, then ${\rm lk}(6) = C_7(5, \boldsymbol{0}, 7, 11, 9,2, \boldsymbol{1})$, ${\rm lk}(7) = C_7(0, \boldsymbol{1}, 8, 4, 11,6, \boldsymbol{5})$, and ${\rm lk}(4) = C_7(5, \boldsymbol{0}, 3, 11, 7,8, \boldsymbol{1})$. Now observe that there are 4 triangular faces incident at 11, which contradicts the fact that $f_{seq}(11)=(3^3.4^2)$. 
	
On the other hand, if $f_{seq}(2)=(4^4)$, then ${\rm lk}(2) = C_8(3, \boldsymbol{0}, 1, \boldsymbol{5},6,\boldsymbol{x_4}, x_5, \boldsymbol{x_6})$, where $x_5 \in \{8,9\}$. If $x_5=8$, then completing successively, we get ${\rm lk}(2) = C_8(3, \boldsymbol{0}, 1, \boldsymbol{5},6,\boldsymbol{4}, 8, \boldsymbol{7})$, 
${\rm lk}(8) = C_8(1, \boldsymbol{0}, 7, \boldsymbol{3},2,\boldsymbol{6}, 4, \boldsymbol{5})$, 
${\rm lk}(3) = C_8(0, \boldsymbol{1}, 2, \boldsymbol{8},7,\boldsymbol{6}, 4, \boldsymbol{5})$, 
${\rm lk}(7) = C_8(0, \boldsymbol{1}, 8, \boldsymbol{2},3,\boldsymbol{4}, 6, \boldsymbol{5})$,
${\rm lk}(6) = C_8(2$, $\boldsymbol{1}, 5, \boldsymbol{0},7,\boldsymbol{3}, 4, \boldsymbol{8})$, and 
${\rm lk}(4) = C_8(3, \boldsymbol{0}, 5, \boldsymbol{1},8,\boldsymbol{2}, 6, \boldsymbol{7})$. This gives a semi equivelar map of type $[4^4]$. 
	
If $x_5=9$, then $x_4=10$ and $x_6=11$. Completing successively, we get  ${\rm lk}(2) = C_8(3, \boldsymbol{0}, 1, \boldsymbol{5},6,\boldsymbol{10}$, $9, \boldsymbol{11})$, ${\rm lk}(3) = C_8(0, \boldsymbol{1}, 2, \boldsymbol{9},11,\boldsymbol{10}, 4, \boldsymbol{5})$, ${\rm lk}(4) = C_8(3, \boldsymbol{0}, 5, \boldsymbol{1},8,\boldsymbol{9}, 10, \boldsymbol{11})$, ${\rm lk}(8) = C_8(1, \boldsymbol{0}, 7, \boldsymbol{11},9$, $\boldsymbol{10}, 4, \boldsymbol{5})$, ${\rm lk}(7) = C_8(0, \boldsymbol{1}, 8, \boldsymbol{9},11,\boldsymbol{10}, 6, \boldsymbol{5})$, ${\rm lk}(10) = C_8(4, \boldsymbol{3}, 11, \boldsymbol{7},6,\boldsymbol{2}, 9, \boldsymbol{8})$, ${\rm lk}(9) = C_8(2, \boldsymbol{3}, 11, \boldsymbol{7},8$, $\boldsymbol{4}, 10, \boldsymbol{6})$, ${\rm lk}(11) = C_8(3, \boldsymbol{2}, 9, \boldsymbol{8},7,\boldsymbol{6}, 10, \boldsymbol{4})$, and ${\rm lk}(6) = C_8(2, \boldsymbol{1}, 5, \boldsymbol{0},7,\boldsymbol{11}, 10, \boldsymbol{9})$. This gives a semi equivelar map of type $[4^4]$. 
	
\smallskip
	
\noindent{\bf Subcase 2.2} If $(x_1,x_2,x_3)=(4,6,5)$, then ${\rm lk}(1) = C_8(2, \boldsymbol{3}, 0, \boldsymbol{7}, 8, \boldsymbol{4}, 6, \boldsymbol{5})$ and ${\rm lk}(6) = C_8(1, \boldsymbol{2}, 5, \boldsymbol{0}$, $7, \boldsymbol{3}, 4, \boldsymbol{8})$. It is easy to see that ${\rm lk}(4) = C_8(3, \boldsymbol{0}, 5, \boldsymbol{2}$, $8, \boldsymbol{1}, 6, \boldsymbol{7})$. Completing successively, we get ${\rm lk}(5) = C_8(0, \boldsymbol{3}, 4, \boldsymbol{8}$, $2, \boldsymbol{1}, 6, \boldsymbol{7})$, ${\rm lk}(2) = C_8(1, \boldsymbol{0}, 3, \boldsymbol{7}$, $8, \boldsymbol{4}, 5, \boldsymbol{6})$, ${\rm lk}(7) = C_8(0, \boldsymbol{1}, 8, \boldsymbol{2},3,\boldsymbol{4}, 6, \boldsymbol{5})$, ${\rm lk}(3) = C_8(0, \boldsymbol{1}, 2, \boldsymbol{8},7,\boldsymbol{6}, 4, \boldsymbol{5})$, ${\rm lk}(8) = C_8(1, \boldsymbol{0}, 7, \boldsymbol{3},2,\boldsymbol{5}, 4, \boldsymbol{6})$. This gives a semi equivelar map of type $[4^4]$. 

\smallskip

\noindent{\bf Subcase 2.3} If $(x_1,x_2,x_3)=(4,9,6)$, then ${\rm lk}(1) = C_8(2, \boldsymbol{3}, 0, \boldsymbol{7}, 8, \boldsymbol{4}, 9, \boldsymbol{6})$. This implies ${\rm lk}(6) = C_8(2, \boldsymbol{1}, 9, \boldsymbol{x_4},7,\boldsymbol{0}, 5, \boldsymbol{x_5})$ or ${\rm lk}(6) = C_8(2, \boldsymbol{1}, 9$, $\boldsymbol{x_4},5,\boldsymbol{0}, 7, \boldsymbol{x_5})$. In case ${\rm lk}(6) = C_8(2, \boldsymbol{1}, 9, \boldsymbol{x_4},7,\boldsymbol{0}, 5$, $\boldsymbol{x_5})$, we see that $(x_4,x_5) \in \{(4,10), (10,11)\}$. If $(x_4,x_5)=(4,10)$, then considering ${\rm lk}(5)$, we see that ${\rm lk}(9)$ can not be completed. If $(x_4,x_5)=(10,11)$, then it is easy to see that ${\rm lk}(7) = C_8(0, \boldsymbol{1}, 8, \boldsymbol{4},10,\boldsymbol{9}, 6, \boldsymbol{5})$. Now completing successively we get ${\rm lk}(8) = C_8(1, \boldsymbol{0}, 7, \boldsymbol{10},11,\boldsymbol{5}, 4, \boldsymbol{9})$, ${\rm lk}(5) = C_8(0, \boldsymbol{3}, 4, \boldsymbol{8}, 11, \boldsymbol{2}, 6, \boldsymbol{7})$, ${\rm lk}(4) = C_8(3, \boldsymbol{0}, 5, \boldsymbol{11}, 8, \boldsymbol{1}, 9, \boldsymbol{10})$. Then $M$ is a semi-equivelar map of type $[4^4]$. On the other hand if ${\rm lk}(6) = C_8(2, \boldsymbol{1}, 9, \boldsymbol{x_4},5,\boldsymbol{0}, 7, \boldsymbol{x_5})$, then $x_4=10$ and $x_5=11$. Now considering ${\rm lk}(5)$ and ${\rm lk}(4)$, we see that ${\rm lk}(3)$ can not be completed. 

\smallskip
		
\noindent{\bf Subcase 2.4} If $(x_1,x_2,x_3)=(4,9,10)$, then ${\rm lk}(4) = C_8(3, \boldsymbol{0}, 5, \boldsymbol{x_4},8,\boldsymbol{1}, 9, \boldsymbol{x_5})$ or ${\rm lk}(4) = C_8(3, \boldsymbol{0}, 5$, $\boldsymbol{x_4},9,\boldsymbol{1}, 8, \boldsymbol{x_5})$. In case ${\rm lk}(4) = C_8(3, \boldsymbol{0}, 5, \boldsymbol{x_4},8,\boldsymbol{1}, 9, \boldsymbol{x_5})$, we see $(x_4,x_5) \in \{(10,6), (10,11),(11,6)\}$. If $(x_4,x_5)=(10,6)$, then ${\rm lk}(3) = C_8(0, \boldsymbol{1}, 2, \boldsymbol{7},6,\boldsymbol{9}, 4, \boldsymbol{5})$, ${\rm lk}(2) = C_8(1, \boldsymbol{0}, 3, \boldsymbol{6}, 7, \boldsymbol{8}, 10, \boldsymbol{9})$, and ${\rm lk}(7) = C_8(0, \boldsymbol{1}, 8, \boldsymbol{10},2,\boldsymbol{3}, 6, \boldsymbol{5})$. Now observe that ${\rm lk}(8)$ can not be completed. Similarly, we see that links of all the vertices can not be completed for  $(x_4,x_5) \in \{(10,11),(11,6)\}$. On the other hand when ${\rm lk}(4) = C_8(3, \boldsymbol{0}, 5$, $\boldsymbol{x_4},9,\boldsymbol{1}, 8, \boldsymbol{x_5})$ then $(x_4, x_5) \in \{(11,10), (11,6)\}$. If $(x_4,x_5)=(11,10)$, then ${\rm lk}(3) = C_8(0, \boldsymbol{1}, 2, \boldsymbol{9},10,\boldsymbol{8}, 4, \boldsymbol{5})$, which implies $C_5(0,1,9,10,3) \subseteq {\rm lk}(2)$. Similarly, for $(x_4,x_5)=(11,6)$, considering ${\rm lk}(3)$ and ${\rm lk}(2)$, we see that ${\rm lk}(8)$ can not be completed.  

\smallskip
	
\noindent{\bf Subcase 2.5} If $(x_1,x_2,x_3)=(5,4,6)$, then ${\rm lk}(4) = C_8(3, \boldsymbol{0}, 5, \boldsymbol{8},1,\boldsymbol{2}, 6, \boldsymbol{x_4})$, where we see that $x_4 \in \{7,9\}$. If $x_4=9$, then ${\rm lk}(3) = C_8(0, \boldsymbol{1}, 2, \boldsymbol{10},9,\boldsymbol{6}, 4, \boldsymbol{5})$, now observe that ${\rm lk}(2)$ can not be completed. If $x_4=7$, then completing successively, we get ${\rm lk}(6) = C_8(2, \boldsymbol{1}, 4, \boldsymbol{3},7,\boldsymbol{0}, 5, \boldsymbol{8})$, ${\rm lk}(5) = C_8(0, \boldsymbol{3}, 4, \boldsymbol{1}, 8, \boldsymbol{2}, 6, \boldsymbol{7})$, ${\rm lk}(3) = C_8(0, \boldsymbol{1}, 2, \boldsymbol{8},7,\boldsymbol{6}, 4, \boldsymbol{5})$, ${\rm lk}(2) = C_8(1, \boldsymbol{0}, 3, \boldsymbol{7}, 8, \boldsymbol{5}, 6, \boldsymbol{4})$, ${\rm lk}(7) = C_8(0$, $\boldsymbol{1}, 8, \boldsymbol{2},3,\boldsymbol{4}, 6, \boldsymbol{5})$, and ${\rm lk}(8) = C_8(1, \boldsymbol{0}, 7, \boldsymbol{3},2,\boldsymbol{6}, 5, \boldsymbol{4})$. This gives a semi equivelar map of type $[4^4]$. 

\smallskip

\noindent{\bf Subcase 2.6} If $(x_1,x_2,x_3)=(5,4,9)$, then ${\rm lk}(4) = C_8(3, \boldsymbol{0}, 5, \boldsymbol{8},1,\boldsymbol{2}, 9, \boldsymbol{x_4})$, where $x_4 \in \{6,10\}$. If $x_4=6$, then it is easy to see that 
${\rm lk}(3) = C_8(0, \boldsymbol{1}, 2, \boldsymbol{7},6,\boldsymbol{9}, 4, \boldsymbol{5})$ and 
${\rm lk}(7) = C_8(0,\boldsymbol{1}, 8, \boldsymbol{9},2,\boldsymbol{3}, 6, \boldsymbol{5})$, now completing successively, we get ${\rm lk}(8) = C_8(1, \boldsymbol{0}, 7, \boldsymbol{2},9,\boldsymbol{6}, 5, \boldsymbol{4})$,  
${\rm lk}(5) = C_8(0, \boldsymbol{3}, 4, \boldsymbol{1}, 8, \boldsymbol{9}, 6, \boldsymbol{7})$, 
${\rm lk}(6) = C_8(5, \boldsymbol{0}, 7, \boldsymbol{2},3,\boldsymbol{4}, 9, \boldsymbol{8})$, and 
${\rm lk}(9) = C_8(2, \boldsymbol{1}, 4, \boldsymbol{3},6,\boldsymbol{5}, 8, \boldsymbol{7})$. This gives a semi equivelar map of type $[4^4]$. 
	
If $x_4=10$, then ${\rm lk}(3) = C_8(0, \boldsymbol{1}, 2, \boldsymbol{x_5},10,\boldsymbol{9}, 4, \boldsymbol{5})$. This gives $x_5 \in \{6,11\}$. If $x_5=6$, then completing successively we get  ${\rm lk}(3) = C_8(0, \boldsymbol{1}, 2, \boldsymbol{6},10,\boldsymbol{9}, 4, \boldsymbol{5})$, 
${\rm lk}(2) = C_8(1, \boldsymbol{0}, 3, \boldsymbol{10},6,\boldsymbol{7}, 9, \boldsymbol{4})$,
${\rm lk}(7) = C_8(0,\boldsymbol{1}, 8, \boldsymbol{10},9,\boldsymbol{2}, 6, \boldsymbol{5})$, 
${\rm lk}(8) = C_8(1, \boldsymbol{0}, 7, \boldsymbol{9},10,\boldsymbol{6}, 5, \boldsymbol{4})$, 
${\rm lk}(5) = C_8(0, \boldsymbol{3}, 4, \boldsymbol{1}, 8, \boldsymbol{10}, 6, \boldsymbol{7})$,
${\rm lk}(6) = C_8(5, \boldsymbol{0}, 7, \boldsymbol{9},2,\boldsymbol{3}, 10, \boldsymbol{8})$,
${\rm lk}(9) = C_8(2, \boldsymbol{1}, 4, \boldsymbol{3},10,\boldsymbol{8}, 7, \boldsymbol{6})$, 
${\rm lk}(10) = C_8(3, \boldsymbol{2}, 6, \boldsymbol{5},8,\boldsymbol{7}, 9, \boldsymbol{4})$. This gives a semi equivelar map of type $[4^4]$. 

If $x_5=11$, then  ${\rm lk}(3) = C_8(0, \boldsymbol{1}, 2, \boldsymbol{11},10,\boldsymbol{9}, 4, \boldsymbol{5})$. Now completing successively, we get 
${\rm lk}(2) = C_8(1, \boldsymbol{0}, 3, \boldsymbol{10},11,\boldsymbol{6}, 9, \boldsymbol{4})$,
${\rm lk}(9) = C_8(2, \boldsymbol{1}, 4, \boldsymbol{3},10,\boldsymbol{7}, 6, \boldsymbol{11})$,
${\rm lk}(6) = C_8(5, \boldsymbol{0}, 7, \boldsymbol{10},9,\boldsymbol{2}, 11, \boldsymbol{8})$,
${\rm lk}(5) = C_8(0, \boldsymbol{3}, 4, \boldsymbol{1}, 8, \boldsymbol{11}, 6, \boldsymbol{7})$,
${\rm lk}(7) = C_8(0,\boldsymbol{1}, 8, \boldsymbol{11},10,\boldsymbol{9}, 6, \boldsymbol{5})$,
${\rm lk}(10) = C_8(3, \boldsymbol{2}, 11, \boldsymbol{8},7,\boldsymbol{6}, 9, \boldsymbol{4})$, and 
${\rm lk}(11) = C_8(2, \boldsymbol{3}, 10, \boldsymbol{7},8,\boldsymbol{5}, 6, \boldsymbol{9})$. This gives a semi equivelar map of type $[4^4]$. 

\smallskip

\noindent{\bf Subcase 2.7} If $(x_1,x_2,x_3)=(6,4,5)$, then it is easy to see that ${\rm lk}(4) = C_8(3, \boldsymbol{0}, 5, \boldsymbol{2},1,\boldsymbol{8}, 6, \boldsymbol{7})$. Now completing successively, we get 
${\rm lk}(6) = C_8(5, \boldsymbol{0}, 7, \boldsymbol{3},4,\boldsymbol{1}, 8, \boldsymbol{2})$,
${\rm lk}(3) = C_8(0, \boldsymbol{1}, 2, \boldsymbol{8},7,\boldsymbol{6}, 4, \boldsymbol{5})$,
${\rm lk}(2) = C_8(1, \boldsymbol{0}, 3, \boldsymbol{7},8,\boldsymbol{6}, 5, \boldsymbol{4})$,
${\rm lk}(8) = C_8(1, \boldsymbol{0}, 7, \boldsymbol{3},2,\boldsymbol{5}, 6, \boldsymbol{4})$,
${\rm lk}(5) = C_8(0, \boldsymbol{3}, 4, \boldsymbol{1}, 2, \boldsymbol{8}, 6, \boldsymbol{7})$,
${\rm lk}(7) = C_8(0,\boldsymbol{1}, 8, \boldsymbol{2},3,\boldsymbol{4}, 6, \boldsymbol{5})$.  This gives a semi equivelar map of type $[4^4]$. 

\noindent{\bf Subcase 2.8} If $(x_1,x_2,x_3)=(6,9,10)$, then ${\rm lk}(6) = C_8(5, \boldsymbol{x_4}, 8, \boldsymbol{1},9,\boldsymbol{x_5}, 7, \boldsymbol{0})$, where we see easily that $(x_4,x_5) \in \{(10,4), (10,11), (11,4)\}$. But for each value of $(x_4,x_5)$, we get a vertex whose link can not be completed.

\noindent{\bf Subcase 2.9} If $(x_1,x_2,x_3)=(9,4,5)$, then ${\rm lk}(4) = C_8(3, \boldsymbol{0}, 5, \boldsymbol{2},1,\boldsymbol{8}, 9, \boldsymbol{x_4})$, where $x_4 \in \{6,10\}$. If $x_4=6$ then considering ${\rm lk}(4)$ and ${\rm lk}(3)$, we see that ${\rm lk}(6)$ can not be completed. If $x_4=10$, then ${\rm lk}(4) = C_8(3, \boldsymbol{0}, 5, \boldsymbol{2},1,\boldsymbol{8}, 9, \boldsymbol{10})$, this implies ${\rm lk}(3) = C_8(2, \boldsymbol{1}, 0, \boldsymbol{5},4,\boldsymbol{9}, 10, \boldsymbol{x_5})$. It is easy to see that $x_5=11$. Then successively, we get ${\rm lk}(3) = C_8(2, \boldsymbol{1}, 0, \boldsymbol{5},4,\boldsymbol{9}, 10, \boldsymbol{11})$, ${\rm lk}(2) = C_8(1, \boldsymbol{0}, 3, \boldsymbol{10},11,\boldsymbol{6}, 5, \boldsymbol{4})$, ${\rm lk}(5) = C_8(0, \boldsymbol{3}, 4, \boldsymbol{1}, 2, \boldsymbol{11}, 6, \boldsymbol{7})$, ${\rm lk}(4) = C_8(1, \boldsymbol{2}, 5, \boldsymbol{0}, 3, \boldsymbol{10}, 9, \boldsymbol{8})$. Then $f_{seq}(6) = (3^3.4^2)$ or $(4^4)$. 

In case $f_{seq}(6)=(3^3.4^2)$, we get ${\rm lk}(6) = C_7(5, \boldsymbol{2}, 11, x_6, x_7,7, \boldsymbol{0})$. Here $(x_6,x_7) \in \{(8,9), (9,10)\}$. If $(x_6,x_7)=(8,9)$, then considering ${\rm lk}(9)$, we see that ${\rm lk}(7)$ can not be completed. If $(x_6,x_7)=(9,10)$, then 
considering ${\rm lk}(9)$ and ${\rm lk}(11)$, we see that $\deg(8)>5$. 

On the other hand when $f_{seq}(6)=(4^4)$, then ${\rm lk}(6) = C_8(7, \boldsymbol{0}, 5, \boldsymbol{2}, 11, \boldsymbol{x_6}, x_7, \boldsymbol{x_8})$. Observe that $x_7=9$, this implies $x_6=8$ and $x_8=10$. Now completing successively, we get  
${\rm lk}(9) = C_8(6, \boldsymbol{7}, 10, \boldsymbol{3},4,\boldsymbol{1}, 8, \boldsymbol{11})$,
${\rm lk}(7) = C_8(0,\boldsymbol{1}, 8, \boldsymbol{11},10,\boldsymbol{9}, 6, \boldsymbol{5})$,
${\rm lk}(10) = C_8(3, \boldsymbol{2}, 11, \boldsymbol{8},7,\boldsymbol{6}, 9, \boldsymbol{4})$, and
${\rm lk}(11) = C_8(2, \boldsymbol{3}, 10, \boldsymbol{7},8,\boldsymbol{9}, 6, \boldsymbol{5})$. This gives a semi equivelar map of type $[4^4]$. 

\smallskip

\noindent{\bf Subcase 2.10} If $(x_1,x_2,x_3)=(9,6,5)$, then ${\rm lk}(6) = C_8(7, \boldsymbol{0}, 5, \boldsymbol{2},1,\boldsymbol{8}, 9, \boldsymbol{x_4})$, where $x_4 \in \{4,10\}$. 

If $x_4=4$, then completing successively, we get 
${\rm lk}(6) = C_8(5, \boldsymbol{0}, 7, \boldsymbol{4},9,\boldsymbol{8}, 1, \boldsymbol{2})$, 
${\rm lk}(7) = C_8(0,\boldsymbol{1}, 8$, $\boldsymbol{3},4,\boldsymbol{9}, 6, \boldsymbol{5})$, 
${\rm lk}(8) = C_8(1, \boldsymbol{0}, 7, \boldsymbol{4},3,\boldsymbol{2}, 9, \boldsymbol{6})$,
${\rm lk}(3) = C_8(2, \boldsymbol{1}, 0, \boldsymbol{5},4,\boldsymbol{7}, 8, \boldsymbol{9})$, 
${\rm lk}(2) = C_8(1, \boldsymbol{0}, 3, \boldsymbol{8},9,\boldsymbol{4}$, $5, \boldsymbol{6})$,
${\rm lk}(5) = C_8(0, \boldsymbol{3}, 4, \boldsymbol{9}, 2, \boldsymbol{1}, 6, \boldsymbol{7})$,
${\rm lk}(4) = C_8(9, \boldsymbol{2}, 5, \boldsymbol{0}, 3, \boldsymbol{8}, 7, \boldsymbol{6})$. This gives a semi equivelar map of type $[4^4]$. 

If $x_4=10$, then ${\rm lk}(6) = C_8(5, \boldsymbol{0}, 7, \boldsymbol{10},9,\boldsymbol{8}, 1, \boldsymbol{2})$. This implies ${\rm lk}(7) = C_8(0,\boldsymbol{1}, 8, \boldsymbol{x_5},10,\boldsymbol{9}, 6, \boldsymbol{5})$, where $x_5 \in \{4,11\}$. If $x_5=4$, it is easy to see that ${\rm lk}(8) = C_8(1, \boldsymbol{0}, 7, \boldsymbol{10},4,\boldsymbol{3}, 9, \boldsymbol{6})$, now completing successively we get 
${\rm lk}(9) = C_8(6, \boldsymbol{7}, 10, \boldsymbol{2},3,\boldsymbol{4}, 8, \boldsymbol{1})$, 
${\rm lk}(3) = C_8(2, \boldsymbol{1}, 0, \boldsymbol{5},4,\boldsymbol{8}, 9, \boldsymbol{10})$,
${\rm lk}(2) = C_8(1, \boldsymbol{0}, 3, \boldsymbol{9},10,\boldsymbol{4},5, \boldsymbol{6})$,
${\rm lk}(5) = C_8(0, \boldsymbol{3}, 4, \boldsymbol{10}, 2, \boldsymbol{1}, 6, \boldsymbol{7})$,
${\rm lk}(6) = C_8(5, \boldsymbol{0}, 7, \boldsymbol{10},9,\boldsymbol{8}, 1, \boldsymbol{2})$,
${\rm lk}(4) = C_8(10, \boldsymbol{2}, 5, \boldsymbol{0}, 3, \boldsymbol{9}, 8, \boldsymbol{7})$. This gives a semi equivelar map of type $[4^4]$. 

If $x_5=11$, it is easy to see that ${\rm lk}(8) = C_8(1, \boldsymbol{0}, 7, \boldsymbol{10},11,\boldsymbol{4}, 9, \boldsymbol{6})$. Now completing successively, we get 
${\rm lk}(9) = C_8(6, \boldsymbol{7}, 10, \boldsymbol{3},4,\boldsymbol{11}, 8, \boldsymbol{1})$, 
${\rm lk}(4) = C_8(3, \boldsymbol{0}, 5, \boldsymbol{2}, 11, \boldsymbol{8}, 9, \boldsymbol{10})$,
${\rm lk}(5) = C_8(0, \boldsymbol{3}, 4, \boldsymbol{11}, 2, \boldsymbol{1}$, $6, \boldsymbol{7})$,
${\rm lk}(6) = C_8(5, \boldsymbol{0}, 7, \boldsymbol{10},9,\boldsymbol{8}, 1, \boldsymbol{2})$,
${\rm lk}(3) = C_8(2, \boldsymbol{1}, 0, \boldsymbol{5},4,\boldsymbol{9}, 10, \boldsymbol{11})$,
${\rm lk}(2) = C_8(1, \boldsymbol{0}, 3, \boldsymbol{10},11,\boldsymbol{4}$, $5, \boldsymbol{6})$. This gives a semi equivelar map of type $[4^4]$. 

\smallskip
 
\noindent{\bf Subcase 2.11} If $(x_1,x_2,x_3)=(9,10,4)$, then ${\rm lk}(4) = C_8(3, \boldsymbol{0}, 5, \boldsymbol{x_4},2,\boldsymbol{1}, 10, \boldsymbol{x_5})$, where we see that $(x_4,x_5) \in \{(9,6), (9,11), (11,6)\}$. If $(x_4,x_5)=(9,6)$, then ${\rm lk}(5)$ can not be completed. If $(x_4,x_5)=(9,11)$, then ${\rm lk}(3)$ can not be completed. If $(x_4,x_5)=(11,6)$, then considering ${\rm lk}(5)$ and ${\rm lk}(3)$, we see that ${\rm lk}(6)$ can not be completed. 

\smallskip

\noindent{\bf Subcase 2.12} If $(x_1,x_2,x_3)=(9,10,6)$, then ${\rm lk}(6) = C_8(10, \boldsymbol{1}$, $2, \boldsymbol{x_4},7,\boldsymbol{0}, 5, \boldsymbol{x_5})$ or ${\rm lk}(6) = C_8(2, \boldsymbol{1}, 10, \boldsymbol{x_4},7,\boldsymbol{0}, 5, \boldsymbol{x_5})$. If ${\rm lk}(6) = C_8(10, \boldsymbol{1}, 2, \boldsymbol{x_4},7,\boldsymbol{0}, 5, \boldsymbol{x_5})$, then $(x_4,x_5) \in \{(4,11), (9,11)\}$. But for both the cases of $(x_4,x_5)$, we get a vertex whose link can not be completed. On the other hand if ${\rm lk}(6) = C_8(2, \boldsymbol{1}, 10, \boldsymbol{x_4},7,\boldsymbol{0}, 5, \boldsymbol{x_5})$, then $(x_4,x_5) \in \{(4,9), (4,11)\}$. If $(x_4,x_5)=(4,11)$, then considering ${\rm lk}(6)$ and ${\rm lk}(7)$, we see that ${\rm lk}(8)$ can not be completed. If $(x_4,x_5)=(4,9)$, then completing successively we get
${\rm lk}(6) = C_8(5, \boldsymbol{0}, 7, \boldsymbol{4},10,\boldsymbol{1}, 2, \boldsymbol{9})$,
${\rm lk}(5) = C_8(0, \boldsymbol{3}, 4, \boldsymbol{10}, 9, \boldsymbol{2}, 6, \boldsymbol{7})$,
${\rm lk}(9) = C_8(2, \boldsymbol{3}, 8, \boldsymbol{1},10,\boldsymbol{4}, 5, \boldsymbol{6})$, 
${\rm lk}(3) = C_8(2, \boldsymbol{1}, 0, \boldsymbol{5},4,\boldsymbol{7}, 8, \boldsymbol{9})$,
${\rm lk}(4) = C_8(3, \boldsymbol{0}, 5, \boldsymbol{9}, 10, \boldsymbol{6}, 7, \boldsymbol{8})$,
${\rm lk}(2) = C_8(1, \boldsymbol{0}, 3, \boldsymbol{8},9,\boldsymbol{5},6, \boldsymbol{10})$. This gives a semi equivelar map of type $[4^4]$.

\smallskip

\noindent{\bf Subcase 2.13.} If $(x_1,x_2,x_3)=(9,10,11)$, then $f_{seq}(8)=(3^3.4^2)$ or $f_{seq}(8)=(4^4)$. 

\smallskip

\noindent{\bf Subcase 2.13.1.} Let $f_{seq}(8)=(3^3.4^2)$. Then ${\rm lk}(8)=C_7(1,\boldsymbol{0},7, x_5,x_4,9, \boldsymbol{10})$. Observe that $(x_4, x_5) \in \{(2,3), (2,11), (3,2), (3,4), (5,4), (6,4), (11,4), (11,2)\}$. So, we need not consider the last case. 

\noindent{\bf Claim 1.} $(x_4,x_5) = (2,3), (3,2), (3,4)$ or $(11,4)$. 

If $(x_4,x_5)=(2,11)$, then ${\rm lk}(8)=C_7(1,\boldsymbol{0},7, 11,2,9, \boldsymbol{10})$. This implies 
${\rm lk}(2)=C_7(1,\boldsymbol{0},3, 9,8,11$, $\boldsymbol{10})$ and ${\rm lk}(3)=C_7(0,\boldsymbol{1},2, 9,x_6,4$, $\boldsymbol{5})$, where $x_6 \in \{6, 10\}$. If $x_6=6$, we get ${\rm lk}(3)=C_7(0,\boldsymbol{1},2, 9,6,4$, $\boldsymbol{5})$ and ${\rm lk}(9)=C_7(10,\boldsymbol{1},8, 2,3,6, \boldsymbol{x_7})$, where $x_7 \in \{5,7\}$. In case $x_7=5$, considering ${\rm lk}(9)$,  ${\rm lk}(6)$ and ${\rm lk}(7)$, we get ${\rm lk}(4)=C_7(5,\boldsymbol{0},3, 6,7,11, \boldsymbol{x_8})$, but observe that $x_8$ can not be replaced by any vertex in $V(M)$. If $x_7=7$, then considering ${\rm lk}(9)$ and ${\rm lk}(6)$, we get $C_4(0,3,6,5) \subseteq {\rm lk}(4)$. On the other hand for $x_6=10$, considering ${\rm lk}(3)$ and ${\rm lk}(10)$, we get $C_5(1,8,2,3,10) \subseteq {\rm lk}(9)$.

If $(x_4,x_5)=(5,4)$, then ${\rm lk}(8)=C_7(1,\boldsymbol{0},7, 4,5,9, \boldsymbol{10})$. This implies
${\rm lk}(5)=C_7(0,\boldsymbol{3},4, 8,9,6, \boldsymbol{7})$. Now considering ${\rm lk}(6)$, we get quadrangular face $[0,5,6,7]$ which is adjacent with triangular faces $[5,6,9]$ and $[6,7,x_6]$. This implies ${\rm lk}(6)$ is not possible. 

 If $(x_4,x_5)=(6,4)$, then ${\rm lk}(8)=C_7(1,\boldsymbol{0},7, 4,6,9, \boldsymbol{10})$. This implies
${\rm lk}(7)=C_7(0,\boldsymbol{1},8, 4,x_6,6, \boldsymbol{5})$, where $x_6 \in \{3,9\}$. If $x_6=3$, then considering ${\rm lk}(4)$, we get 4 triagular faces $[6,4,8]$, $[6,8,9]$, $[6,2,3]$ and $[6,3,7]$ at 6, which is not allowed. If $x_6=9$, then ${\rm lk}(9)=C_7(10,\boldsymbol{1},8, 6,7,4, \boldsymbol{x_7})$, where $x_7 \in \{3,5\}$. But for both the values of $x_7$, considering ${\rm lk}(4)$, we see $\deg(8)>5$. Thus the claim.

\smallskip

\noindent{\bf Subcase 2.13.1.1.}  If $(x_4,x_5)=(2,3)$, then ${\rm lk}(8)=C_7(1,\boldsymbol{0},7, 3,2,9, \boldsymbol{10})$, ${\rm lk}(2)=C_7(1,\boldsymbol{0},3, 8,9,11$, $\boldsymbol{10})$, 
${\rm lk}(3)=C_7(0,\boldsymbol{1},2, 8,7,4, \boldsymbol{5})$,
${\rm lk}(7)=C_7(0,\boldsymbol{1},8, 3,4,6, \boldsymbol{5})$. It is easy to see that 
${\rm lk}(4)=C_7(5,\boldsymbol{0},3,7,6,11, \boldsymbol{10})$, 
${\rm lk}(11)=C_7(10,\boldsymbol{1},2, 9,6,4, \boldsymbol{5})$,
${\rm lk}(9)=C_7(10,\boldsymbol{1},8, 2,11,6, \boldsymbol{5})$, 
${\rm lk}(10) = C_8(1$, $\boldsymbol{2}, 11, \boldsymbol{4}, 5, \boldsymbol{6}, 9, \boldsymbol{8})$,
${\rm lk}(5) = C_8(0, \boldsymbol{3}, 4, \boldsymbol{11}, 10, \boldsymbol{9}, 6, \boldsymbol{7})$, and 
${\rm lk}(6)=C_7(5,\boldsymbol{0},7, 11,4,9, \boldsymbol{10})$. Then $M \cong \mathcal{F}_6(T)$ by the map $(1,7)(2,6)(3,5)$.

\smallskip

\noindent{\bf Subcase 2.13.1.3.} If $(x_4,x_5)=(3,2)$, then ${\rm lk}(8)=C_7(1,\boldsymbol{0},7, 2,3,9, \boldsymbol{10})$, 
${\rm lk}(2)=C_7(1,\boldsymbol{0},3, 8,7,11$, $\boldsymbol{10})$, 
${\rm lk}(7)=C_7(0,\boldsymbol{1},8, 2,11,6, \boldsymbol{5})$
${\rm lk}(3)=C_7(0,\boldsymbol{1},2, 8,9,4, \boldsymbol{5})$
. It is easy to see that 
${\rm lk}(9)=C_7(10,\boldsymbol{1},8, 3,4,6, \boldsymbol{5})$,
${\rm lk}(6)=C_7(5,\boldsymbol{0},7, 11,4,9, \boldsymbol{10})$, 
${\rm lk}(5) = C_8(0, \boldsymbol{3}, 4, \boldsymbol{11}, 10, \boldsymbol{9}, 6, \boldsymbol{7})$,
${\rm lk}(4)=C_7(5$, $\boldsymbol{0},3,9,6,11, \boldsymbol{10})$, 
${\rm lk}(11)=C_7(10,\boldsymbol{1},2, 7,6,4, \boldsymbol{5})$,
${\rm lk}(10) = C_8(1, \boldsymbol{2}, 11, \boldsymbol{4}, 5, \boldsymbol{6}, 9, \boldsymbol{8})$. Then $M \cong \mathcal{F}_1(K)$ by the map $(0,7,8,9,2,11,4,5)(1,10,3,6)$.

\smallskip

\noindent{\bf Subcase 2.13.1.2.} If $(x_4,x_5)=(3,4)$, then ${\rm lk}(8)=C_7(1,\boldsymbol{0},7, 4,3,9, \boldsymbol{10})$. This implies
${\rm lk}(3)=C_7(0,\boldsymbol{1},2, 9,8,4, \boldsymbol{5})$, ${\rm lk}(7)=C_7(0,\boldsymbol{1},8, 4,11,6, \boldsymbol{5})$ and ${\rm lk}(4)=C_7(5,\boldsymbol{0},3,8,7,11, \boldsymbol{x_6})$, where $x_6 =10$. Completing successively, we get ${\rm lk}(2)=C_7(1,\boldsymbol{0},3,9,6,11, \boldsymbol{10})$, ${\rm lk}(9)=C_7(10,\boldsymbol{1},8,3,2,6, \boldsymbol{5})$, ${\rm lk}(6)=C_7(5,\boldsymbol{0},7,11,2,9, \boldsymbol{10})$, ${\rm lk}(10) = C_8(1, \boldsymbol{2}, 11, \boldsymbol{4}, 5, \boldsymbol{6}, 9, \boldsymbol{8})$,
${\rm lk}(5) = C_8(0, \boldsymbol{3}, 4, \boldsymbol{11}, 10, \boldsymbol{9}, 6, \boldsymbol{7})$. Then $M \cong \mathcal{F}_3(T)$ by the map $(0,7,8,10,3,6,1,11,4,5)$.

\smallskip

\noindent{\bf Subcase 2.13.1.4} If $(x_4,x_5)=(11,4)$, then ${\rm lk}(8)=C_7(1,\boldsymbol{0},7, 4,11,9, \boldsymbol{10})$,
  ${\rm lk}(7)=C_7(0,\boldsymbol{1},8, 4,3$, $6, \boldsymbol{5})$, 
   ${\rm lk}(3)=C_7(0,\boldsymbol{1},2, 6,7,4, \boldsymbol{5})$. It is easy to see that 
   ${\rm lk}(4)=C_7(5,\boldsymbol{0},3,7,8,11, \boldsymbol{10})$. Now completing successively, we get
   ${\rm lk}(11)=C_7(10,\boldsymbol{1},2, 9,8,4, \boldsymbol{5})$, ${\rm lk}(2)=C_7(1,\boldsymbol{0},3, 6,9,11, \boldsymbol{10})$, ${\rm lk}(9)=C_7(10,\boldsymbol{1},8, 11,2,6, \boldsymbol{5})$, ${\rm lk}(6)=C_7(5,\boldsymbol{0},7, 3,2,9, \boldsymbol{10})$, ${\rm lk}(5) = C_8(0, \boldsymbol{3}, 4, \boldsymbol{11}, 10, \boldsymbol{9}, 6, \boldsymbol{7})$, ${\rm lk}(10) = C_8(1, \boldsymbol{2}, 11, \boldsymbol{4}, 5, \boldsymbol{6}, 9, \boldsymbol{8})$. Then $M \cong \mathcal{F}_1(K)$ by the map $(0,10)(1,3,9,5,7,11)(4,8)$.

 \smallskip
 
\noindent{\bf Subcase 2.13.2} Let $f_{seq}(8)=(4^4)$. Then ${\rm lk}(8) = C_8(7, \boldsymbol{0}, 1, \boldsymbol{10}, 9, \boldsymbol{x_4}, x_5, \boldsymbol{x_6})$, where we see that $x_5 \in \{2,3,4,11\}$. 

 \noindent{\bf Subcase 2.13.2.1} Let $x_5=2$. Then $(x_4,x_5) \in \{(3,11), (11,3)\}$. If $(x_4,x_6)=(3,11)$, then completing successively, we get 
 ${\rm lk}(8) = C_8(1, \boldsymbol{0}$, $7, \boldsymbol{11}, 2, \boldsymbol{3}, 9, \boldsymbol{10})$, 
  ${\rm lk}(2) = C_8(1, \boldsymbol{0}$, $3, \boldsymbol{9}, 8, \boldsymbol{7}, 11, \boldsymbol{10})$,
  ${\rm lk}(3) = C_8(0, \boldsymbol{1}, 2, \boldsymbol{8}, 9, \boldsymbol{6}, 4, \boldsymbol{5})$, 
 ${\rm lk}(9) = C_8(3, \boldsymbol{2}$, $8, \boldsymbol{1}, 10, \boldsymbol{5}, 6, \boldsymbol{4})$, 
 ${\rm lk}(10) = C_8(1, \boldsymbol{2}$, $11, \boldsymbol{4}, 5, \boldsymbol{6}, 9, \boldsymbol{8})$, 
  ${\rm lk}(5) = C_8(0, \boldsymbol{3}, 4, \boldsymbol{11}, 10, \boldsymbol{9}, 6, \boldsymbol{7})$, 
 ${\rm lk}(6) = C_8(4, \boldsymbol{3}, 9, \boldsymbol{10}, 5, \boldsymbol{0}, 7, \boldsymbol{11})$. This gives a semi equivelar map of type $[4^4]$. 
 
 If $(x_4,x_5)=(11,3)$, then 
 ${\rm lk}(8) = C_8(1, \boldsymbol{0}$, $7, \boldsymbol{3}, 2, \boldsymbol{11}, 9, \boldsymbol{10})$, 
 ${\rm lk}(2) = C_8(1, \boldsymbol{0}$, $3, \boldsymbol{7}, 8, \boldsymbol{9}, 11, \boldsymbol{10})$,
 ${\rm lk}(3) = C_8(0, \boldsymbol{1}, 2, \boldsymbol{8}, 7, \boldsymbol{6}, 4, \boldsymbol{5})$, and 
 ${\rm lk}(7) = C_8(0, \boldsymbol{1}, 8, \boldsymbol{2}, 3, \boldsymbol{4}, 6, \boldsymbol{5})$. This implies $f_{seq}(4)=(3^3.4^2)$ or $(4^4)$.

 \smallskip
 
 \noindent{\bf Subcase 2.13.2.1.1} Let $f_{seq}(4)=(3^3.4^2)$.  Then ${\rm lk}(4)=C_7(3,\boldsymbol{7},6,x_7,x_8,5, \boldsymbol{0})$. It is easy to see that $(x_7,x_8) \in \{(9,10), (9$, $11), (10,9), (10,11), (11,10), (11,9)\}$. Here $(10,11) \cong (9,10)$ by the map $(0,3,7)(1,2,8)(4,6,5)(9,10,11)$, $(11,10) \cong (9,11)$ by the map $(2,8)(3,7)(4,6)(9,11)$, and $(11,9)\cong (9,10)$ by the map $(0,7,3)(1,8,2)(4,5,6)(9,11,10)$. Thus we search for the following cases.
 
 If $(x_7,x_8)=(9, 10)$, then completing successively, we get 
 ${\rm lk}(4)=C_7(3,\boldsymbol{0},5,10,9,6, \boldsymbol{7})$, 
 ${\rm lk}(9)=C_7(8,\boldsymbol{1},10, 4,6,11, \boldsymbol{2})$,
 ${\rm lk}(6)=C_7(7,\boldsymbol{0},5, 11,9,4, \boldsymbol{3})$,
 ${\rm lk}(5)=C_7(0,\boldsymbol{3},4, 10,11,6, \boldsymbol{7})$,
 ${\rm lk}(10)=C_7(1,\boldsymbol{2}$, $11, 5,4,9, \boldsymbol{8})$,
 ${\rm lk}(11)=C_7(2,\boldsymbol{1},10, 5,6,9, \boldsymbol{8})$. Then $M \cong \mathcal{F}_8(K)$ by the map $(0,6,10,8,3,5,9,2)(1$, $7,4,11)$. 
 
 \smallskip
 
 If $(x_7,x_8)=(9, 11)$, then completing successively, we get 
 ${\rm lk}(4)=C_7(3,\boldsymbol{0},5,11,9,6, \boldsymbol{7})$, 
 ${\rm lk}(9)=C_7(8,\boldsymbol{1},10, 6,4,11, \boldsymbol{2})$,
 ${\rm lk}(6)=C_7(7,\boldsymbol{0},5, 10,9,4, \boldsymbol{3})$,
 ${\rm lk}(5)=C_7(0,\boldsymbol{3},4, 11,10,6, \boldsymbol{7})$,
 ${\rm lk}(10)=C_7(1,\boldsymbol{2}$, $11, 5,6,9, \boldsymbol{8})$,
 ${\rm lk}(11)=C_7(2,\boldsymbol{1},10, 5,4,9, \boldsymbol{8})$. Then $M \cong \mathcal{F}_7(T)$ by the map $(0,5,9,2,7,4,11,8,3$, $6,10,1)$. 
 
  \smallskip
 
 If $(x_7,x_8)=(10, 9)$, then completing successively, we get 
 ${\rm lk}(4)=C_7(3,\boldsymbol{0},5,9,10,6, \boldsymbol{7})$, 
 ${\rm lk}(9)=C_7(8,\boldsymbol{1},10, 4,5,11, \boldsymbol{2})$,
 ${\rm lk}(5)=C_7(0,\boldsymbol{3},4, 9,11,6, \boldsymbol{7})$, 
 ${\rm lk}(11)=C_7(2,\boldsymbol{1},10, 6,5,9, \boldsymbol{8})$, 
 ${\rm lk}(6)=C_7(7,\boldsymbol{0}$, $5, 11,10,4, \boldsymbol{3})$,
 ${\rm lk}(10)=C_7(1,\boldsymbol{2}, 11, 6,4,9, \boldsymbol{8})$.
 Then $M \cong \mathcal{F}_9(T)$ by the map $(0,5,11,2,3,4,9,8,7$, $6,10,1)$. 

 \smallskip

 \noindent{\bf Subcase 2.13.2.1.2} Let $f_{seq}(4)=(4^4)$. Then 
 ${\rm lk}(4) = C_8(5, \boldsymbol{0}, 3, \boldsymbol{7}, 6, \boldsymbol{x_7}, x_8, \boldsymbol{x_9})$. Here we see that $(x_7,x_8,x_9) \in \{(10,9,11), (11,9,10), (9, 10,11), (11,10,9), (9,11, 10)\}$. If $(x_7,x_8,x_9)=(10,9,11)$, then completing successively, we get 
 ${\rm lk}(9) = C_8(4, \boldsymbol{5}, 11, \boldsymbol{2}, 8, \boldsymbol{1}, 10, \boldsymbol{6})$, 
 ${\rm lk}(10) = C_8(1, \boldsymbol{2}, 11, \boldsymbol{5}, 6, \boldsymbol{4}, 9$, $\boldsymbol{8})$, 
 ${\rm lk}(6) = C_8(4, \boldsymbol{3}, 7, \boldsymbol{0}, 5, \boldsymbol{11}, 10, \boldsymbol{9})$, 
 ${\rm lk}(5) = C_8(0, \boldsymbol{3}, 4, \boldsymbol{9}, 11, \boldsymbol{10}, 6, \boldsymbol{7})$. This gives a semi-equivelar map of type $[4^4]$. Similarly for the remaining cases of $(x_7,x_8,x_9)$, we get semi equivelar maps of type $[4^4]$. 
 
  \smallskip

\noindent{\bf Subcase 2.13.2.2} Let $x_5=3$. Then $(x_4,x_6) \in \{(2,4), (4,2)\}$. If $(x_4,x_6)=(2,4)$, then completing successively, we get
 ${\rm lk}(8) = C_8(1, \boldsymbol{0}, 7, \boldsymbol{4}, 3, \boldsymbol{2}, 9, \boldsymbol{10})$, 
 ${\rm lk}(3) = C_8(0, \boldsymbol{1}, 2, \boldsymbol{9}, 8, \boldsymbol{7}, 4, \boldsymbol{5})$,
 ${\rm lk}(2) = C_8(1, \boldsymbol{0}, 3, \boldsymbol{8}, 9, \boldsymbol{6}, 11, \boldsymbol{10})$,
 ${\rm lk}(9) = C_8(2, \boldsymbol{3}, 8, \boldsymbol{1}, 10, \boldsymbol{5}, 6, \boldsymbol{11})$, 
 ${\rm lk}(10) = C_8(1, \boldsymbol{2}, 11, \boldsymbol{4}, 5, \boldsymbol{6}, 9, \boldsymbol{8})$,
 ${\rm lk}(5) = C_8(0, \boldsymbol{3}, 4, \boldsymbol{11}, 10, \boldsymbol{9}, 6, \boldsymbol{7})$,
 ${\rm lk}(6) = C_8(5, \boldsymbol{0}, 7, \boldsymbol{4}, 11, \boldsymbol{2}, 9, \boldsymbol{10})$, 
 ${\rm lk}(11) = C_8(2, \boldsymbol{1}, 10, \boldsymbol{5}, 4, \boldsymbol{7}, 6, \boldsymbol{9})$, and 
 ${\rm lk}(4) = C_8(3, \boldsymbol{0}, 5, \boldsymbol{10}, 11, \boldsymbol{6}, 7$, $\boldsymbol{8})$. This gives a semi equivelar map of type $[4^4]$. 
 
  \smallskip
 
  If $(x_4,x_6)=(4,2)$, then successively we get
 ${\rm lk}(8) = C_8(1, \boldsymbol{0}, 7, \boldsymbol{2}, 3, \boldsymbol{4}, 9, \boldsymbol{10})$, 
 ${\rm lk}(3) = C_8(0, \boldsymbol{1}, 2, \boldsymbol{7}$, $8, \boldsymbol{9}, 4, \boldsymbol{5})$,
 ${\rm lk}(2) = C_8(1, \boldsymbol{0}, 3, \boldsymbol{8}, 7, \boldsymbol{6}, 11, \boldsymbol{10})$,
 ${\rm lk}(7) = C_8(0, \boldsymbol{1}, 8, \boldsymbol{3}, 2, \boldsymbol{11}, 6, \boldsymbol{5})$. Then $f_{seq}(4)=(3^3.4^2)$ or $(4^4)$. 
 
  \smallskip

 \noindent{\bf Subcase 2.13.2.2.1} Let $f_{seq}(4)=(3^3.4^2)$. Then  
 ${\rm lk}(4)=C_7(3,\boldsymbol{8},9,x_7,x_8,5, \boldsymbol{0})$, where $(x_7,x_8) \in \{(6,11), (11,10)\}$. In case $(x_7,x_8)=(6,11)$, considering successively ${\rm lk}(4)$, ${\rm lk}(6)$, ${\rm lk}(9)$, and ${\rm lk}(5)$, we see that $\deg(4)>5$. Similarly, if $(x_7,x_8)=(11,10)$, then considering  successively, ${\rm lk}(4)$, ${\rm lk}(10)$, ${\rm lk}(5)$, and ${\rm lk}(9)$, we see that $\deg(4)>5$. 
  
  \smallskip
 
\noindent{\bf Subcase 2.13.2.2.2} Let $f_{seq}(4) = (4^4)$. Then 
 ${\rm lk}(4) = C_8(5, \boldsymbol{0}, 3, \boldsymbol{8}, 9, \boldsymbol{x_7}, x_8$, $\boldsymbol{x_9})$. Observe that $x_8=11$ and $(x_7,x_9) \in \{(6,10), (10,6)\}$. If $(x_7,x_9) = (10,6)$, then 
 ${\rm lk}(11) = C_8(2, \boldsymbol{1}, 10, \boldsymbol{5}, 4, \boldsymbol{9}, 6, \boldsymbol{7})$, which implies $C_5(1,2,11,4,9) \subseteq {\rm lk}(10)$. If $(x_7,x_9)=(6,10)$, then 
 ${\rm lk}(4) = C_8(3, \boldsymbol{0}, 5, \boldsymbol{10}, 11, \boldsymbol{6}, 9$, $\boldsymbol{8})$, 
 ${\rm lk}(11) = C_8(2, \boldsymbol{1}, 10, \boldsymbol{5}, 4, \boldsymbol{9}, 6, \boldsymbol{7})$, 
 ${\rm lk}(6) = C_8(5, \boldsymbol{0}, 7, \boldsymbol{2}, 11, \boldsymbol{4}, 9, \boldsymbol{10})$, and
 ${\rm lk}(9) = C_8(4, \boldsymbol{3}, 8, \boldsymbol{1}, 10, \boldsymbol{5}, 6$, $\boldsymbol{11})$. This gives a semi equivelar map of type $[4^4]$. 
 
  \smallskip

\noindent{\bf Subcase 2.13.2.3} Let $x_5=4$. Then $(x_4,x_6) \in \{(3,11), (5,11), (11,3)\}$. 
 
 If $(x_4,x_6)=(3,11)$, then completing successively, we get
 ${\rm lk}(8) = C_8(1, \boldsymbol{0}, 7, \boldsymbol{11}, 4, \boldsymbol{3}, 9, \boldsymbol{10})$,
 ${\rm lk}(4) = C_8(3, \boldsymbol{0}, 5, \boldsymbol{10}, 11, \boldsymbol{7}, 8$, $\boldsymbol{9})$,
 ${\rm lk}(5) = C_8(0, \boldsymbol{3}, 4, \boldsymbol{11}, 10, \boldsymbol{9}, 6, \boldsymbol{7})$, 
 ${\rm lk}(10) = C_8(1, \boldsymbol{2}, 11, \boldsymbol{4}, 5, \boldsymbol{6}, 9, \boldsymbol{8})$,
 ${\rm lk}(9) = C_8(3, \boldsymbol{2}, 6, \boldsymbol{5}, 10, \boldsymbol{1}, 8, \boldsymbol{4})$, 
 ${\rm lk}(6) = C_8(5, \boldsymbol{0}, 7, \boldsymbol{11}, 2, \boldsymbol{3}, 9, \boldsymbol{10})$.
 This gives a semi equivelar map of type $[4^4]$.
 
  \smallskip
 
 If $(x_4,x_6)=(5,11)$, then completing successively, we get
 ${\rm lk}(8) = C_8(1, \boldsymbol{0}, 7, \boldsymbol{11}, 4, \boldsymbol{5}, 9, \boldsymbol{10})$,
 ${\rm lk}(7) = C_8(0, \boldsymbol{1}, 8, \boldsymbol{4}, 11, \boldsymbol{10}, 6$, $\boldsymbol{5})$,
 ${\rm lk}(5) = C_8(0, \boldsymbol{3}, 4, \boldsymbol{8}, 9, \boldsymbol{10}, 6, \boldsymbol{7})$. 
 This implies $C_5(0,5,9,10,11,7) \subseteq {\rm lk}(6)$.
 
  \smallskip 
 
 If $(x_4,x_6)=(11,3)$, then completing successively, we get
 ${\rm lk}(8) = C_8(1, \boldsymbol{0}, 7, \boldsymbol{3}, 4, \boldsymbol{11}, 9, \boldsymbol{10})$,
 ${\rm lk}(4) = C_8(3, \boldsymbol{0}, 5, \boldsymbol{10}, 11, \boldsymbol{9}, 8$, $\boldsymbol{7})$,
 ${\rm lk}(5) = C_8(0, \boldsymbol{3}, 4, \boldsymbol{11}, 10, \boldsymbol{9}, 6, \boldsymbol{7})$, 
 ${\rm lk}(9) = C_8(11, \boldsymbol{2}, 6, \boldsymbol{5}, 10, \boldsymbol{1}, 8, \boldsymbol{4})$, 
 ${\rm lk}(6) = C_8(5, \boldsymbol{0}, 7, \boldsymbol{3}, 2, \boldsymbol{11}, 9, \boldsymbol{10})$.
 This gives a semi equivelar map of type $[4^4]$.
 
 \smallskip
 
\noindent{\bf Subcase 2.13.2.4}  If $x_5=11$, then $(x_4,x_6) \in \{(2,4), (4,2)\}$. 
 
 If $(x_4,x_6)=(2,4)$, then completing successively, we get
 ${\rm lk}(8) = C_8(1, \boldsymbol{0}, 7, \boldsymbol{4}, 11, \boldsymbol{2}, 9, \boldsymbol{10})$,
 ${\rm lk}(11) = C_8(2, \boldsymbol{1}, 10, \boldsymbol{5}, 4, \boldsymbol{7}, 8, \boldsymbol{9})$,
 ${\rm lk}(10) = C_8(1, \boldsymbol{2}, 11, \boldsymbol{4}, 5, \boldsymbol{6}, 9, \boldsymbol{8})$,
 ${\rm lk}(5) = C_8(0, \boldsymbol{3}, 4, \boldsymbol{11}, 10, \boldsymbol{9}, 6, \boldsymbol{7})$,
 ${\rm lk}(4) = C_8(3, \boldsymbol{0}, 5, \boldsymbol{10}, 11, \boldsymbol{8}, 7$, $\boldsymbol{6})$,
 ${\rm lk}(3) = C_8(0, \boldsymbol{1}, 2, \boldsymbol{9}, 6, \boldsymbol{7}, 4, \boldsymbol{5})$.
 This gives a semi equivelar map of type $[4^4]$.
 
 \smallskip
 
 If $(x_4,x_6)=(4,2)$, then completing successively, we get
 ${\rm lk}(8) = C_8(1, \boldsymbol{0}, 7, \boldsymbol{2}, 11, \boldsymbol{4}, 9, \boldsymbol{10})$,
 ${\rm lk}(7) = C_8(0, \boldsymbol{1}, 8, \boldsymbol{11}, 2, \boldsymbol{3}, 6, \boldsymbol{5})$,
 ${\rm lk}(2) = C_8(1, \boldsymbol{0}, 3, \boldsymbol{6}, 7, \boldsymbol{8}, 11, \boldsymbol{10})$,
 ${\rm lk}(11) = C_8(2, \boldsymbol{1}, 10, \boldsymbol{5}, 4, \boldsymbol{9}, 8, \boldsymbol{7})$,
 ${\rm lk}(4) = C_8(3, \boldsymbol{0}, 5, \boldsymbol{10}, 11, \boldsymbol{8}, 9$, $\boldsymbol{6})$,
 ${\rm lk}(10) = C_8(1, \boldsymbol{2}, 11, \boldsymbol{4}, 5, \boldsymbol{6}, 9, \boldsymbol{8})$.
 This gives a semi equivelar map of type $[4^4]$.
  This proves the lemma. \hfill $\Box$
 
 \smallskip

 \noindent{\bf Proof of Theorem \ref{t2}:} Let $M_1$ and $M_2$ be 2-semi-equivelar maps with the vertex sets $V(M_1)$ and $V(M_2)$ respectively. Clearly, $M_1 \ncong M_2$ if: (i) one is on the torus and other is on the Klein bottle or (ii) their types are distinct or (iii) $|V(M_1)| \neq |V(M_2)|$, where $|V(M_i)|$ denotes the cardinality of $V(M_i)$, for $i \in \{1,2\}$. Now, the proof follows from Theorem \ref{t1}, Lemma \ref{l0} and Lemmas \ref{l2}-\ref{l7}. \hfill $\Box$

	\section{Conclusion} \label{s5}
The 2-semi equivelar maps are generalization of Johnson solids, as are 1-semi equivelar maps of Platonic solids and Archimedean solids. In this article, 2-semi-equivelar maps with curvature 0 have been studied for the surfaces of Euler characteristic 0. It has been obtained that there are exactly 16 types 2-semi-equivelar maps on these surfaces. Further, enumerating the maps for these types on at most 12 vertices, we have obtained 31 2-semi-equivelar maps. Out of which, 18 are on the torus and remaining 13 are on the Klein bottle. The study motivates us to determine  all possible types 2-semi-equivelar maps on the surfaces of Euler characteristic 0, and other closed surfaces. 




	\section{Acknowledgment} The first author expresses his thanks to IIIT Allahabad for providing the facility and resources to carry this research work.

\end{document}